# NONPARAMETRIC HETEROGENEITY TESTING FOR MASSIVE DATA

By Junwei Lu[‡], Guang Cheng[*,§] and Han Liu[†,‡]

*Princeton University*[‡] *and Purdue University*[§]

A massive dataset often consists of a growing number of (potentially) heterogeneous sub-populations. This paper is concerned about testing various forms of heterogeneity in the framework of nonparametric varying-coefficient models. By taking advantage of parallel computing, we design a set of heterogeneity testing procedures with reasonable scalability. In theory, their null limiting distributions are derived as being nearly Chi-square with diverging degrees of freedom as long as the number of sub-populations, denoted as $s$, does not grow too fast. Moreover, we characterize a lower bound on $s$ that is needed to gain sufficient power in testing, so-called "blessing of aggregation" phenomenon. As a by-product, homogeneity testing is also considered with a test statistic being aggregated over all sub-populations. Numerical results are presented to support our theory.

**1. Introduction.** Massive datasets are generally collected from multiple and distributed sources at different time points by different experimental methods. Hence, it is commonly believed that a massive dataset consists of a growing number of (potentially) heterogeneous sub-populations. Heterogeneity testing thus becomes a necessary pre-processing step. In this paper, we aim to develop computationally feasible and theoretically valid heterogeneity testing in presence of a growing number of sub-populations.

Varying coefficient model (Hastie and Tibshirani, 1993) is known to provide

[*]Associate Professor, Department of Statistics, Purdue University, West Lafayette, IN 47906. E-mail: chengg@purdue.edu. Tel: +1 (765) 496-9549. Fax: +1 (765) 494-0558. Research Sponsored by NSF CAREER Award DMS-1151692, DMS-1418042, Simons Fellowship in Mathematics, Office of Naval Research (ONR N00014-15-1-2331) and a grant from Indiana Clinical and Translational Sciences Institute. Guang Cheng was on sabbatical at Princeton while this work was finalized; he would like to thank the Princeton ORFE department for its hospitality and support.

[†]Assistant Professor. Research Sponsored by NSF IIS1408910, NSF IIS1332109, NIH R01MH102339, NIH R01GM083084, and NIH R01HG06841.

*MSC 2010 subject classifications:* Primary 62G10, 62G08; secondary 62E20

a flexible framework for modelling heterogeneous data as follows:

$$Y = \widetilde{\boldsymbol{Z}}^T \boldsymbol{\beta}_0(X) + \epsilon, \tag{1.1}$$

where $X$ is a scalar, $\widetilde{\boldsymbol{Z}}^T = (1, \boldsymbol{Z}^T)$ and $\boldsymbol{\beta}_0(\cdot)$ includes both commonality and heterogeneity components as illustrated by the following example. Cleveland et al. (1993) were interested in exploring how ethanol exhaustion $\mathrm{NO}_x$ for a engine depends on the fuel-air mixture equivalence ratio $E$ and the compression ratio of the engine $C$. They found that $\mathrm{NO}_x$ is linear to $C$ when $E$ is fixed and nonlinear to $E$, which can be modelled as

$$\mathrm{NO}_x = \beta_0(E) + \beta_1(E)C + \epsilon.$$

Nowadays, it becomes common to collect a huge amount of data from various types of engines in distributed places. In this case, one can reasonably assume that the slope $\beta_1(E)$ potentially varies with the engine type, while the intercept $\beta_0(E)$ is a shared factor across all vehicles. In reality, it is of interest in checking whether the effect of compression ratio $C$ is the same for a selection of engine types, and also estimating the common factor $\beta_0$ by utilizing the data collected from all engines.

In view of the above discussions, we assume the following varying coefficient model (as a special case of (1.1)) in the $j$-th sub-population:

$$Y = g_0(X) + Z\beta_{0j}(X) + \epsilon, \text{ for } j = 1, \ldots, s, \tag{1.2}$$

where $s$ denotes the total number of sub-populations. We remark that the homogeneity function $g_0$ and heterogeneity function $\beta_{0j}$ do not necessarily share the same smoothness, and also that a multivariate extension of $Z$ (and thus $\beta_{0j}$) is possible (but omitted for simplicity). The random covariate $Z$ represents the magnitude of heterogeneity, and is assumed to be independent of $X$ and $\epsilon$. Note that allowing $s \to \infty$ makes our work substantially different from existing nonparametric literature, e.g., Pardo-Fernández et al. (2007); Munk and Dette (1998); Srihera and Stute (2010); Chen and Zhong (2010) and Zhu et al. (2012).

A main purpose of this paper is to develop large-scale heterogeneity testing with computational feasibility and theoretical guarantee. As a first step, we estimate both homogeneity function $g_0$ and heterogeneity functions $\beta_{0j}$'s through a two-step procedure in a general kernel ridge regression framework (Shawe-Taylor and Cristianini, 2004). This requires a careful design and study of a general product Hilbert space. In theory, these two estimators are shown to achieve the so-called "oracle efficiency" in the sense that their (joint) asymptotic behaviors are exactly the same as those for oracle estimators (defined in Remark 3.5). This oracle rule holds as long as the number of sub-populations does not grow too fast. The obtained estimators will be used to construct heterogeneity and homogeneity testings, and their oracle efficiency plays a key role in establishing nice theoretical properties of those testing methods.

In this paper, we consider the following hypotheses for detecting various forms of heterogeneity:

(a) Simple test: $H_0 : \beta_{0j} = \beta_0$.
(b) Pairwise test: $H_0 : \beta_{0j} = \beta_{0j'}$, for some $j \neq j'$.
(c) Simultaneous test: $H_0 : \beta_{01} = \cdots = \beta_{0s}$.

A set of test statistics is proposed for (a) – (c) based on the (penalized) likelihood ratio principle (Shang and Cheng, 2013). These testing methods are particularly designed to accommodate a growing number of sub-populations in the sense that they are much more computationally affordable than those commonly used nonparametric multiple-sample tests such as Pardo-Fernández et al. (2007). A relevant heterogeneity testing is considered in Zhao et al. (2016), but only applied to parametric components in partially linear models.

The diverging $s$ brings challenges in analyzing theoretical properties of heterogeneity testing. In particular, we are interested in knowing how fast $s$ is allowed to grow such that these testing methods are still consistent. For this purpose, we provide sufficient conditions on $s$ under which null limiting distributions for (a) – (c) are derived as being nearly Chi-square with diverging degrees of freedom. Furthermore, their power performances are examined under a sequence of local alternatives, measured by separation rates. To obtain *minimal* separation rates, we need to borrow sufficient information across all the sub-populations to reduce the estimation error on $g_0$, which is now treated as a nuisance parameter. Interestingly, this requires $s$ to diverge fast enough, called as "blessing of aggregation" phenomenon.

As a by-product, we also consider homogeneity testing on $g_0$ such as $H_0 : g_0 = g_0^*$ for some specified $g_0^*$. Since $g_0$ is shared among all sub-populations, we certainly want to construct our test statistic based on the entire massive data. This leads to our construction of an *aggregated* test statistics over $s$ sub-populations. Our power analysis further illustrates that this aggregation mechanism yields higher testing sensitivity to a sequence of local alternatives as $s$ grows properly. This again reflects "blessing of aggregation," but from another perspective. Under a different upper bound on $s$, we derive a similar null limiting distribution as in the case of heterogeneity testing. We note that the effects of aggregation have recently been examined from many other viewpoints such as Chen and Xie (2012); Li et al. (2013); Kleiner et al. (2014); Lee et al. (2015); Battey et al. (2015); Wang et al. (2015); Huang and Huo (2015). However, none of these work takes into account of data heterogenetiy, which is a salient feature for massive data.

In the end, we would like to comment on the selection of regularization parameters applied to $g_0$ and $\beta_{0j}$'s. No matter for the purpose of testing or estimation, it is not surprising that the former parameter should be chosen in the scale of entire sample size, while the latter in the scale specific to the sample to which the heterogeneity function belongs. The only exception to this general principle is homogeneity testing in which both regularization parameters need to be in the scale of entire sample size. This is because we need to under-regularize heterogeneity estimators such that their biases play an asymptotically ignorable role in the homogeneity testing statistic.

The rest of this paper is organized as follows. Section 2 reviews basic knowledge on reproducing kernel Hilbert space and further defines a type of product Hilbert space needed in our theoretical analysis. Basic model assumptions and a preliminary two step estimation procedure (together with oracle rule) are introduced in Section 3. Our general theory on heterogeneity and homogeneity testing are presented in Section 4, while more specific examples are discussed in Section 5. Simulations are given in Section 6.

**Notations:** $\mathbb{N} = \{0, 1, 2, \ldots\}$ represents a set of all nonnegative integers. Denote $\delta_{jk}$ as the Kronecker delta, where $\delta_{jk} = 1$ if $j = k$ and $\delta_{jk} = 0$ if $j \neq k$. For positive sequences $a_n$ and $b_n$, we write $a_n \lesssim b_n$ if there exists a constant $c > 0$ independent to $n$ such that $a_n \leq cb_n$ for all $n \in \mathbb{N}$ and $a_n \gtrsim b_n$ if $a_n \leq c'b_n$ for some universal constant $c'$. We denote $a_n \asymp b_n$ if both $a_n \lesssim b_n$ and $a_n \gtrsim b_n$. We also denote $a_n \lesssim b_n$ as $a_n = O(b_n)$. Write $x \vee y = \max(x, y)$, $x \wedge y = \min(x, y)$, and denote $\sup_{x \in \mathcal{X}} |f(x)|$ as $\|f\|_{\sup}$.

**2. Background.** In this section, we provide some necessary background on the reproducing kernel Hilbert space (RKHS).

2.1. *Reproducing Kernel Hilbert Space.* Let $\mathcal{X}$ be a subset of $\mathbb{R}$ and $\mathbb{P}$ be a probability measure defined on the $\mathcal{X}$. Given an RKHS $\mathcal{H} \subseteq L^2(\mathbb{P})$ with an inner product $\langle \cdot, \cdot \rangle_{\mathcal{H}}$ and the induced norm $\| \cdot \|_{\mathcal{H}}$, there exists a symmetric and square integrable function $K : \mathcal{X} \times \mathcal{X} \to \mathbb{R}$ satisfying

$$\langle f, K(x, \cdot) \rangle_{\mathcal{H}} = f(x), \text{ for all } f \in \mathcal{H},$$

where $K_x(\cdot) := K(x, \cdot)$. We call $K$ as the reproducing kernel of $\mathcal{H}$. If $K$ is also a continuous function on $\mathcal{X} \times \mathcal{X}$, Mercer's theorem (Riesz and Sz.-Nagy, 1955) guarantees that the kernel has the following expansion

$$K(x, t) = \sum_{k=0}^{\infty} \mu_k \phi_k(x) \phi_k(t),$$

where $\{\mu_k\}_{k=0}^{\infty}$ is a non-negative descending sequence of eigenvalues and $\{\phi_k\}_{k=0}^{\infty}$ are eigen-functions forming an orthonormal basis for $L^2(\mathbb{P})$. In particular, the eigen-system satisfies for any $j, k \in \mathbb{N}$,

$$\langle \phi_j, \phi_k \rangle_{L^2(\mathbb{P})} = \delta_{jk} \quad \text{and} \quad \langle \phi_j, \phi_k \rangle_{\mathcal{H}} = \delta_{jk}/\mu_k.$$

The underlying smoothness of $f \in \mathcal{H}$ can be characterized by the decaying rate of the eigenvalues $\{\mu_k\}_{k=0}^{\infty}$. We focus on three typical decaying rates in this paper. As will be demonstrated later, the nonparametric inference for kernel ridge regression crucially depends on the decaying rate.

**Finite rank:** There exists some $r \in \mathbb{N}$ such that $\mu_k = 0$ for $k > r$. This gives arise to the so-called finite rank kernel (with rank $r$). For example, the linear kernel $K(\boldsymbol{z}_1, \boldsymbol{z}_2) = \langle \boldsymbol{z}_1, \boldsymbol{z}_2 \rangle_{\mathbb{R}^d}$ has rank $d$, and generates a $d$-dimensional linear function space. The eigen-functions are given by $\phi_\ell(\boldsymbol{z}) = z_\ell$ for $\ell = 1, \ldots, d$. For one dimensional $z$, the polynomial kernel $K(z_1, z_2) = (1+z_1 z_2)^m$ has rank $m+1$, and induces a space of polynomial functions with degree at most $m$. The eigen-functions are given by $\phi_\ell = z^{\ell-1}$ for $\ell = 1, \ldots, m+1$.

**Exponentially decaying:** There exist some $\alpha, p > 0$ such that $\mu_k \asymp \exp(-\alpha k^p)$. Exponentially decaying kernels include the Gaussian kernel $K(x, t) = \exp(-\|x - t\|_2^2/\sigma^2)$. The corresponding eigen-functions are

$$(2.1) \qquad \phi_k(x) = (\sqrt{5}/4)^{1/4} (2^{k-1} k!)^{-1/2} e^{-(\sqrt{5}-1)x^2/4} H_k\big((\sqrt{5}/2)^{1/2} x\big),$$

and the eigenvalues are $\mu_k = \ell^{2k+1}$, for $k \in \mathbb{N}$, where $\ell = (\sqrt{5} - 1)/2$ and $H_k(\cdot)$ is the $k$-th Hermite polynomial (Sollich and Williams, 2005).

**Polynomially decaying:** There exists some $m > 0$ such that $\mu_k \asymp k^{-2m}$. The polynomially decaying function class includes Sobolev space and Besov Space (Birman and Solomjak, 1967). We denote an $m$-th order Sobolev space as $\mathcal{H}^m([0,1])$[1]. For its homogeneous subspace $\mathcal{H}_0^m([0,1])$[2], we have trigonometric eigen-functions as follows:

$$(2.2) \qquad \phi_k(x) = \begin{cases} 1, & k = 0 \\ \sqrt{2} \cos(2\pi \ell x), & k = 2\ell \text{ for } \ell = 1, 2, \ldots, \\ \sqrt{2} \sin(2\pi \ell x), & k = 2\ell - 1 \text{ for } \ell = 1, 2, \ldots, \end{cases}$$

and eigenvalues are $\mu_{2\ell-1} = \mu_{2\ell} = (2\pi\ell)^{-2m}$ for $\ell \geq 0$ and $\mu_0 = \infty$. The corresponding Sobolev kernels can be found in Gu (2013).

2.2. *Product Hilbert Space.* In this section, we define a type of product Hilbert space that is needed in analyzing our main model (1.2).

Suppose the true functions $g_0 \in \mathcal{H}_1$ and $\beta_{0j} \in \mathcal{H}_2$ for all $1 \leq j \leq s$, where

---

[1] $\mathcal{H}^m([0,1]) = \{g : [0,1] \mapsto \mathbb{R} | \ g^{(j)}$ is abs. cont. for $j = 0, \ldots, m-1$, and $g^{(m)} \in L_2([0,1])\}$.
[2] $\mathcal{H}_0^m([0,1]) \subseteq \mathcal{H}^m([0,1])$ is a class of periodic functions satisfying a set of additional boundary conditions: $g^{(j)}(0) = g^{(j)}(1)$ for $j = 0, \ldots, m-1$.

$\mathcal{H}_1$ and $\mathcal{H}_2$ are two (possibly different) RKHS's. Define a product space as

$$\mathcal{H}^2 = \Big\{ f \mid f(u) = g(x) + z\beta(x), \text{for } u = (x,z) \in \mathcal{X} \times \mathcal{Z} := \mathcal{U}, g \in \mathcal{H}_1, \beta \in \mathcal{H}_2 \Big\}.$$

Denote $f = (g, \beta)$ if and only if $f(u) = g(x) + z\beta(x)$ for any $u = (x,z)$. For any $f_1 = (g_1, \beta_1), f_2 = (g_2, \beta_2) \in \mathcal{H}^2$, we define the inner product on $\mathcal{H}^2$:

$$(2.3) \qquad \langle f_1, f_2 \rangle = \mathbb{E}\big[f_1(U)f_2(U)\big] + \lambda_1 \langle g_1, g_2 \rangle_{\mathcal{H}_1} + \lambda_2 \langle \beta_1, \beta_2 \rangle_{\mathcal{H}_2},$$

where the random variable $U := (X, Z)$ takes values on $\mathcal{U} = \mathcal{X} \times \mathcal{Z}$. We denote the induced norm $\|f\|^2 = \langle f, f \rangle$. Moreover, for any $g_1, g_2 \in \mathcal{H}_1, \beta_1, \beta_2 \in \mathcal{H}_2$, we simply denote $\langle g_1, g_2 \rangle = \langle (g_1, 0), (g_2, 0) \rangle$, $\langle \beta_1, \beta_2 \rangle = \langle (0, \beta_1), (0, \beta_2) \rangle$ and similar for the norm $\|\cdot\|$. The inner product is designed as (2.3) for facilitating the asymptotic analysis of our estimators defined later. For example, the second Fréchet derivative of the objective function for estimation turns out to be identity under this inner product; see Appendix A.2 for more details. The following result shows that $(\mathcal{H}^2, \langle \cdot, \cdot \rangle)$ is indeed a RKHS.

**Proposition 2.1.** There exists a kernel function $K^\lambda : \mathcal{U} \times \mathcal{U} \to \mathbb{R}$ and an operator $W_\lambda : \mathcal{H}^2 \mapsto \mathcal{H}^2$ such that for any $f_1 = (g_1, \beta_1), f_2 = (g_2, \beta_2) \in \mathcal{H}^2$ and $u = (x, z) \in \mathcal{U}$, we have

$$\langle K_u^\lambda, f_1 \rangle = g_1(x) + z\beta_1(x) \text{ and } \langle W_\lambda f_1, f_2 \rangle = \lambda_1 \langle g_1, g_2 \rangle_{\mathcal{H}_1} + \lambda_2 \langle \beta_1, \beta_2 \rangle_{\mathcal{H}_2},$$

where $K_u^\lambda(\cdot) = K^\lambda(u, \cdot)$.

**3. Two-Step Estimation.** In this section, we estimate the nonparametric heterogenous model (1.2) in two steps. The obtained estimators for homogeneity and heterogeneity functions will be used in the main result of this paper: heterogeneity/homogeneity testing (Section 4). In theory, these two estimators are shown to achieve the so-called "oracle efficiency" in the sense that their (joint) asymptotic behaviors are exactly the same as those for oracle estimators; see Remark 3.5. This oracle rule holds as long as the number of sub-populations does not grow too fast and the smoothing parameters in the kernel ridge estimation are properly chosen.

In the $j$-th sub-population, we observe $\{(Y_i, X_i, Z_i)\}_{i \in I_j}$ (with $|I_j| = n$)

$$Y_i = g_0(X_i) + Z_i \beta_{0j}(X_i) + \epsilon_i,$$

where $g_0 \in \mathcal{H}_1$ and $\beta_{0j} \in \mathcal{H}_2$. Denote $n \times s$ as the entire sample size $N$. With a bit abuse of notation, we write $f_{0j} = (g_0, \beta_{0j})$. For simplicity, we assume $(X_i, Z_i, \epsilon_i)$'s are independent and follow the same distribution across all sub-populations, denoted as $(X, Z, \epsilon)$.

3.1. *Model Assumptions.* In this section, we introduce three regularity assumptions: Assumption 3.1 on the model specifications; Assumption 3.2 on the reproducing kernel $K$; Assumption 3.3 on the complexity of $\mathcal{H}_1, \mathcal{H}_2$.

**Assumption 3.1** (Covariate and Noise). The random variables $(X, Z, \epsilon)$ are pairwise independent with $\mathbb{E}(Z) = \mathbb{E}(\epsilon) = 0$. We assume that $|Z| \leq c_z$ almost surely and $\epsilon$ is a subgaussian random variable such that there exists a $\sigma_\epsilon > 0$ and $\mathbb{E}[\exp(t\epsilon)] \leq \exp(t^2 \sigma_\epsilon^2/2)$ for any $t \in \mathbb{R}$. Denote $\mathrm{Var}(Z) = \sigma_z^2$ and $\mathrm{Var}(\epsilon) = \sigma^2$.

**Assumption 3.2** (Kernel and Eigenfunction). The kernel functions of $\mathcal{H}_1$, $\mathcal{H}_2$, denoted by $K_1(x,t)$, $K_2(x,t)$, are bounded on the diagonal by some constant $c_b$ such that
$$\sup_{x \in \mathcal{X}} |K_1(x,x) \vee K_2(x,x)| \leq c_b.$$
Their eigen-functions $\{\phi_k\}_{k=0}^\infty$ and $\{\psi_k\}_{k=0}^\infty$ are uniformly bounded on $\mathcal{X}$, namely there exists a constant $c_\phi$ such that
$$(3.1) \qquad \sup_{j \in \mathbb{N}} \{\|\phi_j\|_{\sup} \vee \|\psi_j\|_{\sup}\} \leq c_\phi.$$

We define two useful constants, i.e., $c_\mathcal{H} = c_b(1 + c_z)$ and $c_K = c_\phi \sqrt{1 + c_z}$, which will be used in stating our inference results later.

Assumption 3.2 is commonly used in the literature and is satisfied for a variety of kernels (Shang and Cheng, 2013; Zhang et al., 2013; Lafferty and Lebanon, 2005; Guo, 2002). For example, the polynomial kernel $K(x,t) = (1 + xt)^d$ is uniformly bounded as long as $\mathcal{X}$ is a compact subset of $\mathbb{R}$. The trigonometric basis (2.2) in the periodic Sobolev space $H_0^m([0,1])$ trivially satisfies this assumption. In fact, Proposition 2.2 in Shang and Cheng (2013) proved that the eigen-functions of the more general $H^m([0,1])$ are bounded under a mild condition on covariate density. More generally, Lemma 3 in Minh et al. (2006) says that the eigen-functions of any kernel function expressed in the form $K(x,t) = g(\langle x,t \rangle)$, where $g : [-1,1] \to \mathbb{R}$ is any continuous function, satisfy $\sup_{k \in \mathbb{N}} \|\phi_k\|_{\sup} < \infty$ if $\mathcal{X} = S^1 := \{x \in \mathbb{R}^2 \,|\, \|x\|_2 = 1\}$.

Before stating Assumption 3.3, we need to introduce a few notations related to the covering number. Suppose $\mathcal{F}$ is a function space. An $\epsilon$-net for $(\mathcal{F}, d)$ with $d$ being a metric is a set $\{f_1, \ldots, f_M\} \subseteq \mathcal{F}$ such that for any $f \in \mathcal{F}$, there exists an $f_j$ for some $j = 1, \ldots, M$ satisfying $d(f, f_j) \leq \epsilon$. The $\epsilon$-covering number is defined as
$$\mathcal{N}(\mathcal{F}; d, \epsilon) = \min\{|\mathcal{N}| : \mathcal{N} \text{ is an } \epsilon\text{-net for } (\mathcal{F}, d)\}.$$
A subspace in $\mathcal{H}^2$ of particular interest is defined as
$$(3.2) \qquad \mathcal{H}_0 = \{f = (g, \beta) \in \mathcal{H}^2 : \|f\|_{\sup} \leq 1, \|g\|_{\mathcal{H}_1} \leq 1, \|\beta\|_{\mathcal{H}_2} \leq 1\},$$

whose entropy integral is defined as

$$\omega(\delta) = \int_0^\delta \sqrt{\log(1 + \mathcal{N}(\mathcal{H}_0; \|\cdot\|_{\sup}, \epsilon))} d\epsilon. \tag{3.3}$$

**Assumption 3.3.** (Complexity of $\mathcal{H}_0$) For any fixed $a > 0$, $\omega(a\delta)/\delta$ is always upper bounded by some non-increasing function, denoted as $\bar{\omega}_a(\delta)$, for $\delta \in (0,1)$.

Assumption 3.3 characterizes the complexity of $\mathcal{H}_0$ through its entropy integral $\omega(\delta)$, and is satisfied for function classes considered in this paper. Suppose $\mathcal{H}_1 = \mathcal{H}_2 = \mathcal{H}$. When the eigenvalues of $\mathcal{H}$ decay polynomially, the entropy integral can be written as $\omega(\delta) = c\delta^q$, where $q$ is some constant on $(0,1)$. When $\mathcal{H}$ is finite rank or exponentially decaying, the entropy integral has the form $\omega(\delta) = c\delta(\log(1/\delta))^q$ for some $0 < q < 1$. Obviously, Assumption 3.3 holds in both cases. When two RKHS's are different, $\omega(\cdot)$ can be bounded by a summation of the above two forms, which still satisfies Assumption 3.3.

3.2. *Oracle Rule.* The two-step estimation procedure is described in this subsection. A nice "oracle rule" holds for both homogeneity and heterogeneity estimators once the number of sub-populations does not grow too fast. Our general theory is further applied to three specific RKHS families.

**Step 1.** (Individual Estimation) In the $j$-th sub-population, we estimate $(g_0, \beta_{0j})$ as

$$\left(\widehat{g}_{n,\lambda}^{(j)}, \widehat{\beta}_{n,\lambda}^{(j)}\right) := \widehat{f}_{n,\lambda}^{(j)} = \operatorname*{argmin}_{f \in \mathcal{H}^2} \mathcal{L}_{n,\lambda}^{(j)}(f), \tag{3.4}$$

where the loss function is

$$\mathcal{L}_{n,\lambda}^{(j)}(f) := \frac{1}{2n} \sum_{i \in I_j} (Y_i - g(X_i) - Z_i\beta(X_i))^2 + \frac{\lambda_1}{2}\|g\|_{\mathcal{H}_1}^2 + \frac{\lambda_2}{2}\|\beta\|_{\mathcal{H}_2}^2. \tag{3.5}$$

**Step 2.** (Aggregation for Commonality) Average over all the estimators for homogeneity function:

$$\bar{g}_{N,\lambda} = \frac{1}{s}\sum_{j=1}^s \widehat{g}_{n,\lambda}^{(j)}. \tag{3.6}$$

Several quantities need to be defined before stating our oracle result. Denote $\{\mu_{1,j}\}_{j=0}^\infty$ and $\{\mu_{2,j}\}_{j=0}^\infty$ as the eigenvalues of $\mathcal{H}_1$ and $\mathcal{H}_2$, and

$$\gamma_1 = \gamma_1(\lambda_1) := \sum_{j=0}^\infty (1 + \lambda_1/\mu_{1,j})^{-1}, \quad \gamma_2 = \gamma_2(\lambda_2) := \sum_{j=0}^\infty (\sigma_z^2 + \lambda_2/\mu_{2,j})^{-1},$$

$$\gamma = \gamma_1(\lambda_1) + \gamma_2(\lambda_2), \; h_1 = 1/\gamma_1, h_2 = 1/\gamma_2, h = 1/\gamma.$$

Also, we define a quantity related to the statistical rate

$$(3.7) \qquad r_n := (\log N)/\sqrt{nh} + \sqrt{\lambda_1} + \sqrt{\lambda_2}.$$

and the measure on the capacity of RKHS:

$$(3.8) \qquad \xi_\lambda = (c_K^2 \gamma \lambda)^{-1/2} \omega\big((c_K^2 \gamma \lambda)^{1/2}\big),$$

where $\lambda = \lambda_1 \wedge \lambda_2$ and $c_K$ is specified in Assumption 3.2. The conditions on $(s, \lambda_1, \lambda_2)$ required in Theorem 3.4 are determined by the above quantities.

The following theorem depicts an asymptotic normality of $(\bar{g}_{N,\lambda}, \widehat{\beta}_{n,\lambda}^{(j)})$, irrespective of their different convergence rates.

**Theorem 3.4.** Suppose Assumptions 3.1 – 3.3 hold. Given any fixed $x_1, x_2 \in \mathcal{X}$, we assume $\lambda = o(1)$ and

$$(3.9) \qquad \begin{aligned} \lim_{n \to \infty} \sigma^2 h_1 \sum_{k=0}^\infty \left(\frac{\phi_k(x_1)}{1+\lambda_1/\mu_{1,k}}\right)^2 = \sigma_g^2 < \infty, \\ \lim_{n \to \infty} \sigma^2 \sigma_z^2 h_2 \sum_{k=0}^\infty \left(\frac{\psi_k(x_2)}{\sigma_z^2+\lambda_2/\mu_{2,k}}\right)^2 = \sigma_\beta^2 < \infty. \end{aligned}$$

If $s, \lambda_1, \lambda_2$ satisfy

$$(3.10) \qquad \log\big(N\xi_\lambda/(h\lambda)\big) = o(\log^2 N),$$

$$(3.11) \qquad (nh)^{-1}(\sqrt{n} \cdot \xi_\lambda + 1) r_n \log N = o(N^{-1/2}),$$

$$(3.12) \qquad r_n^{-1} \cdot \omega\big(r_n\sqrt{\lambda/h}\big)(nh\lambda)^{-1/2} \log N = o(1),$$

we have the following joint asymptotic normality result

$$(3.13) \qquad \begin{pmatrix} \sqrt{Nh_1}\,\big(\bar{g}_{N,\lambda}(x_1) - g_0^u(x_1)\big) \\ \sqrt{nh_2}\,\big(\widehat{\beta}_{n,\lambda}^{(j)}(x_2) - \beta_{0j}^u(x_2)\big) \end{pmatrix} \rightsquigarrow N\left(\mathbf{0}, \begin{pmatrix} \sigma_g^2 & 0 \\ 0 & \sigma_\beta^2 \end{pmatrix}\right),$$

where the re-centered functions $(g_0^u, \beta_{0j}^u) = f_{0j} + W_\lambda f_{0j} := f_{0j}^u$, for any $j = 1, \ldots, s$. Recall that the bias operator $W_\lambda$ is defined in Proposition 2.1.

Note that when $\mathcal{H}_1 = H_0^m([0,1])$, $\sigma_g^2$ has an explicit expression given in Shang and Cheng (2013). The asymptotic independence between $\bar{g}_{N,\lambda}$ and $\widehat{\beta}_{n,\lambda}^{(j)}$ is not surprising since the sub-sample used to estimate the heterogeneity is asymptotically ignorable comparing with the entire sample size used to estimate the homogeneity, i.e., $n/N = s^{-1} \to 0$. We remark that (3.10) controls some tail probability such that the union bound argument can be applied to a diverging number of sub-populations, while (3.11) and (3.12) control the nonparametric estimation error in each sub-population.

**Remark 3.5.** An important consequence of Theorem 3.4 is the following oracle rule. Define the oracle estimate for $g_0$ and $\beta_{0j}$ as follows:

$$\widehat{g}_{\text{or}} := \operatorname{argmin}_{g \in \mathcal{H}_1} N^{-1} \sum_{j=1}^{s} \sum_{i \in I_j} \big(Y_i - g(X_i) - Z_i \beta_{0j}(X_i)\big)^2 + \lambda_1 \|g\|_{\mathcal{H}_1}^2,$$
$$\widehat{\beta}_{\text{or},j} := \operatorname{argmin}_{\beta \in \mathcal{H}_2} n^{-1} \sum_{i \in I_j} \big(Y_i - g_0(X_i) - Z_i \beta(X_i)\big)^2 + \lambda_2 \|\beta\|_{\mathcal{H}_2}^2.$$

These two estimators are "oracle" in the sense that we use true homogeneity/heterogeneity function in the estimation. The proof of Theorem 3.4 implies that $(\widehat{g}_{\text{or}}, \widehat{\beta}_{\text{or},j})$ shares the same joint distribution as $(\bar{g}_{N,\lambda}, \widehat{\beta}_{n,\lambda}^{(j)})$ (after similar re-centering and re-scaling). This type of "oracle efficiency" indicates that our two-step estimation efficiency cannot be further improved.

Corollary B.3 in Appendix B gives more explicit conditions in terms of $(s, \lambda_1, \lambda_2)$ for three RKHS families. From Corollary B.3, it is easy to see that $\lambda_1$ should be chosen in the scale of $N$, while $\lambda_2$ in the scale of $n$. In other words, we need to sub-regularize the homogeneity estimation in each sub-population. Without heterogeneity part, Zhang et al. (2013) discovered a similar phenomenon for estimating $g_0$ optimally. We show that this phenomenon persists at an inferential level (in a more general context).

**4. Penalized Likelihood Ratio Theory.** In this section, we apply the penalized likelihood ratio principle advocated in (Shang and Cheng, 2013) to various types of heterogeneity/homogeneity testing. For simplicity, we only focus on the case that $\mathcal{H}_1 = \mathcal{H}_2 = \mathcal{H}$ in this section.

4.1. *Heterogeneity Testing.* In this section, we are interested in testing the following three hypotheses:

(a) Simple test: $H_0 : \beta_{0j} = \beta_0$.
(b) Pairwise test: $H_0 : \beta_{0j} = \beta_{0j'}$, for some $j \neq j'$.
(c) Simultaneous test: $H_0 : \beta_{01} = \cdots = \beta_{0s}$.

In all these tests, the homogeneity function is treated as nuisance. The simultaneous test (c) is concerned about whether the entire dataset is homogenous or not. We are particularly interested in the largest number of sub-populations that can be tested in (c) with theoretical guarantee, e.g., minimax optimality. This is technically challenging given that $s$ diverges.

Let us begin with the simple test (a). A penalized likelihood ratio test (PLRT) statistic is constructed as follows:

$$(4.1) \qquad \text{PLRT}_{n,\lambda} = \mathcal{L}_{n,\lambda}^{(j)}\big(\bar{g}_{N,\lambda}, \widehat{\beta}_{n,\lambda}^{(j)}\big) - \mathcal{L}_{n,\lambda}^{(j)}\big(\bar{g}_{N,\lambda}, \beta_0\big)$$

by plugging in the aggregated estimate $\bar{g}_{N,\lambda}$. Recall that $\mathcal{L}_{n,\lambda}^{(j)}$ is defined in

(3.5). As for the pairwise test (b), we propose

$$\text{PLRT}_{n,\lambda}^{(j,j')} = \left[\mathcal{L}_{n,\lambda}^{(j)}(\bar{g}_{N,\lambda}, \widehat{\beta}_{n,\lambda}^{(j)}) - \mathcal{L}_{n,\lambda}^{(j)}(\bar{g}_{N,\lambda}, \widehat{\beta}_{n,\lambda}^{(j')})\right]$$
$$+ \left[\mathcal{L}_{n,\lambda}^{(j')}(\bar{g}_{N,\lambda}, \widehat{\beta}_{n,\lambda}^{(j')}) - \mathcal{L}_{n,\lambda}^{(j')}(\bar{g}_{N,\lambda}, \widehat{\beta}_{n,\lambda}^{(j)})\right]$$
$$:= T_1 + T_2.$$

In fact, either $T_1$ or $T_2$ is a valid test statistic. Here, we take their summation in order to increase its power in practice. To better understand $T_1$, we decompose it as $T_1 = T_{11} - T_{12}$ under the null $\beta_{0j} = \beta_{0j'}$, where

$$T_{11} = \mathcal{L}_{n,\lambda}^{(j)}(\bar{g}_{N,\lambda}, \widehat{\beta}_{n,\lambda}^{(j)}) - \mathcal{L}_{n,\lambda}^{(j)}(\bar{g}_{N,\lambda}, \beta_{0j}),$$
$$T_{12} = \mathcal{L}_{n,\lambda}^{(j)}(\bar{g}_{N,\lambda}, \widehat{\beta}_{n,\lambda}^{(j')}) - \mathcal{L}_{n,\lambda}^{(j)}(\bar{g}_{N,\lambda}, \beta_{0j'}).$$

Note that $T_{11}$ and $T_{12}$ are both PLRT statistic for some simple test (4.1). Hence, the null limiting distribution of $\text{PLRT}_{n,\lambda}^{(j,j')}$ follows from that of $\text{PLRT}_{n,\lambda}$.

For the most challenging test (c), we can similarly construct (also see Remark 4.4)

(4.2) $\quad \text{PLRT}_{n,\lambda}^{s} = \sum_{j=1}^{s-1} \mathcal{L}_{n,\lambda}^{(j)}(\bar{g}_{N,\lambda}, \widehat{\beta}_{n,\lambda}^{(j)}) - \mathcal{L}_{n,\lambda}^{(j)}(\bar{g}_{N,\lambda}, \widehat{\beta}_{n,\lambda}^{(j+1)}).$

Alternatively, a maximum type test also works here

$$\max_{1 \leq j \neq j' \leq s} \left|\mathcal{L}_{n,\lambda}^{(j)}((\bar{g}_{N,\lambda}, \widehat{\beta}_{n,\lambda}^{(j)})) - \mathcal{L}_{n,\lambda}^{(j)}((\bar{g}_{N,\lambda}, \widehat{\beta}_{n,\lambda}^{(j')}))\right|,$$

but with slightly more computation as $s$ grows. We want to point out that $\text{PLRT}_{n,\lambda}^{s}$ is more computationally feasible than the commonly used nonparametric multiple-sample test (Pardo-Fernández et al., 2007), which requires to calculate an additional estimator based on the entire sample. This additional computation may not be affordable in the setting of massive data. In addition, the nonparametric multiple-sample tests in the literature (without nuisance $g_0$) such as Pardo-Fernández et al. (2007), Munk and Dette (1998) and Srihera and Stute (2010) all require the number of tested functions to be fixed. However, our test method remains valid even when $s$ diverges to infinity at a proper rate; see Theorems 4.1 and 4.2 below.

We are ready to present null limiting distributions for tests (a) – (c). In particular, a statistic $T_n$ is said to be nearly $\chi_{b_n}^2$, denoted as $T_n \overset{a}{\sim} \chi_{b_n}^2$, if $(2b_n)^{-1/2}(T_n - b_n)$ weakly converges to $N(0,1)$ for some sequence $b_n \to \infty$.

**Theorem 4.1.** Suppose that Assumptions 3.1 – 3.3, and (3.10) – (3.12) hold, and $\mathbb{E}[\epsilon^4|X] \leq C$, a.s., for some constant $C$. If $nh_2\lambda_2 = o(1)$ and

$(nh_2^2)^{-1} = o(1)$, we have

$$-2nr_T \cdot \text{PLRT}_{n,\lambda} \overset{a}{\sim} \chi^2_{u_n} \quad \text{and} \quad -\frac{1}{2}nr_T \cdot \text{PLRT}_{n,\lambda}^{(j,j')} \overset{a}{\sim} \chi^2_{u_n},$$

where $u_n = h_2^{-1}\sigma_T^4/\rho_T^2$, $r_T = \sigma_T^2/\rho_T^2$, and

$$\sigma_T^2 = \sigma^2 \sigma_z^2 h_2 \sum_{k=0}^{\infty} (\sigma_z^2 + \lambda_2/\mu_{2,k})^{-1}, \rho_T^2 = \sigma^4 \sigma_z^4 h \sum_{k=0}^{\infty} (\sigma_z^2 + \lambda_2/\mu_{2,k})^{-2}.$$

If we further assume $(Nh_2^2)^{-1} = o(1)$, then

(4.3) $$-\frac{4}{3}nr_T \cdot \text{PLRT}_{n,\lambda}^s \overset{a}{\sim} \chi^2_{\widetilde{u}_N},$$

where $\widetilde{u}_N = (4s/3)h_2^{-1}\sigma_T^4/\rho_T^2$.

We next study the power performance of the above nonparametric testing in terms of their separation rates. Define a function space

(4.4) $$\mathcal{F}_\zeta \equiv \{\beta \in \mathcal{H} | \text{Var}(\beta(X)^2) \leq \zeta \mathbb{E}^2[\beta(X)^2] \text{ and } \|\beta\|_{\mathcal{H}}^2 \leq \zeta\},$$

for some constant $\zeta > 0$, and two separation rates as

(4.5) $$\eta_n^2 = \lambda_2 + \sqrt{1/(n^2 h_2)} + r_N^2,$$

(4.6) $$\widetilde{\eta}_N^2 = \lambda_2 + \sqrt{s/(n^2 h_2)} + r_N^2,$$

where $r_N = \log N/\sqrt{Nh_1} + \sqrt{\lambda_1} + \sqrt{\lambda_2}$.

**Theorem 4.2.** Suppose that the same conditions as Theorem 4.1 hold.

- Under the alternative sequence $H_{1n} : \beta_{0j} = \beta_n$, we have for any $c_n \to \infty$,

$$\inf_{\substack{g_0 \in \mathcal{F}_\zeta \\ \delta_{0j}^* \geq c_n \eta_n^2}} \inf_{\beta_0 - \beta_n \in \mathcal{F}_\zeta} \mathbb{P}_{H_{1n}}\left((2u_n)^{-1/2}\left|-2nr_T \cdot \text{PLRT}_{n,\lambda} - u_n\right| > z_\alpha\right) \to 1,$$

where $\delta_{0j}^* := \|\beta_0 - \beta_n\|^2$.

- Under the alternative sequence $H_{1n} : \beta_{0j} \neq \beta_{0j'}$, we have for any $c_n \to \infty$,

$$\inf_{\substack{g_0 \in \mathcal{F}_\zeta \\ \Delta_2^* \geq c_n \eta_n^2}} \inf_{\beta_{0j} - \beta_{0j'} \in \mathcal{F}_\zeta} \mathbb{P}_{H_{1n}}\left((2u_n)^{-1/2}\left|-\frac{1}{2}nr_T \cdot \text{PLRT}_{n,\lambda}^{(j,j')} - u_n\right| > z_\alpha\right) \to 1,$$

where $\Delta_2^* := \|\beta_{0j} - \beta_{0j'}\|^2$.

- Under the alternative sequence $H_{1n} : \beta_{0j} \neq \beta_{0j'}$ for at least one pair

$(j, j')$, we have for any $c_n \to \infty$

$$\inf_{\substack{\{(\beta_{0j}-\beta_{0(j+1)})\}_{j=1}^{s-1} \subseteq \mathcal{F}_\zeta \\ g_0 \in \mathcal{F}_\zeta, \Delta_s^* \geq c_n \widetilde{\eta}_N^2}} \mathbb{P}_{H_{1n}}\left((2\widetilde{u}_N)^{-1/2}\left|-\frac{4}{3}nr_T \cdot \text{PLRT}_{n,\lambda}^s - \widetilde{u}_N\right| > z_\alpha\right) \to 1,$$

where $\Delta_s^* := \sum_{j=1}^{s-1} \|\beta_{0j} - \beta_{0(j+1)}\|^2$ and "$\{(\beta_{0j} - \beta_{0(j+1)})\}_{j=1}^{s-1} \subseteq \mathcal{F}_\zeta$" means that each $(\beta_{0j} - \beta_{0(j+1)}) \in \mathcal{F}_\zeta$ for $j = 1, \ldots, s-1$.

When $g_0$ is known, we can show that Theorem 4.2 still holds for hypotheses (a) and (b) by simply replacing (4.5) by $\eta_n^2 \gtrsim \lambda_2 + \sqrt{1/(n^2 h_2)}$. In other words, $r_N^2$ in (4.5) is due to the estimation of the nuisance parameter $g_0$ based on the entire sample. However, this term can be dominated by the first two terms in the RHS of (4.5) if we can borrow sufficient information across all the sub-populations to reduce the estimation error, i.e., $s$ diverges fast enough. This explains why we need to set a lower bound on $s$ in order to obtain the minimal possible separation rate. We call this lower bound phenomenon as "blessing of aggregation". See below for more details.

Following Theorem 4.2, we derive the minimal possible values of $\eta_n$ for tests (a) and (b), denoted as $\eta_n^*$, in three RKHS families.

- **Finite rank kernel (with rank $r$)**

$$\eta_n^* = O(r^{1/4}/\sqrt{n})$$

given that $\lambda_1 \asymp r/N$, $\lambda_2 \asymp \sqrt{r}/n$ and $\sqrt{r} \lesssim s \lesssim r^{-5/2} N^{1/2} (\log N)^{-4}$;

- **Exponentially decaying kernel (with $\mu_k \asymp e^{-\alpha k^p}$)**

$$\eta_n^* = O((\log n)^{1/(4p)}/\sqrt{n})$$

given that $\lambda_1 \asymp (\log N)^{1/p}/N$, $\lambda_2 \asymp (\log n)^{1/(2p)}/n$ and $(\log n)^{1/(2p)} \lesssim s \lesssim N^{1/2}(\log N)^{-\frac{7p+1}{2p}}$;

- **Polynomially decaying kernel (with $\mu_k \asymp k^{-2m}$)**

$$\eta_n^* = O(n^{-2m/(4m+1)})$$

given that $\lambda_1 \asymp N^{-2m/(2m+1)}$, $\lambda_2 \asymp n^{-4m/(4m+1)}$ and $s \asymp n^d$ for $(4m+1)^{-1} < d < (8m^2 + 4m)/(12m^2 - 2m + 1)$.

We remark that in the polynomial case, the minimal separation rate $n^{-2m/(4m+1)}$ turns out to be the minimax lower bound established by Ingster (1993) in the special case that $s = 1$, $Z \equiv 1$ and $g_0 \equiv 0$. The upper bounds of $s$ in the above three cases guarantee that the higher order terms of test statistics are asymptotically ignorable under the same conditions of Corollary B.3 in Appendix B. Rather, the lower bounds make sure that the nuisance estimator $\bar{g}_{N,\lambda}$ is close enough to the truth.

As for the test (c), we obtain similar conclusions as follows.

- **Finite rank kernel (with rank $r$)**
$$\widetilde{\eta}_N^* = O((sr)^{1/4}/\sqrt{n})$$
given that $\lambda_1 \asymp r/N, \lambda_2 \asymp \sqrt{r}/n$, and $\sqrt{r} \lesssim s \lesssim r^{-5/2} N^{1/2} (\log N)^{-4}$;

- **Exponentially decaying kernel (with $\mu_k \asymp e^{-\alpha k^p}$)**
$$\widetilde{\eta}_N^* = O(s^{1/4}(\log n)^{1/(4p)}/\sqrt{n})$$
given that $\lambda_1 \asymp (\log N)^{1/p}/N, \lambda_2 \asymp (\log n)^{1/(2p)}/n$ and $(\log n)^{1/(2p)} \lesssim s \lesssim N^{1/2}(\log N)^{-\frac{7p+1}{2p}}$;

- **Polynomially decaying kernel (with $\mu_k \asymp k^{-2m}$)**
$$\widetilde{\eta}_N^* = O(s^{m/(4m+1)} n^{-2m/(4m+1)})$$
given that $\lambda_1 \asymp N^{-2m/(2m+1)}$, $\lambda_2 \asymp (s/n^2)^{2m/(4m+1)}$ and $s \asymp n^d$ for $(2m+1)^{-1} < d < (8m^2 + 4m)/(12m^2 - 2m + 1)$.

We point out that only the $\lambda_2$ in the polynomial case is related to $s$. Also, the lower bound of $s$ in this case is larger than those in tests (a) and (b) due to the fact that a faster convergent $\bar{g}_{N,\lambda}$ is needed here.

**Remark 4.3.** Another test of interest is $H_0 : \beta_{0j} = c$ for some possibly unknown constant $c$. Under $H_0$, our main model (1.2) becomes a partially linear model estimated through

$$(\widehat{g}_c^{(j)}, \widehat{c}^{(j)}) = \underset{g,c}{\operatorname{argmin}} \left\{ \frac{1}{n} \sum_{i \in I_j} (Y_i - g(X_i) - cZ_i)^2 + \lambda_1 \|g\|_{\mathcal{H}}^2 \right\}.$$

Note that the joint asymptotic behavior of $(\widehat{g}_c^{(j)}, \widehat{c}^{(j)})$ for each fixed $j$ can be found in Cheng and Shang (2015). Similarly, our test statistic follows as

$$\text{PLRT}_{n,\lambda}^{\text{const}} = \mathcal{L}_{n,\lambda}^{(j)}(\bar{g}_{N,\lambda}, \widehat{\beta}_{n,\lambda}^{(j)}) - \mathcal{L}_{n,\lambda}^{(j)}(\bar{g}_{N,\lambda}, \widehat{c}^{(j)}),$$

where $\widehat{\beta}_{n,\lambda}^{(j)}$ and $\bar{g}_{N,\lambda}$ are defined in (3.4) and (3.6), respectively. By slightly modifying the proof of Theorem 4.1 (see Supplementary S.3.3), we have

$$-2nr'_T \cdot \text{PLRT}_{n,\lambda}^{\text{const}} \overset{a}{\sim} \chi^2_{u'_n}, \text{ where } r'_T = \sigma_T^2/(2\rho_P^2 + \rho_T^2),$$

$u'_n = h_2^{-1} \sigma_T^4/(2\rho_P^2 + \rho_T^2)$ and $\rho_P^2 = \sigma^4 h_1 \sum_{k=0}^{\infty} (1 + \lambda_1/\mu_{1,k})^{-2}$. Similar arguments also apply to composite hypotheses such as $H_0 : \beta_{0j}$ belongs a class of polynomial functions (with a specified degree).

**Remark 4.4.** To be consistent with the form of $\text{PLRT}_{n,\lambda}^{(j,j')}$, we can modify

PLRT$_{n,\lambda}^s$ by adding an additional term $\mathcal{L}_{n,\lambda}^{(s)}(\bar{g}_{N,\lambda}, \widehat{\beta}_{n,\lambda}^{(s)}) - \mathcal{L}_{n,\lambda}^{(s)}((\bar{g}_{N,\lambda}, \widehat{\beta}_{n,\lambda}^{(1)}))$ to the RHS of (4.2). Given this modification, all the above theoretical results remain valid after some minor changes, e.g., the degrees of freedom in (4.9).

4.2. *Homogeneity Testing.* In this section, we turn our attention to homogeneity testing

(4.7) $$H_0 : g_0 = g_0^*,$$

where $g_0^*$ is given. Since $g_0$ is shared among all sub-populations, we define an averaged penalized likelihood ratio test (APLRT) statistic as

(4.8) $$\text{APLRT}_{N,\lambda} = \frac{1}{s}\sum_{j=1}^{s}\left\{\mathcal{L}_{n,\lambda}^{(j)}(\bar{g}_{N,\lambda}, \widehat{\beta}_{n,\lambda}^{(j)}) - \mathcal{L}_{n,\lambda}^{(j)}(g_0^*, \widehat{\beta}_{n,\lambda}^{(j)})\right\}.$$

Note that the use of $\widehat{\beta}_{n,\lambda}^{(j)}$ in the second $\mathcal{L}_{n,\lambda}^{(j)}$ of (4.8) greatly saves computation cost than that using the restricted estimate of $\beta_{0j}$ under $H_0$, especially for a large $s$. We prove in Theorem 4.6 that the *aggregated* homogeneity test across $s$ sub-populations is more sensitive to local alternative sequences (in terms of a smaller separation rate) than that based on a single sub-population as $s \to \infty$. Also see more discussions after Theorem 4.6. And if we replace the *aggregated* estimate $\bar{g}_{N,\lambda}$ by $\widehat{g}_{n,\lambda}^{(j)}$ in each summand of (4.8), the separation rate will also be slowed down; see Remark 4.7. These discussions justify the form of (4.8), and again echo with the "blessing of aggregation" phenomenon.

The following theorem gives the null limiting distribution of APLRT$_{N,\lambda}$.

**Theorem 4.5.** Suppose that Assumptions 3.1 – 3.3, and (3.10) – (3.12) hold, and $E[\epsilon^4|X] \leq C$, a.s., for some constant $C$. If $nh_1\lambda_1 = o(1)$ and $(Nh_1^2)^{-1} = o(1)$, then we have

(4.9) $$-2Nr_P \cdot \text{APLRT}_{N,\lambda} \overset{a}{\sim} \chi_{u_N}^2$$

where $u_N = h_1^{-1}\sigma_P^4/\rho_P^2$, $r_P = \sigma_P^2/\rho_P^2$, and

$$\sigma_P^2 = \sigma^2 h_1 \sum_{k=0}^{\infty}(1 + \lambda_1/\mu_{1,k})^{-1}, \rho_P^2 = \sigma^4 h_1 \sum_{k=0}^{\infty}(1 + \lambda_1/\mu_{1,k})^{-2}.$$

Our next theorem derives a minimal separation rate between null and the alternative sequence under which the proposed homogeneity test statistic APLRT$_{N,\lambda}$ can still detect the latter. Define a separation rate as

(4.10) $$\eta_N^2 = \lambda_1 + \lambda_2 + (Nh_1^{1/2})^{-1}.$$

**Theorem 4.6.** Suppose the same conditions as Theorem 4.5. Under the alternative sequence $H_{1N} : g_0 = g_N$, if $r_n = o(1)$, we have for any $c_N \to \infty$,
(4.11)
$$\inf_{\substack{\{\beta_{0j}\}_{j=1}^s \subseteq \mathcal{F}_\zeta}} \inf_{\substack{\Delta g_N \in \mathcal{F}_\zeta \\ \|\Delta g_N\|^2 \geq c_N \eta_N^2}} \mathbb{P}_{H_{1N}} \left( \frac{|-2Nr_P \text{APLRT}_{N,\lambda} - u_N|}{\sqrt{2u_N}} > z_\alpha \right) \to 1,$$

where $\Delta g_N = g_N - g_0^*$. Here, "$\{\beta_{0j}\}_{j=1}^s \subseteq \mathcal{F}_\zeta$" means that each $\beta_{0j} \in \mathcal{F}_\zeta$ for $j = 1, \ldots, s$.

In what follows, we derive explicit conditions of $(\lambda_1, \lambda_2, s)$ for achieving the minimal separation rate, denoted as $\eta_N^*$, for three specific RKHS families.

- **Finite rank kernel (with rank $r$)**
$$\eta_N^* = O(r^{1/4}/\sqrt{N})$$
given that $\lambda_1, \lambda_2 \asymp \sqrt{r}/N$ and $s \lesssim r^{-5/2} N^{1/2} (\log N)^{-4}$;

- **Exponentially decaying kernel (with $\mu_k \asymp e^{-\alpha k^p}$)**
$$\eta_N^* = O\big([(\log N)^{1/(2p)}/N]^{1/2}\big)$$
given that $\lambda_1, \lambda_2 \asymp (\log N)^{1/(2p)}/N$ and $s \lesssim N^{1/2} (\log N)^{-\frac{7p+1}{2p}}$;

- **Polynomially decaying kernel (with $\mu_k \asymp k^{-2m}$)**
$$\eta_N^* = O(N^{-2m/(4m+1)})$$
given that $\lambda_1, \lambda_2 \asymp N^{-4m/(4m+1)}$ and $s \asymp n^d$ for $d < \frac{4m^2 - 7m + 1}{8m^2 + 2m}$.

Note that only an upper bound for $s$ is required in comparison with heterogeneity testing. The forms of the minimal separation rates $\eta_N^*$'s in the above three cases illustrate that the test statistic constructed based on the entire sample (with size $N$) is practically more powerful than that on any sub-sample (with size $n$), especially when $s \to \infty$. In the polynomial case, we remark that the minimal rate $\eta_N^*$ derived in presence of an increasing number of heterogeneity functions is actually minimax optimal acccording to Ingster (1993). We are able to achieve this rate optimality since the estimation error of each heterogeneity function vanishes, i.e., $r_n = o(1)$, under proper conditions on $s$. Note that $\lambda_2$ needs to be scaled in $N$ (rather than $n$) for obtaining the minimal separation rate. This is because we want to match with the magnitude of $\lambda_1$ in (4.10) in order to reduce the biases of heterogeneity estimators.

**Remark 4.7.** An alternative test statistic is formed by replacing $\bar{g}_{N,\lambda}$ with $\widehat{g}_{n,\lambda}^{(j)}$ in each summand of (4.8). Similar arguments in Theorem 4.5 imply the theoretical validity of this new test. However, in this case the separation rate, derived as $(\eta'_N)^2 = \lambda_1 + \lambda_2 + (s(nh_1)^{1/2})^{-1}$, becomes slower than that of APLRT$_{N,\lambda}$ in (4.10) due to the slower converging rate of $\widehat{g}^{(j)}_{n,\lambda}$ than $\bar{g}_{N,\lambda}$.

**5. Examples.** In this section, we study two examples in full details: the smoothing spline regression and the Gaussian kernel regression. Specifically, we compute the values of $(\sigma_g^2, \sigma_\beta^2)$ required in Theorem 3.4 and Corollary B.3, and the values of $(\rho_P^2, \rho_T^2, \sigma_P^2, \sigma_T^2)$ required in Theorems 4.1 and 4.5. Define two variance estimators as $\widehat{\sigma}^2 = N^{-1} \sum_{j=1}^{s} \sum_{i \in I_j} (Y_i - \bar{g}_{N,\lambda}(X_i) - \widehat{\beta}^{(j)}_{n,\lambda}(X_i))^2$ and $\widehat{\sigma}_z^2 = N^{-1} \sum_{i=1}^{N} Z_i^2$.

5.1. *Smoothing Spline Regression.* For $j = 1, \ldots, s$, we assume that $g_0, \beta_{0j} \in H_0^m[0,1]$ in our main model (1.2). Here, the periodic Sobolev space $H_0^m[0,1]$ is an RKHS space with $\|g\|_{\mathcal{H}}^2 = \int_0^1 (g^{(m)}(x))^2 dx$. For simplicity, we assume $X \sim U[0,1]$. In the above setup, the eigen-system is given in (2.2).

Based on Lemma S.9 in the Supplementary Material, we have $\lambda_1 = h_1^{2m}$, $\lambda_2 = h_2^{2m}$, and also we require the following constants:

Local CI: $\sigma_g^2 = \sigma^2 I_2/\pi, \sigma_\beta^2 = \sigma^2 I_2/(\sigma_z^{2-1/m}\pi)$,

Heter. test: $\sigma_T^2 = \pi^{-1}\sigma^2\sigma_z^{1/m}I_1, \rho_T^2 = \pi^{-1}\sigma^4\sigma_z^{1/m}I_2, r_T = \sigma^{-2}I_1/I_2$,

$u_n = \sigma_z^{1/m}h_2^{-1}I_1^2/(I_2\pi), \widetilde{u}_N = (4s/3)\sigma_z^{1/m}h_2^{-1}I_1^2/(I_2\pi)$,

Homo. test: $\rho_P^2 = \sigma^4 I_2/\pi, \sigma_P^2 = \sigma^2 I_1/\pi, r_P = \sigma^{-2}I_1/I_2, u_N = h_1^{-1}I_1^2/(I_2\pi)$,

where $I_1 = \pi/[2m\sin(\pi/(2m))]$ and $I_2 = (2m-1)\pi/[4m^2\sin(\pi/(2m))]$. We plug $\widehat{\sigma}$ and $\widehat{\sigma}_z$ into the above formulas as the estimate in practice.

5.2. *Gaussian Kernel Regression.* In this section, we explore the RKHS generated by a Gaussian kernel $K(x,t) = \exp(-|x-y|^2/2)$, denoted as $\mathcal{H}_G$. Again, we assume $g_0, \beta_{0j} \in \mathcal{H}_G$, $\epsilon \sim N(0, \sigma^2)$ and $X \sim N(0,1)$. The eigen-system of Gaussian kernel is given in (2.1). According to Krasikov (2004), we can get that

$$c_K = \sup_{k \in \mathbb{N}} \|\phi_k\|_{\sup} \leq \frac{2e^{15/8}(\sqrt{5}/4)^{1/4}}{3\sqrt{2\pi}2^{1/6}} \leq 1.336.$$

Hence, Assumption 3.2 holds.

From Lemma S.9, we have $h_1 = [\log(1/\lambda_1)]^{-1}$ and $h_2 = [\log(1/\lambda_2)]^{-1}$. Rewrite the eigen-value as $\mu_k = \ell \exp(2\ln\ell \cdot k)$, where recall that $\ell =$

$(\sqrt{5} - 1)/2$. Again, Lemma S.9 gives

$$\text{Local CI: } \sigma_g^2 \leq 1.7178\sigma^2, \sigma_\beta^2 \leq 1.7178\sigma^2/\sigma_z^2,$$
$$\text{Heter. test: } \sigma_T^2 = 2\sigma^{-2}\log(1/\ell), \rho_T^2 = 2\sigma^4\log(1/\ell), r_T = \sigma^{-2},$$
$$u_n = 2h_2^{-1}\log(1/\ell), \widetilde{u}_N = (8s/3)h_2^{-1}\log(1/\ell),$$
$$\text{Homo. test: } \sigma_P^2 = 2\sigma^{-2}\log(1/\ell), \rho_P^2 = 2\sigma^4\log(1/\ell),$$
$$r_P = \sigma^{-2}, u_N = 2h_1^{-1}\log(1/\ell).$$

The upper bounds for $\sigma_g^2$ and $\sigma_\beta^2$ are provided. Hence, we can construct *conservative* confidence intervals based on it.

**6. Simulations.** In this section, we provide numerical evidences to support our theory in the smoothing spline case. By the representer theorem (Schölkopf et al., 2001), the solution to (3.4) takes the form

$$\widehat{g}_{n,\lambda}^{(j)}(\cdot) = \sum_{k \in I_j} \alpha_k K(X_k, \cdot) \text{ and } \widehat{\beta}_{n,\lambda}^{(j)}(\cdot) = \sum_{k \in I_j} \eta_k K(X_k, \cdot).$$

Let $\mathbf{K}_j = (K(X_u, X_v))_{u,v \in I_j}, \mathbf{Y}_j = (Y_{i_1}, \ldots, Y_{i_n})^T$ and $\mathbf{Z}_j = \text{diag}(Z_{i_1}, \ldots, Z_{i_n})$ for $i_1, \ldots, i_n \in I_j$. We have

$$(6.1) \quad \left(\widehat{\boldsymbol{\alpha}}_j^T, \widehat{\boldsymbol{\eta}}_j^T\right)^T = \operatorname*{argmin}_{\boldsymbol{\alpha}, \boldsymbol{\eta}} n^{-1}\|\mathbf{Y}_j - \mathbf{K}_j\boldsymbol{\alpha} - \mathbf{Z}_j\mathbf{K}_j\boldsymbol{\eta}\|^2 + \lambda_1 \boldsymbol{\alpha}^T \mathbf{K}_j \boldsymbol{\alpha} + \lambda_2 \boldsymbol{\eta}^T \mathbf{K}_j \boldsymbol{\eta},$$

where $\widehat{\boldsymbol{\alpha}}_j = (\widehat{\alpha}_{j,i_1}, \ldots, \widehat{\alpha}_{j,i_n})^T, \widehat{\boldsymbol{\eta}}_j = (\widehat{\eta}_{j,i_1}, \ldots, \widehat{\eta}_{j,i_n})^T$. Solving (6.1) directly yields

$$\left(\widehat{\boldsymbol{\alpha}}_j^T, \widehat{\boldsymbol{\eta}}_j^T\right)^T = \left(\overline{\mathbf{Z}}_j^T \overline{\mathbf{K}}_j + n\boldsymbol{\Lambda}\right)^{-1} \overline{\mathbf{Z}}_j^T \mathbf{Y}_j,$$

where $\overline{\mathbf{Z}}_j = (\mathbf{I}_n, \mathbf{Z}_j), \overline{\mathbf{K}}_j = (\mathbf{K}_j, \mathbf{Z}_j \mathbf{K}_j)$ and $\boldsymbol{\Lambda} = \text{diag}(\lambda_1 \mathbf{I}_n, \lambda_2 \mathbf{I}_n)$. Therefore, we obtain the following closed forms

$$\widehat{g}_{n,\lambda}^{(j)}(\cdot) = \sum_{k \in I_j} \widehat{\alpha}_{j,k} K(X_k, \cdot) \text{ and } \widehat{\beta}_{n,\lambda}^{(j)}(\cdot) = \sum_{k \in I_j} \widehat{\eta}_{j,k} K(X_k, \cdot).$$

In the simulations, we set $X \sim U[0,1], Z \sim N(0,1)$ and $\epsilon \sim N(0,\sigma^2)$. Moreover, $\bar{g}_{N,\lambda}$ and $\{\widehat{\beta}_{n,\lambda}^{(j)}\}_{j=1}^s$ are estimated as cubic splines with $\lambda_1$ and $\lambda_2$ being chosen by cross validation.

6.1. *Local Confidence Interval.* We study the empirical performance of point-wise confidence interval (CI) in this section. The homogeneity function $g_0(x) = 0.6B(30, 17) + 0.4B(3, 11) - 1$, with $B(a, b)$ being a beta function, is plotted in Figure 1(a), while the heterogeneity functions are set as $\beta_{0j}(x) = \sin(2\pi j x)$, for $1 \leq j \leq s$. The standard deviation of the noise is chosen as $\sigma = 0.1$. We construct the pointwise CI at 95% level for $g_0$.

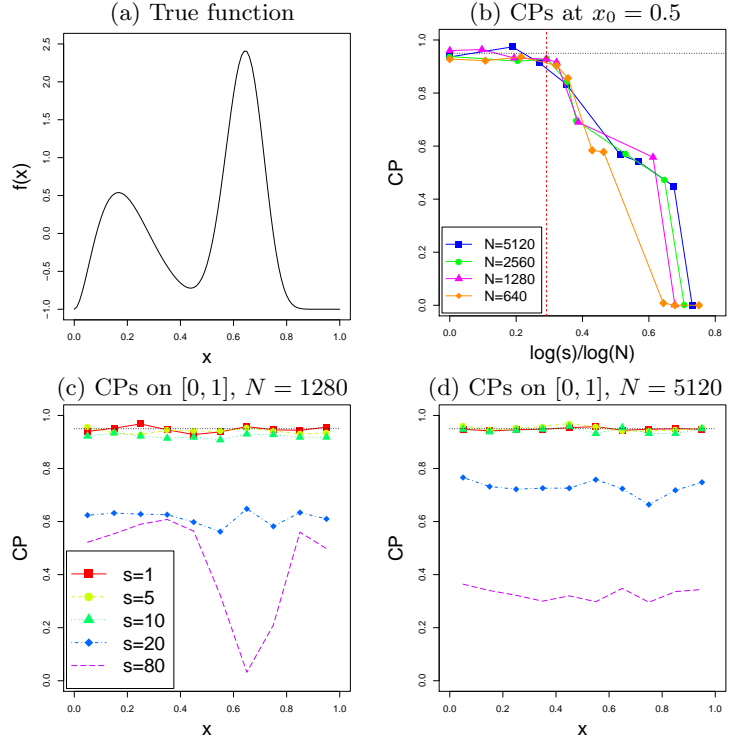

FIG 1. *The covering probabilities (CPs) of pointwise CIs with varying s.*

The covering probabilities of CIs are computed based on 500 repetitions for varying sample sizes $N = \{640, 1280, 2560, 5120\}$.

Figure 1(b) illustrates the covering probabilities of CIs at $x_0 = 0.5$ as $\log s / \log N$ increases. The red vertical line is the theoretical upper bound of $s$. The numerical results in Figure 1(b) demonstrate that the covering probability is close to 95% when $\log s / \log N$ is not very large, and decreases dramatically at the critical points around the red vertical line. Figures 1(c) and 1(d) give covering probabilities at 10 equally spaced points over $[0, 1]$ with the number of sub-populations ranging from 1 to 80. Similar "phase transition" phenomenon occurs in these two figures. For $s = 1, 5, 10$ which are above the critical value, the covering probabilities are around the nominal level, while for $s = 20$ and 80, they drop dramatically.

6.2. *Heterogeneity Testing.* For the heterogeneity testing, we focus on the constant test $H_0 : \beta_{0j} = $ constant and the simultaneous test $H_0 : \beta_{01} = \ldots = \beta_{0s}$ considered in Section 4.1. In the former test, we set $g_0(x) = 6x(x-1) + 1$ and $\beta_{01}(x) = c\cos(2\pi x)$ with $\sigma = 0.5$. The empirical power results are summarized in Table 1 for varying $c = \{0, 0.1, 0.5, 1\}$ and $n = $ $\{400, 600, 800, 1000\}$. In the latter test, we set $g_0(x) = 6x(x-1) + 1$ and $\beta_{01}(x) = \ldots = \beta_{0(s-1)}(x) = \sin(4\pi x)$ and $\beta_{0s}(x) = \sin(4\pi x) + c\cos(2\pi x)$ with $\sigma = 0.5$. The empirical power results are summarized in Table 2 for varying $c = \{0, 0.1, 0.5, 1\}$, $N = \{240, 420, 640, 1000\}$ and $s = 10$ respectively. Tables 1 and 2 both verify theoretical validity of our heterogeneity testing.

TABLE 1
Power performance of $\text{PLRT}_{n,\lambda}^{\text{const}}$ for the constant test

| $n$ | $c = 0.0$ | $c = 0.1$ | $c = 0.5$ | $c = 1.0$ |
|---|---|---|---|---|
| 400 | 0.06 | 0.95 | 0.97 | 1.00 |
| 600 | 0.05 | 0.98 | 0.99 | 1.00 |
| 800 | 0.07 | 0.99 | 1.00 | 1.00 |
| 1000 | 0.05 | 1.00 | 1.00 | 1.00 |

TABLE 2
Power performance of $\text{PLRT}_{n,\lambda}^{s}$ for the heterogeneity test

| $N$ | $c = 0.0$ | $c = 0.1$ | $c = 0.5$ | $c = 1.0$ |
|---|---|---|---|---|
| 640 | 0.07 | 0.75 | 0.83 | 0.98 |
| 1280 | 0.07 | 0.81 | 0.84 | 1.00 |
| 2560 | 0.05 | 0.84 | 0.94 | 1.00 |
| 5120 | 0.05 | 0.88 | 0.96 | 1.00 |

6.3. *Homogeneity Testing.* The homogeneity testing considered in Section 4.2 is studied in this section. Specifically, we set $g_0(x) = 6x(x-1) + 1 + c\cos(2\pi x)$ and $\beta_{0j}(x) = \sin(2\pi j x)$ with $\sigma = 0.1$. Consider a null hypothesis $H_0 : g_0 = 6x(x-1) + 1$. The empirical power results are demonstrated in Figure 2 with varying $c = \{0, 0.1, 0.2, 0.3\}$ and $N = \{640, 1280, 2560, 5120\}$. From these figures, we observe that the size is close to nominal level and the power is close to one when $s$ does not grow too fast w.r.t. $N$. Also, the empirical critical values are close to the theoretical results, denoted by the red vertical lines. As a side remark, we point out that powers under $c = 0.1, 0.2, 0.3$ return back to 1 when $s$ is very large. This is not surprising given very large values of statistics for those $s$. When the number of subpopulations is too large, the testing statistic becomes far deviated from the limiting distribution under both null and alternative hypotheses. Combining with the bad size under the null hypothesis, we can still conclude that the homogeneity tests do not perform satisfactorily when $s$ increases too fast.

## APPENDIX A: PROPERTIES OF RKHS

In this section, we present the proof of Proposition 2.1, and also study the properties of Fréchet derivatives of the loss function $\mathcal{L}_{n,\lambda}^{(j)}$ defined in (3.5).

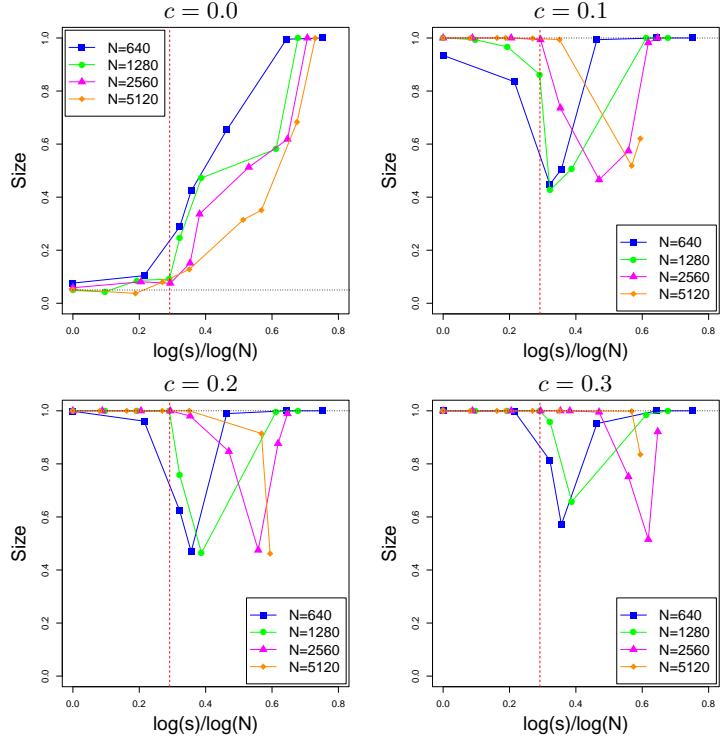

Fig 2. *Power comparison of* $\text{APLRT}_{N,\lambda}$ *for varying* $s$.

**A.1. Proof of Proposition 2.1.** We prove the existence of reproducing kernel by Riesz representation theorem (Riesz and Sz.-Nagy, 1955). To prove $(\mathcal{H}^2, \langle \cdot, \cdot \rangle)$ is an RKHS, it suffices to show that the evaluation functional $\delta_u$ is a bounded operator on any $u \in \mathcal{U}$. Recall that $K_1(x,y)$ and $K_2(x,y)$ are kernels for $\mathcal{H}_1$ and $\mathcal{H}_2$ respectively. We denote $K_{1,x}(\cdot) = K_1(x,\cdot)$ and $K_{2,x}(\cdot) = K_2(x,\cdot)$. For any $f = (g, \beta) \in \mathcal{H}^2$, we can bound the evaluation functional by

$$\begin{aligned}|\delta_u f| &= |g(x) + z\beta(x)| \leq |\langle K_{1,x}, g \rangle_{\mathcal{H}_1}| + c_z|\langle K_{2,x}, \beta \rangle_{\mathcal{H}_2}| \\ &\leq \|K_{1,x}\|_{\mathcal{H}_1}\|g\|_{\mathcal{H}_1} + c_z\|K_{2,x}\|_{\mathcal{H}_2}\|\beta\|_{\mathcal{H}_2} \\ &\leq (\lambda_1 \wedge \lambda_2)^{-1}(\|K_{1,x}\|_{\mathcal{H}_1} + c_z\|K_{2,x}\|_{\mathcal{H}_2})\|f\|,\end{aligned}$$

where the last inequality is due to the definition of the norm $\|\cdot\|$. Thus $\delta_u$ is a bounded operator and therefore $K^\lambda$ exists due to the Riesz representation theorem. The existence of the operator $W_\lambda$ can be proved following the same method of showing boundedness of $W_\lambda$. We present the proof in the Supplementary Material Section S.1.1

In the following proposition, we can further derive explicit expansions of the kernel function $K_u^\lambda$ and the operator $W_\lambda$. The proof is deferred to Supplementary Material Section S.1.1 as well.

**Proposition A.1.** We can write the kernel function as $K_u^\lambda = (P_x^\lambda, zT_x^\lambda)$ whose Fourier expansions are

$$\text{(A.1)} \qquad P_x^\lambda(t) = \sum_{k=0}^{\infty} \frac{\phi_k(x)\phi_k(t)}{1 + \lambda_1/\mu_{1,k}}, \ T_x^\lambda(t) = \sum_{k=0}^{\infty} \frac{\phi_k(x)\phi_k(t)}{\sigma_z^2 + \lambda_2/\mu_{2,k}}.$$

We also have the reproducing property that $\langle P_x^\lambda, g \rangle = g(x)$ and $\langle T_x^\lambda, \beta \rangle = \beta(x)$ for any $g \in \mathcal{H}_1, \beta \in \mathcal{H}_2$. If we write the operator $W_\lambda(g, \beta)(\cdot) = (W_{\lambda,1}g(\cdot), W_{\lambda,2}\beta(\cdot))$, we also have the following Fourier expansions

$$W_{\lambda,1}g(\cdot) = \sum_{k=0}^{\infty} \frac{\lambda_1/\mu_{1,k}}{1 + \lambda_1/\mu_{1,k}} \langle g, \phi_k \rangle_{L^2(\mathbb{P})} \phi_k(\cdot),$$

$$W_{\lambda,2}\beta(\cdot) = \sum_{k=0}^{\infty} \frac{\lambda_2/\mu_{2,k}}{\sigma_z^2 + \lambda_2/\mu_{2,k}} \langle \beta, \phi_k \rangle_{L^2(\mathbb{P})} \phi_k(\cdot).$$

**A.2. Fréchet Derivatives.** In this section, we calculate Fréchet derivatives of $\mathcal{L}_{n,\lambda}^{(j)}$. For notational simplicity, we use $\mathcal{L}_{n,\lambda}$ to represent $\mathcal{L}_{n,\lambda}^{(j)}$. The Fréchet derivative of $\mathcal{L}_{n,\lambda}(f)$ is a linear operator $D\mathcal{L}_{n,\lambda}(f) : \mathcal{H}^2 \to \mathbb{R}$, and the second Féchet derivative of $\mathcal{L}_{n,\lambda}(f)$ is a bilinear operator $D^2\mathcal{L}_{n,\lambda}(f_1, f_2) : \mathcal{H}^2 \times \mathcal{H}^2 \to \mathbb{R}$.

**Proposition A.2.** For any $f \in \mathcal{H}^2$, there exist operators $\nabla\mathcal{L}_{n,\lambda}(f) \in \mathcal{H}^2$ and $\nabla^2\mathcal{L}_{n,\lambda}(f) : \mathcal{H}^2 \to \mathbb{R}$ such that for any $f_1, f_2 \in \mathcal{H}^2$,

$$D\mathcal{L}_{n,\lambda}(f)[f_1] = \langle \nabla\mathcal{L}_{n,\lambda}(f), f_1 \rangle, \quad D^2\mathcal{L}_{n,\lambda}(f)[f_1, f_2] = \langle \nabla^2\mathcal{L}_{n,\lambda}(f)[f_1], f_2 \rangle.$$

In particular, we have the following closed forms of two operators

$$\text{(A.2)} \qquad \nabla\mathcal{L}_{n,\lambda}(f) = -\frac{1}{n}\sum_{i=1}^{n}(Y_i - g(X_i) - Z_i\beta(X_i))K_{U_i}^\lambda + W_\lambda f,$$

$$\text{(A.3)} \qquad \nabla^2\mathcal{L}_{n,\lambda}(f)[f_1] = \frac{1}{n}\sum_{i=1}^{n}\langle K_{U_i}^\lambda, f_1 \rangle K_{U_i}^\lambda + W_\lambda f_1.$$

Furthermore, we have $\nabla\mathcal{L}_{n,\lambda}(\widehat{f}_{n,\lambda}^{(j)}) = 0$ and $\mathbb{E}\nabla^2\mathcal{L}_{n,\lambda}(f_{0j})$ is an identity operator for any $j = 1, \ldots, s$.

We defer the proof to Supplementary Material Section S.1.2.

## APPENDIX B: LOCAL ASYMPTOTIC NORMALITY

In this section, we present the proof of Theorem 3.4 and also state Corollary B.3 that gives more explicit conditions in terms of $(s, \lambda_1, \lambda_2)$ for three important RKHS families.

**B.1. Proof of Theorem 3.4.** To show the joint asymptotic normality, we decompose our estimator into a leading term and a remainder term as

(B.1)
$$\bar{g}_{N,\lambda}(x) - g_0^u(x) = \underbrace{\frac{1}{N}\sum_{i=1}^{N} \epsilon_i P_{X_i}^\lambda(x)}_{\text{Leading term}} + \underbrace{G_N(x)}_{\text{Remainder term}},$$

$$\widehat{\beta}_{n,\lambda}^{(j)}(x) - \beta_{0j}^u(x) = \underbrace{\frac{1}{n}\sum_{i \in I_j} \epsilon_i Z_i T_{X_i}^\lambda(x)}_{\text{Leading term}} + \underbrace{H_n^{(j)}(x)}_{\text{Remainder term}}.$$

Recall that $(g_0^u, \beta_0^u) := f_{0j} + W_\lambda f_{0j}$, and $P_x^\lambda$ and $T_x^\lambda$ are defined in (A.1).

A key step in our proof is the following lemma on the rate of $R_N^{(j)} := (G_N, H_n^{(j)}) \in \mathcal{H}^2$.

**Lemma B.1.** Under the same assumptions as Theorem 3.4, for $\xi_\lambda$ defined in (3.8), we have,
$$\|R_N^{(j)}\| = o_P\left((nh)^{-1}(\sqrt{n}\xi_\lambda + 1)\left((\log n)(nh)^{-1/2} + \lambda_1^{1/2} + \lambda_2^{1/2}\right)\right).$$
for any $j = 1, \ldots, s$.

We defer the proof to Section S.2.1.

Combining Lemma B.1 with (3.11), we have $\|R_N^{(j)}\| = o_P(N^{-1/2})$. Define

(B.2)
$$I_1 := \sqrt{Nh_1} \cdot G_N(x_1) = \sqrt{Nh_1}\langle P_{x_1}^\lambda, G_N\rangle,$$
$$I_2 := \sqrt{nh_2} \cdot H_n^{(j)}(x_2) = \sqrt{nh_2}\langle T_{x_2}^\lambda, H_n^{(j)}\rangle.$$

In order to bound $I_1$ and $I_2$, we need the following lemma on the operator norms of $P_x^\lambda$ and $T_x^\lambda$.

**Lemma B.2.** Under Assumption 3.2, we have $\|P_x^\lambda\|^2 \leq c_K h_1^{-1}$ and $\|T_x^\lambda\|^2 \leq c_K h_2^{-1}$. Recalling that $c_K = c_\phi\sqrt{1+c_z}$ is a universal constant defined in Assumptions 3.1 and 3.2.

The proof is deferred to Supplementary Material Section S.1.3.

To apply the Cramer-Wold's device, we need to bound any linear combination of $I_1$ and $I_2$. Combining (B.2) with Lemma B.2, we have for any $t_1, t_2 \in \mathbb{R}$,

(B.3)
$$|t_1 I_1 + t_2 I_2| \leq |t_1|(Nh_1)^{1/2}\|P_{x_1}^\lambda\|\|G_N\| + |t_2|(nh_2)^{1/2}\|T_{x_1}^\lambda\|\|H_n^{(j)}\|$$
$$\leq \sqrt{c_K}(|t_1| + |t_2|)\sqrt{N}\|R_N^{(j)}\| = o_P(1).$$

Next we turn to calculate the variances of the leading terms. Following the Fourier expansion of $P_x^\lambda$ in Proposition A.1, it can be derived that

$$\text{Var}\left(\sum_{i=1}^{N}\epsilon_i P_{X_i}^\lambda(x_1)\right) = N\mathbb{E}\bigl[\epsilon^2|P_X^\lambda(x_1)|^2\bigr] = \sigma^2 N \sum_{k=0}^{\infty}\left(\frac{\phi_k(x_1)}{1+\lambda_1/\mu_{1,k}}\right)^2,$$

and similarly, we have

$$\text{Var}\left(\sum_{i=1}^{n}\epsilon_i Z_i T_{X_i}^\lambda(x_2)\right) = \sigma^2\sigma_z^2 n \sum_{k=0}^{\infty}\left(\frac{\phi_k(x_2)}{\sigma_z^2+\lambda_2/\mu_{2,k}}\right)^2.$$

As $\mathbb{E}[Z_i] = 0$, the covariance between two leading terms is

$$\text{Cov}\left(\sum_{i=1}^{N}\epsilon_i P_{X_i}^\lambda(x_1), \sum_{i\in I_j}\epsilon_i Z_i T_{X_i}^\lambda(x_2)\right) = 0, \text{ for any } j=1,\ldots,s.$$

Therefore, for any $(t_1, t_2) \in \mathbb{R}^2$, by (B.3) and CLT, we have

$$(t_1, t_2)\begin{pmatrix}\sqrt{Nh_1}\left(\bar{g}_{N,\lambda}(x_1) - g_0^u(x_1)\right)\\ \sqrt{nh_2}\left(\widehat{\beta}_{n,\lambda}^{(j)}(x_2) - \beta_{0j}^u(x_2)\right)\end{pmatrix} \rightsquigarrow N\bigl(0, t_1^2\sigma_g^2 + t_2^2\sigma_h^2\bigr),$$

where $\sigma_g^2, \sigma_h^2$ are defined in (3.9). Applying the Cramer-Wold's device, the joint normality in (3.13) follows.

**B.2. Corollary B.3.** In this corollary, we assume $\mathcal{H}_1 = \mathcal{H}_2$, and that the estimation bias $W_\lambda f_{0j}$ is removed through under-smoothing for simplicity.

**Corollary B.3.** Suppose Assumptions 3.1 and 3.2 hold and (3.9) is satisfied. We have the following joint asymptotic distribution:

(B.4) $$\begin{pmatrix}\sqrt{Nh_1}\left(\bar{g}_{N,\lambda}(x_1) - g_0(x_1)\right)\\ \sqrt{nh_2}\left(\widehat{\beta}_{n,\lambda}^{(j)}(x_2) - \beta_{0j}(x_2)\right)\end{pmatrix} \rightsquigarrow N\left(\mathbf{0}, \begin{pmatrix}\sigma_g^2 & 0\\ 0 & \sigma_\beta^2\end{pmatrix}\right).$$

given that $(s, \lambda_1, \lambda_2)$ satisfy the following criteria:

1. **Finite Rank:** $\mu_k = 0$ for $k \geq r$.

$$\lambda_1 = o\bigl(\sqrt{r/N}\bigr) \text{ and } \log(\lambda_1^{-1}) = o\left(\log^2 N\right),$$
$$\lambda_2 = o\bigl(\sqrt{r/n}\bigr) \text{ and } \log(\lambda_2^{-1}) = o\left(\log^2 n\right),$$

(B.5) $$s = o\left(\frac{N^{1/2}}{r^{5/2}[\log(\lambda_1^{-1})]^{1/2}(\log N)^3}\right).$$

2. **Exponentially Decaying:** $\mu_k = O(e^{-\alpha k^p})$.

$$\lambda_1 = o\left(\sqrt{1/N}\right) \text{ and } \log(\lambda_1^{-1}) = o\left(\log^2 N\right),$$

$$\lambda_2 = o\left(\sqrt{1/n}\right) \text{ and } \log(\lambda_2^{-1}) = o\left(\log^2 n\right),$$

(B.6) $$s = o\left(\frac{N^{1/2} h^{3/2}}{[\log(h\lambda_1^{-1})]^{(p+1)/2p}(\log N)^3}\right).$$

3. **Polynomially Decaying:** $\mu_k = O(k^{-2m})$. Given some $d$ satisfying $2m/(4m+1) < d < 4m^2/(8m-1)$, we have

$$\lambda_1 \asymp N^{-d}, \lambda_2 \asymp n^{-d}$$

(B.7) $$s = o\left(N^{\frac{1}{2} - \frac{8m-1}{8m^2}d}\right).$$

**Remark B.4.** Corollary B.3 also holds even when $\mathcal{H}_1 \neq \mathcal{H}_2$. For example, if $\mu_{1,k} \asymp e^{-\alpha k^p}$ in $\mathcal{H}_1$ and $\mu_{2,k} \asymp k^{-2m}$ in $\mathcal{H}_2$, we need to make the same choice of $s$ as (B.7), which is irrelevant to $p$. This is mainly because the estimation error of $\widehat{\beta}_{n,\lambda}^{(j)}$ dominates that of $\bar{g}_{N,\lambda}$.

PROOF. Since $\mathcal{H}_1 = \mathcal{H}_2$, we simply denote them by $\mathcal{H} := \mathcal{H}_1 = \mathcal{H}_2$. To apply Theorem 3.4, we first discuss the magnitude of the asymptotic bias $f_{0j}^u - f_{0j} = W_\lambda f_{0j} = (W_{\lambda,1} g_0, W_{\lambda,2} \beta_{0j})$. By (S.23) and Lemma S.1 in the Supplementary Material, we have for any $x_0 \in \mathcal{X}$,

$$\sqrt{Nh_1} W_{\lambda,1} g_0(x_0) \leq \sqrt{Nh_1} \|P_x^\lambda\| \|W_{\lambda,1} g_0\| = o(\sqrt{N\lambda_1}),$$
$$\sqrt{nh_2} W_{\lambda,2} \beta_{0j}(x_0) \leq \sqrt{nh_2} \|T_x^\lambda\| \|W_{\lambda,2} \beta_{0j}\| = o(\sqrt{n\lambda_2}).$$

Therefore, if $\lambda_1 = o(\sqrt{1/N})$ and $\lambda_2 = o(\sqrt{1/n})$, the asymptotic bias terms are eliminated, i.e., $\sqrt{Nh_1} W_{\lambda,1} g_0(x_0) = o(1)$ and $\sqrt{nh_2} W_{\lambda,2} \beta_{0j}(x_0) = o(1)$.

We next bound the covering number of $\mathcal{H}_0$ defined in (3.2). Define $\mathbb{B}_\mathcal{H}(1) = \{g \in \mathcal{H} : \|g\|_\mathcal{H} \leq 1\}$. Suppose $\{g_1, \ldots, g_{K_1}\}$ and $\{\beta_1, \ldots, \beta_{K_2}\}$ are the $\epsilon/2$-net and $\epsilon/(2c_z)$-net of $B_\mathcal{H}(1)$ in terms of $\|\cdot\|_{\sup}$, respectively. For any $f = (g, \beta) \in \mathcal{H}_0$, there exist $g_i, \beta_j$ such that $\|g_i - g\|_{\sup} \leq \epsilon/2$ and $\|\beta_j - \beta\|_{\sup} \leq \epsilon/(2c_z)$ and therefore we have $\|f - (g_i, \beta_j)\|_{\sup} \leq \epsilon$. This implies that $\{g_1, \ldots, g_{K_1}\} \times \{\beta_1, \ldots, \beta_{K_2}\}$ is an $\epsilon$-net of $\mathcal{H}_0$. We thus have

$$N(\mathcal{H}_0; \|\cdot\|_{\sup}, \epsilon) \leq N(B_\mathcal{H}(1); \|\cdot\|_{\sup}, \epsilon/2) N(B_\mathcal{H}(1); \|\cdot\|_{\sup}, \epsilon/(2c_z))$$
(B.8) $$\leq N(B_\mathcal{H}(1); \|\cdot\|_\mathcal{H}, \epsilon/(2c_b)) N(B_\mathcal{H}(1); \|\cdot\|_\mathcal{H}, \epsilon/(2c_z c_b)),$$

where the last inequality follows from the fact that $\|g\|_{\sup} \leq c_b \|g\|_\mathcal{H}$ for any $g \in \mathcal{H}$. Recall that $c_b$ is defined in Assumption 3.2.

In what follows, we will determine the upper bound of $s$ by checking (3.10), (3.11) and (3.12) in three concrete kinds of eigenvalue decaying rates.

- **Finite Rank**

Before verifying (3.10), (3.11) and (3.12), we first bound $h$, $\omega(\delta)$ in (3.3) and $\xi_\lambda$ in (3.8). Since $\lambda_1 = o(1)$ and $\lambda_2 = o(1)$, we have $\gamma_1 = \sum_{k=0}^{r}(1+\lambda_1/\mu_{1,k})^{-1} \in (r/2, r)$ which derives $h_1 = 1/\gamma_1 \in (1/r, 2/r)$. For the same reason, we get $h_1 = 1/\gamma_2 \in (1/r, 2/r)$. We also have $\lambda = \lambda_1 \wedge \lambda_2 = o(\sqrt{r/N})$ as $n \ll N$. According to Lemma S.2 and (B.8), we can quantify the integral $\omega(\delta)$ such that for $\delta \in (0,1)$,

$$\omega(\delta) \lesssim \int_0^\delta r\sqrt{\log(1/\epsilon)}d\epsilon \leq r\Big(\int_0^\delta d\epsilon \int_0^\delta \log(1/\epsilon)d\epsilon\Big)^{1/2}$$

(B.9)
$$= r\sqrt{\delta^2(\log(1/\delta)+1)} \leq Cr\delta\sqrt{\log(1/\delta)},$$

where the second inequality is due to Cauchy-Schwartz inequality. We next bound $\xi_\lambda$. Since $\lambda_1/h = o(1)$, we have

$$\xi_\lambda = (c_K^{-2}h\lambda_1^{-1})^{1/2}\omega((c_K^{-2}h\lambda_1^{-1})^{-1/2}) \lesssim r\sqrt{\log(h\lambda_1^{-1})}.$$

By (B.9), we can construct the upper bound in Assumption 3.3 as $\bar\omega_a(\delta) = Cra\sqrt{\log(1/(a\delta))}$ for some generic constant $C > 0$, which is non-increasing.

Solving (3.11), we have

(B.10) $$s = o\left(\frac{N^{1/2}}{r^{5/2}(\log(h\lambda_1^{-1}))^{1/2}(\log N)^3}\right).$$

Since $\lambda_1 = o(\sqrt{r/N})$ and $\log(1/\lambda_1) = o(\log^2 N)$, we have $\log(N\xi_\lambda/(h\lambda)) = o(\log^2 N)$ such that (3.10) is satisfied. We can verify (3.12) by

$$r_n^{-1}\omega\big((\lambda_1^{-1}h)^{-1/2}r_n\big)(\lambda_1 h)^{-1/2}n^{-1/2}\log N = o(1).$$

- **Exponentially Decaying**

For the exponentially decaying kernel, following Lemma S.2 and (B.8) again, we bound $\omega(\delta)$ by

$$\omega(\delta) \lesssim \int_0^\delta [\log(1/\epsilon)]^{\frac{p+1}{2p}} d\epsilon \lesssim \delta\,[\log(1/\delta)]^{-\frac{p+1}{2p}}.$$

We therefore construct the upper bound required in Assumption 3.3 as $\bar\omega_a(\delta) = Ca\,[\log(1/(a\delta))]^{-\frac{p+1}{2p}}$. It is easy to check $\bar\omega_a(\delta)$ is non-increasing. We also have $\lambda_1 = \lambda_1 \wedge \lambda_2$. We next bound $\xi_\lambda$ as

$$\xi_\lambda = (c_K^{-2}h\lambda_1^{-1})^{1/2}\omega((c_K^{-2}h\lambda_1^{-1})^{-1/2}) \lesssim [\log(h\lambda_1^{-1})]^{\frac{p+1}{2p}}.$$

Following the proof of Lemma S.9 in the Supplementary Material, we have

$$1/h_1 = \gamma_1(\lambda_1) = \sum_{k=0}^{\infty}\frac{1}{1+\lambda_1/\mu_{1,k}} \asymp (\log(1/\lambda_1))^{1/p},$$

and the same for $h_2$. We now verify (3.10) - (3.12). We obtain the upper bound of $s$ by solving (3.11):

$$(B.11) \qquad s = o\left(\frac{N^{1/2}h^{3/2}}{[\log(h\lambda_1^{-1})]^{(p+1)/2p}(\log N)^3}\right).$$

Since $\lambda_1 = o(\sqrt{1/N})$ and $\lambda_2 = o(\sqrt{1/n})$, it is easy to check $\log(N\xi_\lambda(h\lambda_1)^{-1}) = o(\log^2 N)$ such that (3.10) is satisfied. We can easily check (3.12) that

$$r_n^{-1}\omega\big((\lambda_1^{-1}h)^{-1/2}r_n\big)(\lambda_1 h)^{-1/2}n^{-1/2}\log N = o(1).$$

- **Polynomially Decaying**

For the polynomially decaying kernel, similar to the above two cases, we apply Lemma S.2 and (B.8) and get

$$\omega(\delta) \lesssim \int_0^\delta (1/\epsilon)^{1/(2m)}\,d\epsilon \lesssim \delta^{1-1/(2m)}.$$

We then construct a non-increasing function $\bar{\omega}_a(\delta) = Ca^{1-1/2m}\delta^{-1/(2m)}$ as the upper bound required in Assumption 3.3. Also $\xi_\lambda$ can be bounded as

$$\xi_\lambda = (c_K^{-2}h\lambda_1^{-1})^{1/2}\omega((c_K^{-2}h\lambda_1^{-1})^{-1/2}) \lesssim (h\lambda_1^{-1})^{1/(4m)}$$

Moreover, from the proof of Lemma S.9, we have $\gamma(\lambda_1) = \sum_{k=0}^s (1 + \lambda_1/\mu_{1,k})^{-1} \asymp \lambda_1^{-1/(2m)}$ and the same for $\gamma(\lambda_2)$. Solving (3.11), we have

$$(B.12) \qquad s = o\left(N^{1/2}\lambda_1^{\frac{8m-1}{8m^2}}\right),$$

and $\lambda_1^{-1}N^{-\frac{4m^2}{8m-1}} = o(1)$ as $s$ is always larger than 1.

As $\lambda_1 = o(N^{-1/2})$, we have $\log(N\xi_\lambda(h\lambda_1)^{-1}) = o(\log^2 N)$ such that (3.10) is satisfied. We then verify (3.12):

$$r_n^{-1}\omega\big((\lambda_1^{-1}h)^{-1/2}r_n\big)(\lambda_1 h)^{-1/2}n^{-1/2}\log N = o(1)$$

gives $r_n^{1/(2m)} \lesssim \lambda_1^{1/(8m^2)}n^{-1/2}$ and therefore we have

$$n \gtrsim \left(\lambda_1^{-1/(m-2m^2)}(\log N)^{-2/(2m-1)}\right) \vee \lambda_1^{-\frac{2m-1}{4m^2}}.$$

This gives the upper bound $s \lesssim N^{1/2}\lambda_1^{\frac{2m-1}{8m^2}}$, which is weaker than (B.12). □

## APPENDIX C: PROOFS FOR HYPOTHESIS TESTS

In this section, we present partial proofs for Theorems 4.1 and 4.2, in which only $\mathrm{PLRT}_{n,\lambda}^{(j,j')}$ and $\mathrm{PLRT}_{n,\lambda}^s$ are considered. The remaining proofs for $\mathrm{PLRT}_{n,\lambda}$ and $\mathrm{PLRT}_{n,\lambda}^{\mathrm{const}}$ can be found in the Supplementary Material Section S.3.3. In the end, we prove Theorems 4.5 and 4.6.

**C.1. Partial Proofs of Theorems 4.1 and 4.2.**

C.1.1. *Analysis of* $\text{PLRT}_{n,\lambda}^{(j,j')}$. Without loss of generality, we assume $j=1, j'=2$. We can decompose $\text{PLRT}_{n,\lambda}^{(1,2)} = L_{(1)} + L'_{(1)}$, where

$$L_{(1)} := \mathcal{L}_{n,\lambda}^{(1)}\big((\bar{g}_{N,\lambda}, \widehat{\beta}_{n,\lambda}^{(1)})\big) - \mathcal{L}_{n,\lambda}^{(1)}\big((\bar{g}_{N,\lambda}, \widehat{\beta}_{n,\lambda}^{(2)})\big),$$

$$L'_{(1)} := \mathcal{L}_{n,\lambda}^{(2)}\big((\bar{g}_{N,\lambda}, \widehat{\beta}_{n,\lambda}^{(2)})\big) - \mathcal{L}_{n,\lambda}^{(2)}\big((\bar{g}_{N,\lambda}, \widehat{\beta}_{n,\lambda}^{(1)})\big).$$

The high level idea of our proof is to use Taylor expansion to decompose $L_{(1)}$ and $L'_{(1)}$ into the leading asymptotic normal terms and higher order terms.

To implement this idea, we first study the derivative of operators in the above two terms. Denote the second coordinate of the derivative $\nabla \mathcal{L}_{n,\lambda}^{(1)}((g_0, \beta_{01}))$ in (A.2) as

$$\nabla \mathcal{L}_2^{(1)}(\beta_{01}) := -\frac{1}{n} \sum_{i \in I_1} \epsilon_i Z_i T_{X_i}^\lambda + W_{2,\lambda} \beta_{01}.$$

Similarly, we define the second coordinate of $\nabla \mathcal{L}_{n,\lambda}^{(2)}((g_0, \beta_{02}))$ as

$$\nabla \mathcal{L}_2^{(2)}(\beta_{02}) := -\frac{1}{n} \sum_{i \in I_2} \epsilon_i Z_i T_{X_i}^\lambda + W_{2,\lambda} \beta_{02}.$$

Define $\beta_{1,2} = \widehat{\beta}_{n,\lambda}^{(1)} - \widehat{\beta}_{n,\lambda}^{(2)}$. We have the following Taylor expansion

$$L_{(1)} = \langle \nabla \mathcal{L}_{n,\lambda}^{(1)}(\bar{g}_{N,\lambda}, \widehat{\beta}_{n,\lambda}^{(2)}), \beta_{1,2}\rangle + \frac{1}{2}\nabla^2 \mathcal{L}_{n,\lambda}^{(1)}(\bar{g}_{N,\lambda}, \widehat{\beta}_{n,\lambda}^{(2)})[\beta_{1,2}, \beta_{1,2}]$$

$$(\text{C.1}) \quad = \langle \nabla \mathcal{L}_{n,\lambda}^{(1)}(g_0, \widehat{\beta}_{n,\lambda}^{(2)}), \beta_{1,2}\rangle + \frac{1}{2}\mathbb{E}\{\nabla^2 \mathcal{L}_{n,\lambda}^{(1)}(g_0, \beta_{01})\}[\beta_{1,2}, \beta_{1,2}] + \Delta_1/2,$$

where $\Delta_1 = \big(\nabla^2 \mathcal{L}_{n,\lambda}^{(1)}(\bar{g}_{N,\lambda}, \beta_{02}) - \mathbb{E}\{\nabla^2 \mathcal{L}_{n,\lambda}^{(1)}(g_0, \beta_{01})\}\big)[\beta_{1,2}, \beta_{1,2}]$. Following a similar proof as Lemma B.1, we have $\Delta_1 = o_P((n^2 h_2)^{-1/2})$. We refer the detailed derivation to the proof of (S.34) in the Supplementary Material.

Combining (C.1) with the fact that $\mathbb{E}\{\nabla^2 \mathcal{L}_{n,\lambda}^{(1)}(g_0, \beta_{01})\}$ is an identity as shown in Proposition A.2, we have

$$(\text{C.2}) \quad L_{(1)} = \langle \nabla \mathcal{L}_{n,\lambda}^{(1)}(g_0, \widehat{\beta}_{n,\lambda}^{(2)}), \beta_{1,2}\rangle + \frac{1}{2}\|\beta_{1,2}\|^2 + o_P((n^2 h_2)^{-1/2}),$$

By (B.1), we have $\beta_{1,2} = \nabla \mathcal{L}_2^{(2)}(\beta_{02}) - \nabla \mathcal{L}_2^{(1)}(\beta_{01}) + (H_n^{(1)} - H_n^{(2)})$. Applying Lemma B.1 and (3.11), we have $\|H_n^{(1)} - H_n^{(2)}\| = o_P(N^{-1/2})$ and

$$(\text{C.3}) \quad \|\beta_{1,2}\| = \|\nabla \mathcal{L}_2^{(2)}(\beta_{02}) - \nabla \mathcal{L}_2^{(1)}(\beta_{01})\| + o_P(N^{-1/2}).$$

Applying the Taylor expansion again, since $\beta_{01} = \beta_{02}$, we get

$$\mathcal{L}_{n,\lambda}^{(1)}(g_0, \widehat{\beta}_{n,\lambda}^{(2)}) = \langle \nabla \mathcal{L}_{n,\lambda}^{(1)}(g_0, \beta_{01}), \beta_{1,2} \rangle + \nabla^2 \mathcal{L}_{n,\lambda}^{(1)}(g_0, \beta_{01})[\widehat{\beta}_{n,\lambda}^{(2)} - \beta_{02}, \beta_{1,2}]$$

$$= \langle \nabla \mathcal{L}_2^{(1)}(\beta_{01}), \beta_{1,2} \rangle + \langle \widehat{\beta}_{n,\lambda}^{(2)} - \beta_{02}, \beta_{1,2} \rangle + o_P((n^2 h_2)^{-1/2})$$

(C.4) $$= -\|\nabla \mathcal{L}_2^{(1)}(\beta_{01}) - \nabla \mathcal{L}_2^{(2)}(\beta_{02})\|^2 + o_P((n^2 h_2)^{-1/2}),$$

where the second equality follows a similar argument as (C.2) and the third equality is similar to (C.3). The detailed proof follows a similar method in proving Lemma C.1 in the Supplementary Material.

Combining (C.2) – (C.4), we have

(C.5) $$-2nL_{(1)} = n\|\nabla \mathcal{L}_2^{(1)}(\beta_{01}) - \nabla \mathcal{L}_2^{(2)}(\beta_{02})\|^2 + o_P(h_2^{-1/2}),$$

where the bias term $W_{2,\lambda}\beta_{01}$ is ignorable due to the same derivation of (S.36) in the Supplementary Material.

Similarly, we also have $-2nL'_{(1)} = n\|\nabla \mathcal{L}_2^{(1)}(\beta_{01}) - \nabla \mathcal{L}_2^{(2)}(\beta_{02})\|^2 + o_P(h_2^{-1/2})$. Therefore, we have

(C.6) $$-2n\text{PLRT}_{n,\lambda}^{(1,2)} = 2n\|\nabla \mathcal{L}_2^{(1)}(\beta_{01}) - \nabla \mathcal{L}_2^{(2)}(\beta_{02})\|^2 + o_P(h_2^{-1/2}).$$

Notice that $\nabla \mathcal{L}_2^{(1)}(\beta_{01}) - \nabla \mathcal{L}_2^{(2)}(\beta_{02}) \stackrel{d}{=} n^{-1} \sum_{i \in I_1 \cup I_2} \epsilon_i Z_i T_{X_i}^\lambda$. We expand $\|\sum_{i \in I_1 \cup I_2} \epsilon_i Z_i T_{X_i}^\lambda\|_2^2 = W_3(2n) + \Delta_2$, where the cross term can be written as

(C.7) $$W_3(2n) = \sum_{i<j, i,j \in I_1 \cup I_2} 2\epsilon_i \epsilon_j Z_i Z_j T^\lambda(X_i, X_j).$$

The square term $\Delta_2$ has the rate

$$\frac{1}{n} \sum_{i \in I_1 \cup I_2} \epsilon_i^2 Z_i^2 T^\lambda(X_i, X_i) = 2h_2^{-1}\sigma_T^2 + O_P(1).$$

The detailed analysis can be referred to a similar procedure in (S.37) in the Supplementary Material.

The next step is to show the asymptotic normality of $n^{-1}W_3(2n)$. By denoting $W_{ij} = 2\epsilon_i \epsilon_j Z_i Z_j T^\lambda(X_i, X_j)$, we rewrite $W_3(2n) = \sum_{i<j, i,j \in I_1 \cup I_2} W_{ij}$. We next apply Proposition 3.2 in de Jong (1987) by checking several conditions. First, it can be easily checked that $W(N)$ is clean in the sense that $\mathbb{E}(W_{ij}|X_i) = 0$ for all $i < j$. Next, we need to bound several moments of $W_3(2n)$. Define $\sigma(2n)^2 = \text{Var}(W_3(2n))$ and the following quantities

$$G_{\text{I}} = \sum_{i<j} \mathbb{E}[W_{ij}^4], \quad G_{\text{II}} = \sum_{i<j<k} (\mathbb{E}[W_{ij}^2 W_{ik}^2] + \mathbb{E}[W_{ji}^2 W_{jk}^2] + \mathbb{E}[W_{ki}^2 W_{kj}^2]),$$

$$G_{\text{III}} = \sum_{i<j<k<l} (\mathbb{E}[W_{ij}W_{ik}W_{lj}W_{lk}] + \mathbb{E}[W_{ij}W_{il}W_{kj}W_{kl}] + \mathbb{E}[W_{ik}W_{il}W_{jk}W_{jl}]).$$

By Lemma B.2 and $\mathbb{E}[\epsilon^4|X] < \infty$, we have
$$\mathbb{E}[W_{ij}^4] \lesssim \mathbb{E}[\epsilon_i^4 \epsilon_j^4 T^\lambda(X_i, X_j)^4] = O(h_2^{-4}),$$
$$\mathbb{E}[W_{ij}^2 W_{ik}^2] \leq \mathbb{E}[W_{ij}^4] = O(h_2^{-4}),$$
$$\mathbb{E}[W_{ij} W_{ik} W_{lj} W_{lk}] = O(h_2^{-1}).$$

It follows from the results above that
$$G_{\mathrm{I}} = O(n^2 h_2^{-4}),\ G_{\mathrm{II}} = O(n^3 h_2^{-4}) \text{ and } G_{\mathrm{III}} = O(n^4 h_2^{-1}).$$

We next obtain the exact order of $\sigma(2n)^4$. This follows from $\mathbb{E}[W_{ij}^2] = 4h_2^{-1}\rho_T^2$ that $\sigma(2n)^4 = \left(\binom{2n}{2}\mathbb{E}[W_{ij}^2]\right)^2 \sim (8h_2^{-1} n^2 \rho_T^2)^2$. Since $(nh_2^2)^{-1} = o(1)$, we have $G_{\mathrm{I}}, G_{\mathrm{II}}$ and $G_{\mathrm{III}}$ are of smaller order than $\mathrm{Var}^2(W(2n)) = O(n^4 h_2^{-2})$. Therefore, it follows from Proposition 3.2 in de Jong (1987) that

$$\text{(C.8)} \qquad \frac{1}{\sqrt{2/h_2} \cdot 2n\rho_T} W_3(2n) \rightsquigarrow N(0,1).$$

Combining (C.6) with (C.8), we have
$$(8h_2^{-1}\sigma_T^4/\rho_T^2)^{-1/2}(-nr_T \cdot \mathrm{PLRT}_{n,\lambda}^{(j,j')} - 2h_2^{-1}\sigma_T^4/\rho_T^2) \rightsquigarrow N(0,1).$$

We then focus on the power analysis of $\mathrm{PLRT}_{n,\lambda}^{(j,j')}$ in Theorem 4.2. Since we take infimum over all $\beta_{0j} - \beta_{0j'} \in \mathcal{F}_\zeta$, we need a uniform probability bound for the power function. We say a random variable $X_n = O_P(t_n)$ uniformly over $\Delta g_N, \beta_{01}, \ldots, \beta_{0s} \in \mathcal{F}_\zeta$ if for any $\delta \in (0,1)$, there exist positive constants $C_\delta$ and $N_\delta$ such that for any $n \geq N_\delta$, $\inf_{\Delta g_N, \beta_{01}, \ldots, \beta_{0s} \in \mathcal{F}_\zeta} \mathbb{P}(X_n > C_\delta t_n) \geq 1 - \delta$. This definition also applies to the power analysis of other statistics in the paper. Decompose the statistic $L_{(1)} = T_1 + T_2 + T_3$ where

$$\text{(C.9)} \quad \begin{aligned} T_1 &:= \mathcal{L}_{n,\lambda}^{(1)}(\bar{g}_{N,\lambda}, \widehat{\beta}_{n,\lambda}^{(1)}) - \mathcal{L}_{n,\lambda}^{(1)}(\bar{g}_{N,\lambda}, \beta_{01}), \\ T_2 &:= \mathcal{L}_{n,\lambda}^{(1)}(\bar{g}_{N,\lambda}, \beta_{01}) - \mathcal{L}_{n,\lambda}^{(1)}(\bar{g}_{N,\lambda}, \beta_{02}), \\ T_3 &:= \mathcal{L}_{n,\lambda}^{(1)}(\bar{g}_{N,\lambda}, \beta_{02}) - \mathcal{L}_{n,\lambda}^{(1)}(\bar{g}_{N,\lambda}, \widehat{\beta}_{n,\lambda}^{(2)}). \end{aligned}$$

Following similar analysis above we have $(2u_n)^{-1/2}(-2nr_T(T_1 + T_3) - u_n) = O_P(1)$. The detailed analysis is also referred to the similar proof on $\mathrm{PLRT}_{n,\lambda}$ given in Supplementary S.3.3.

Define $r_N = \log N(Nh_1)^{-1/2} + \lambda_1^{1/2} + \lambda_2^{1/2}$. We can prove that
(C.10)
$$-2nT_2 = -n\|\Delta\beta\|^2 + O_P(n\lambda_2 + nr_N\|\Delta\beta\| + n^{1/2}\|\Delta\beta\| + n^{1/2}\|\Delta\beta\|^2),$$

where $\Delta\beta = \beta_{01} - \beta_{02}$. The proof is similar to the power analysis of $\mathrm{PLRT}_{n,\lambda}$ in Supplementary S.3.3.

Combining the above results on $T_1, T_2, T_3$, we have

$$(2u_n)^{-1/2}(-2nr_T \cdot \text{PLRT}_{n,\lambda}^{(j,j')} - u_n)$$
$$\geq n(2u_n)^{-1/2}\|\Delta\beta\|^2(1 + O_P(\lambda_2\|\Delta\beta\|^{-2} + r_N\|\Delta\beta\|^{-1}$$
$$+ n^{-1/2} + n^{-1/2}\|\Delta\beta\|^{-1})) + O_P(1).$$

This result is similar to the Lemma C.2 below on the power analysis of $\text{APLRT}_{N,\lambda}$. Following similar arguments in the proof of Theorem 4.6 in Appendix C.3, we can prove the result on $\text{PLRT}_{n,\lambda}^{(j,j')}$ in Theorem 4.2.

C.1.2. *Analysis of* $\text{PLRT}_{n,\lambda}^s$. Define

$$L_{(k)} = \mathcal{L}_{n,\lambda}^{(k)}(\bar{g}_{N,\lambda}, \widehat{\beta}_{n,\lambda}^{(k)}) - \mathcal{L}_{n,\lambda}^{(k)}(\bar{g}_{N,\lambda}, \widehat{\beta}_{n,\lambda}^{(k+1)}).$$

Therefore, $\text{PLRT}_{n,\lambda}^s = \sum_{k=1}^{s-1} L_{(k)}$. We denote $\beta_{j,k} = \widehat{\beta}_{n,\lambda}^{(j)} - \widehat{\beta}_{n,\lambda}^{(k)}, \widehat{f}_0^{(j)} = (\bar{g}_{N,\lambda}, \widehat{\beta}_{n,\lambda}^{(j)}), \beta_{n,\lambda}^{(j)} = \widehat{\beta}_{n,\lambda}^{(j)} - \beta_{0j}, f_{j,k} = (0, \beta_{j,k}), f_0^{(j)} = (\Delta g_{N,\lambda}, \beta_{n,\lambda}^{(j)})$. Following similar arguments as the proof of Lemma C.1 and (C.5), we have

$$\text{PLRT}_{n,\lambda}^s = W_1 + W_2 + W_3 + o_P(ns^{1/2}h^{-1/2}),$$

where

$$W_1 := \sum_{j=1}^{s-1} \frac{1}{2}\Big(\nabla^2 \mathcal{L}_{n,\lambda}^{(j)}(\widehat{f}_0^{(j)})[f_{j,j+1}, f_{j,j+1}] - \mathbb{E}\{\nabla^2 \mathcal{L}_{n,\lambda}^{(j)}(\widehat{f}_0^{(j)})\}[f_{j,j+1}, f_{j,j+1}]\Big),$$

$$W_2 := \sum_{j=1}^{s-1} \frac{1}{2}\mathbb{E}\{\nabla^2 \mathcal{L}_{n,\lambda}^{(j)}(f_{0j})\}[f_{j,j+1}, f_{j,j+1}],$$

$$W_3 := \sum_{j=1}^{s-1} \langle \nabla^2 \mathcal{L}_{n,\lambda}^{(j)}(f_{0j}) - \nabla^2 \mathcal{L}_{n,\lambda}^{(j+1)}(f_{0,j+1}), f_{j,j+1}\rangle.$$

According to Proposition A.2 and $\|f_{j,j+1}\| \leq \|f_0^{(j)}\| + \|f_0^{(j+1)}\|$, we have

$$|W_1| \leq \frac{1}{2}\sum_{j=1}^{s-1} n^{-1}\Big\|\sum_{i=1}^n f_{j,j+1}(U_i)K_{U_i}^\lambda - \mathbb{E}[f_{j,j+1}(U)K_U^\lambda]\Big\|\|f_{j,j+1}\|$$
$$\leq \sum_{j=1}^{s} n^{-1}\Big\|\sum_{i=1}^n f_0^{(j)}(U_i)K_{U_i}^\lambda - \mathbb{E}[f_0^{(j)}(U)K_U^\lambda]\Big\|\|f_0^{(j)}\|$$
$$+ \sum_{j=1}^{s-1} n^{-1}\Big\|\sum_{i=1}^n f_0^{(j)}(U_i)K_{U_i}^\lambda - \mathbb{E}[f_0^{(j)}(U)K_U^\lambda]\Big\|\|f_0^{(j+1)}\|.$$

Similar to (S.33) in the supplementary, we also have $|W_1| = o_P(sa_n r_n \log N)$. Since $a_n \log N = o(N^{-1/2})$, $n|W_1| = o_P(\sqrt{s/h})$. Applying Proposition A.2,

we have
$$W_2 = \sum_{j=1}^{s} \|f_{j,j+1}\|^2/2.$$

Similar to (C.3), we have
$$W_3 = -\sum_{j=1}^{s-1} \|\nabla \mathcal{L}_2^{(j)}(\beta_{0j}) - \nabla \mathcal{L}_2^{(j+1)}(\beta_{0(j+1)})\|^2 + o_P(s^{1/2}h^{-1/2}).$$

Therefore, it follows from the above results on $W_1, W_2$ and $W_3$ that

$$-2n \cdot \text{PLRT}_{n,\lambda}^s = \sum_{j=1}^{s-1} n\|\nabla \mathcal{L}_2^{(j)}(\beta_{0j}) - \nabla \mathcal{L}_2^{(j+1)}(\beta_{0(j+1)})\|^2 + o_P(s^{1/2}h^{-1/2})$$

$$(\text{C.11}) \qquad = \frac{1}{n}\sum_{j=1}^{s-1} \Big\| \sum_{i \in I_j \cup I_{j+1}} \epsilon_i Z_i T_{X_i}^\lambda \Big\|^2 + o_P(s^{1/2}h^{-1/2}).$$

Applying quadratic expansion to the first term in (C.11) yields
$$\frac{1}{n}\sum_{j=1}^{s-1} \sum_{i \in I_j \cup I_{j+1}} \epsilon_i^2 Z_i^2 T^\lambda(X_i, X_i) + \frac{1}{n}\sum_{1 \le i < j \le N} \widetilde{W}_{ij},$$

where $\widetilde{W}_{ij} = 2\epsilon_i\epsilon_j Z_i Z_j T^\lambda(X_i, X_j)$ if $i,j \in I_k \cup I_{k+1}$ for some $k = 1, \ldots, s-1$ and $W_{ij} = 0$ otherwise. Define another U-statistic
$$W_4(N) = \sum_{1 \le i < j \le N} \widetilde{W}_{ij}$$

with $\text{Var}(W_4(N)) \sim 6sn^2 h_2^{-1} \rho_T^2$. Similar to the analysis of $\text{APLRT}_{N,\lambda}$ in supplementary S.3.3, we can derive
$$\frac{1}{n}\sum_{j=1}^{s-1} \sum_{i \in I_j \cup I_{j+1}} \epsilon_i^2 Z_i^2 T^\lambda(X_i, X_i) = 2sh_2^{-1}\sigma_T^2 + O_P(h_2^{-1}).$$

This gives the following limiting distribution
$$(6sh_2^{-1}\sigma_T^4/\rho_T^2)^{-1/2}(-2n \cdot \text{PLRT}_{n,\lambda}^s - 2sh_2^{-1}\sigma_T^4/\rho_T^2) \rightsquigarrow N(0,1).$$

The power analysis of $\text{PLRT}_{n,\lambda}^s$ is a combination of those of $\text{APLRT}_{N,\lambda}$ and $\text{PLRT}_{n,\lambda}^{(j,j')}$. Define a set $\mathcal{I}_0 = \{j \mid \beta_{0j} = \beta_{0(j+1)}, 1 \le j \le s-1\}$ and $\mathcal{I}_1 = \{1, \ldots, s-1\}\setminus \mathcal{I}_0$. We decompose
$$\text{PLRT}_{n,\lambda}^s = \sum_{k \in \mathcal{I}_0} L_{(k)} + \sum_{k \in \mathcal{I}_1} L_{(k)}.$$

For the first term, due to the null limiting distribution and $|\mathcal{I}_0| < s$, we have
$$(6sh_2^{-1}\sigma_T^4/\rho_T^2)^{-1/2}\big(-2n\sum_{j\in\mathcal{I}_0} L^{(j)} - 2sh_2^{-1}\sigma_T^4/\rho_T^2\big) = O_P(1).$$
Define the difference $\Delta\beta_j = \beta_{0j} - \beta_{0(j+1)}$ on the set $\mathcal{I}_1$. As for the second term, we have the decomposition (by following (C.9))
$$\sum_{j\in\mathcal{I}_1} L_{(j)} = \sum_{j\in\mathcal{I}_1} \big(T_1^{(j)} + T_2^{(j)} + T_3^{(j)}\big),$$
where
(C.12)
$$\begin{aligned}T_1^{(j)} &:= \mathcal{L}_{n,\lambda}^{(j)}\big((\bar{g}_{N,\lambda}, \widehat{\beta}_{n,\lambda}^{(j)})\big) - \mathcal{L}_{n,\lambda}^{(j)}\big((\bar{g}_{N,\lambda}, \beta_{0j})\big),\\ T_2^{(j)} &:= \mathcal{L}_{n,\lambda}^{(j)}\big((\bar{g}_{N,\lambda}, \beta_{0j})\big) - \mathcal{L}_{n,\lambda}^{(j)}\big((\bar{g}_{N,\lambda}, \beta_{0(j+1)})\big),\\ T_3^{(j)} &:= \mathcal{L}_{n,\lambda}^{(j)}\big((\bar{g}_{N,\lambda}, \beta_{0(j+1)})\big) - \mathcal{L}_{n,\lambda}^{(j)}\big((\bar{g}_{N,\lambda}, \widehat{\beta}_{n,\lambda}^{(j+1)})\big).\end{aligned}$$

Using the same proof of Theorem 4.5, we have
$$(6sh_2^{-1}\sigma_T^4/\rho_T^2)^{-1/2}\big(-2n\sum_{j\in\mathcal{I}_1}(T_1^{(j)} + T_3^{(j)}) - 2sh_2^{-1}\sigma_T^4/\rho_T^2\big) = O_P(1).$$
By (C.10), defining $\Delta^* = \sum_{j\in\mathcal{I}_1} \|\Delta\beta_j\|^2$, we also have
$$-2n\sum_{j\in\mathcal{I}_1} T_2^{(j)} = -n\Delta^* + O_P\Big(n\lambda_2 + nh_2^{-1}r_n(\Delta^*)^{1/2} \\ + n^{1/2}(\Delta^*)^{1/2} + n^{1/2}r_N(\Delta^*)^{1/2} + n^{1/2}\Delta^*\Big),$$
Define $u'_N = sh_2^{-1}\sigma_T^4/\rho_T^2$. The results on $T_1^{(j)}, T_2^{(j)}, T_3^{(j)}$ imply that
$$\begin{aligned}(6u'_N)^{-1/2}&(-2n\cdot \mathrm{PLRT}_{n,\lambda}^s - 2u'_N) \\ &\geq n(6u'_N)^{-1/2}\Delta^*\big[1 + O_P(\lambda_2(\Delta^*)^{-1} + r_N(\Delta^*)^{-1/2} \\ &\quad + n^{-1/2} + n^{-1/2}(\Delta^*)^{-1/2})\big] + O_P(1) \\ &\geq n(6u'_N)^{-1/2}\Delta^*\big[1 + O_P(\lambda_2(\Delta^*)^{-1} + r_N(\Delta^*)^{-1/2} \\ &\quad + n^{-1/2} + n^{-1/2}(\Delta^*)^{-1/2})\big] + O_P(1).\end{aligned}$$
We again get a similar result as Lemma C.2. Following a similar proof as that of Theorem 4.6 in Appendix C.3, we can also prove the power result on $\mathrm{PLRT}_{n,\lambda}^s$ in Theorem 4.2.

**C.2. Proof of Theorem 4.5.** Similar to Theorem 3.4, we consider the following decomposition
(C.13)
$$-2N\cdot \mathrm{APLRT}_{N,\lambda} - N\|W_{1,\lambda}g_0\|^2 - h_1^{-1}\sigma_P^2 = \underbrace{N^{-1}W(N)}_{\text{Leading term}} + \underbrace{M_n}_{\text{Remainder term}},$$
where $W(N) = \sum_{i\neq j} \epsilon_i \epsilon_j P^\lambda(X_i, X_j)$. The following lemma derives the convergence rate of $M_n$, whose proof is deferred to supplementary S.3.1.

**Lemma C.1.** Under the same conditions of Theorem 4.5, we have $M_n = o_P(h_1^{-1/2})$.

We next show the asymptotic normality of $N^{-1}W(N)$ by following a similar procedure as $W_3(2n)$ in (C.8). We can again check the moment conditions of $W(N)$ required by the Proposition 3.2 in de Jong (1987) and show that

$$\text{(C.14)} \qquad \frac{1}{\sqrt{2/h_1} \cdot N\rho_P} W(N) \rightsquigarrow N(0,1).$$

Combining (C.13) and (C.14) with Lemma C.1, we have

$$(2u_N)^{-1/2}\left(-2Nr_P \cdot \text{APLRT}_{N,\lambda} - Nr_P\|W_{1,\lambda}g_0\|^2 - u_N\right) \rightsquigarrow N(0,1),$$

where $u_N = h^{-1}\sigma_P^4/\rho_P^2$, $r_P = \sigma_P^2/\rho_P^2$. Since $nh_1\lambda_1 = o(1)$, by (S.22), the bias term $(2u_N)^{-1/2}Nr_P\|W_{1,\lambda}g_0\|^2 = o(1)$, we complete the proof of (4.9).

**C.3. Proof of Theorem 4.6.** We first need the following lemma on the decomposition of testing statistic.

**Lemma C.2.** Suppose the same conditions of Theorem 4.6 hold. We have

$$\left|\frac{-2Nr_P \cdot \text{APLRT}_{N,\lambda} - u_N}{(2u_N)^{1/2}}\right| \geq \frac{N\|\Delta g_N\|^2}{(2u_N)^{1/2}}$$
$$\times \left[1 + O_P\left(\frac{\lambda_1}{\|\Delta g_N\|^2} + \frac{\sqrt{\lambda_2} + N^{-1/2}}{\|\Delta g_N\|} + \frac{1}{N^{1/2}}\right)\right] + O_P(1),$$

where $O_P(\cdot)$ holds uniformly for $f_N \in \mathcal{F}_\zeta$.

We defer the proof of Lemma C.2 to supplementary S.3.2.

For any $\delta \in (0,1)$, it suffices to show $(2u_N)^{-1/2}(-2Nr_P \cdot \text{APLRT}_{N,\lambda} - u_N) \geq z_\alpha$ with probability $1 - \delta$ if there exists constants $N_\delta$ and $c_N$ such that for any $N > N_\delta$, $N^{-1/2}\|\Delta g_N\|^{-1} \leq 1/c_N^2$, $\lambda_1\|\Delta g_N\|^{-2} \leq 1/c_N^2$, $\sqrt{\lambda_2}\|\Delta g_N\|^{-1} \leq 1/c_N^2$ and $\|\Delta g_N\|^2 \geq c_N^2(Nh^{1/2})^{-1}$ based on Lemma C.2. The above four inequalities can be derived from the fact that

$$\|\Delta g_N\| \geq c_N\sqrt{\lambda_1 + \lambda_2 + (Nh^{1/2})^{-1}}.$$

We have for any $N > N_\delta$,

$$\inf_{\{\beta_{0j}\}_{j=1}^s \subseteq \mathcal{F}_\zeta} \inf_{\substack{\Delta g_N \in \mathcal{F}_\zeta \\ \|\Delta g_N\|^2 \geq c_N\eta_N^2}} \mathbb{P}_{H_{1N}}\left(\frac{|-2Nr_P \cdot \text{APLRT}_{N,\lambda} - u_N|}{\sqrt{2u_N}} > z_\alpha\right) \leq 1 - \delta,$$

which completes the proof of this theorem.

*Supplementary material to*

# NONPARAMETRIC HETEROGENOUS INFERENCE FOR MASSIVE DATA

By Junwei Lu¶, Guang Cheng‖,∗∗ and Han Liu¶

*Princeton University¶ and Purdue University‖*

The supplementary material is organized as follows:

- In Section S.1, we prove Proposition A.1, Proposition A.2, Lemma B.2 and some results on the covering number of RKHS.
- In Section S.2, we prove Lemma B.1.
- In Section S.3, we prove Lemma C.1 and Lemma C.2 for APLRT$_{N,\lambda}$. We also prove the remaining part of Theorem 4.1 and Theorem 4.2 for the functional test.
- In Section S.4, we list several technical results useful for our proofs.

## APPENDIX S.1: AUXILIARY LEMMAS FOR RKHS

In this section, we present the proofs of results on the representations of operators in $\mathcal{H}^2$ including Proposition A.1, Proposition A.2 and Lemma B.2. We also state some auxiliary lemmas on the covering number of RKHS. Define another inner product in $\mathcal{H}^2$ as

$$\langle (g_1, \beta_1), (g_2, \beta_2) \rangle_{\mathcal{H}^2} = (\lambda_1 \wedge \lambda_2)^{-1}(\lambda_1 \langle g_1, g_2 \rangle_{\mathcal{H}_1} + \lambda_2 \langle \beta_1, \beta_2 \rangle_{\mathcal{H}_2})$$

for any $(g_1, \beta_1), (g_2, \beta_2) \in \mathcal{H}^2$. The induced norm is denoted as $\|\cdot\|_{\mathcal{H}^2}$.

**S.1.1. Proof of Proposition A.1.** As the existence of $K_u^\lambda$ is proven in Section A.1, we now turn to solve its Fourier expansion. We denote $K_u^\lambda = (\bar{P}_u^\lambda, \bar{T}_u^\lambda)$. For any $f = (g, \beta) \in \mathcal{H}^2$, due to the reproducing property of $K_u^\lambda$, we have $\langle K_u^\lambda, f \rangle = g(x) + z\beta(x)$. On the other hand, by the definition

∗∗Associate Professor, Department of Statistics, Purdue University, West Lafayette, IN 47906. E-mail: chengg@purdue.edu. Tel: +1 (765) 496-9549. Fax: +1 (765) 494-0558. Research Sponsored by NSF CAREER Award DMS-1151692, DMS-1418042, Simons Fellowship in Mathematics, Office of Naval Research (ONR N00014-15-1-2331) and a grant from Indiana Clinical and Translational Sciences Institute. Guang Cheng was on sabbatical at Princeton while this work was finalized; he would like to thank the Princeton ORFE department for its hospitality and support.

of $\langle \cdot, \cdot \rangle$, it can be derived that

$$
\begin{aligned}
g(x) + z\beta(x) &= \langle K_u^\lambda, f \rangle \\
&= \mathbb{E}[(\bar{P}_u^\lambda(X) + Z\bar{T}_u^\lambda(X))(g(X) + Z\beta(X))] \\
&\quad + \lambda_1 \langle \bar{P}_u^\lambda, g \rangle_\mathcal{H} + \lambda_2 \langle \bar{T}_u^\lambda, \beta \rangle_\mathcal{H} \\
&= \mathbb{E}[\bar{P}_u^\lambda(X) g(X)] + \lambda_1 \langle \bar{P}_u^\lambda, g \rangle_\mathcal{H} \\
&\quad + \sigma_z^2 \mathbb{E}[\bar{T}_u^\lambda(X) \beta(X)] + \lambda_2 \langle \bar{T}_u^\lambda, \beta \rangle_\mathcal{H}.
\end{aligned}
\tag{S.1}
$$

Since (S.1) is true for any $g \in \mathcal{H}, h \in \mathcal{H}$, we have

$$
\begin{cases}
\mathbb{E}[\bar{P}_u^\lambda(X) g(X)] + \lambda_1 \langle \bar{P}_u^\lambda, g \rangle_\mathcal{H} = g(x), \\
\sigma_z^2 \mathbb{E}[\bar{T}_u^\lambda(X) \beta(X)] + \lambda_2 \langle \bar{T}_u^\lambda, \beta \rangle_\mathcal{H} = z\beta(x).
\end{cases}
\tag{S.2}
$$

Suppose we have the Fourier expansions $\bar{P}_u^\lambda = \sum_{s=1}^\infty p_s(u) \phi_s$ and $\bar{T}_u^\lambda = \sum_{s=1}^\infty q_s(u) \phi_s$. We set $g, \beta$ in (S.2) as $\phi_s$ and $\phi_s$ and plug the Fourier expansions of $\bar{P}_u^\lambda$ and $\bar{T}_u^\lambda$ into (S.2). We thus have

$$
p_s(u) = \frac{\phi_s(x)}{1 + \lambda_1/\mu_s} \quad \text{and} \quad q_s(u) = z \frac{\phi_s(x)}{\sigma_z^2 + \lambda_2/\mu_s}.
$$

This shows that $\bar{P}_u^\lambda = P_x^\lambda$ and $\bar{T}_u^\lambda = zT_x^\lambda$, which completes the computation of $K_u^\lambda$'s Fourier expansion. We can easily check that the obtained $K_u^\lambda = (P_x^\lambda, zT_x^\lambda)$ has the reproducing property for any $f \in \mathcal{H}^2$ through (S.1). From (S.1), we also have $\langle P_x^\lambda, g \rangle = g(x)$ and $\langle T_x^\lambda, \beta \rangle = \beta(x)$ for any $g, \beta \in \mathcal{H}$.

We now derive the Fourier expansion of $W_\lambda = (W_{\lambda,1}, W_{\lambda,2})$. Since the bilinear operator $\langle \cdot, \cdot \rangle_{\mathcal{H}^2}$ is bounded in the norm $\|\cdot\|$, according to the representation theorem of bilinear operator in Hilbert space, the existence of such a positive definite self-adjoint operator $W_\lambda : \mathcal{H}^2 \to \mathcal{H}^2$ is guaranteed. On one hand, by the definition of $W_\lambda$ we have

$$
\langle W_\lambda(\phi_s, \phi_t), (\phi_u, \phi_v) \rangle = (\lambda_1/\mu_s)\delta_{s,u} + (\lambda_2/\mu_t)\delta_{t,v}.
\tag{S.3}
$$

On the other hand, by the definition of $\langle \cdot, \cdot \rangle$, if we suppose the Fourier expansions are $W_{\lambda,1}\phi_s = \sum_{k=1}^\infty a_{s,k}\phi_k$ and $W_{\lambda,2}\phi_t = \sum_{k=1}^\infty b_{t,k}\phi_k$, we also have

$$
\langle W_\lambda(\phi_s, \phi_t), (\phi_u, \phi_v) \rangle = a_{s,u}(1 + \lambda_1/\mu_u) + b_{t,v}(\sigma_z^2 + \lambda_2/\mu_v).
\tag{S.4}
$$

Combing (S.3) and (S.4), we can solve

$$
W_{\lambda,1}\phi_s = \frac{\lambda_1/\mu_s}{1 + \lambda_1/\mu_s}\phi_s \text{ and } W_{\lambda,2}\phi_s = \frac{\lambda_2/\mu_s}{\sigma_z^2 + \lambda_2/\mu_s}\phi_s.
$$

By linearity, we derive the expansions of $W_{\lambda,1}g$ and $W_{\lambda,2}h$ stated in the Lemma. It is also easy to verify that $\langle (W_{\lambda,1}g_1, W_{\lambda,2}h_1), (g_2, \beta_2) \rangle = \lambda_1 \langle g_1, g_2 \rangle_{\mathcal{H}_1} + \lambda_2 \langle \beta_1, \beta_2 \rangle_{\mathcal{H}_2}$ for all $g_1, g_2, h_1, h_2 \in \mathcal{H}$, which completes the proof of the propo-sition.

**S.1.2. Proof of Proposition A.2.** We first prove the existence of the representations. By the Riesz representation theorem for both linear and bilinear operators, it suffices to check that $D\mathcal{L}_{n,\lambda}(f)$ and $D^2\mathcal{L}_{n,\lambda}(f)$ are bounded. By definition, for any $f_1 \in \mathcal{H}^2$ we have

$$|D\mathcal{L}_{n,\lambda}(f)[f_1]| \leq \frac{1}{n}\sum_{i=1}^n \left|(Y_i - f(U_i))f_1(U_i)\right| + \|f\|\|f_1\|$$

$$\leq \frac{1}{n}\sum_{i=1}^n \left(\left|(Y_i - f(U_i))\right|\|K_{U_i}^\lambda\| + \|f\|\right)\|f_1\|$$

Similarly, for any $f_1, f_2 \in \mathcal{H}^2$,

$$|D^2\mathcal{L}_{n,\lambda}(f)[f_1, f_2]| \leq \frac{1}{n}\sum_{i=1}^n |f_1(U_i)f_2(U_i)| + \|f_1\|\|f_2\|$$

$$\leq \frac{1}{n}\sum_{i=1}^n \left(\|K_{U_i}^\lambda\|^2 + 1\right)\|f_1\|\|f_2\|.$$

The representation of specific operators are given in Proposition A.2.

According to the optimality of $\widehat{f}_{n,\lambda}^{(j)}$ and the corresponding KKT condition (Luenberger, 1969), $\nabla\mathcal{L}_{n,\lambda}(\widehat{f}_{n,\lambda}^{(j)}) = 0$. The definition of $D\mathcal{L}_{n,\lambda}(f)[\widetilde{f}]$ implies that for any $\widetilde{f} \in \mathcal{H}$,

$$D\mathcal{L}_{n,\lambda}(f_{0j})[\widetilde{f}] = -\frac{1}{n}\sum_{i=1}^n (Y_i - f_{0j}(U_i))\widetilde{f}(X_i) + \langle W_\lambda f_{0j}, \widetilde{f}\rangle$$

$$= -\frac{1}{n}\sum_{i=1}^n \epsilon_i\langle K_{U_i}^\lambda, \widetilde{f}\rangle + \langle W_\lambda f_{0j}, \widetilde{f}\rangle$$

$$= \Big\langle -\frac{1}{n}\sum_{i=1}^n \epsilon_i K_{U_i}^\lambda + W_\lambda f_{0j}, \widetilde{f}\Big\rangle.$$

By a similar method, for any $f_1 = (g_1, \beta_1), f_2 = (g_2, \beta_2) \in \mathcal{H}^2$,

$$\mathbb{E}D^2\mathcal{L}_{n,\lambda}(f_{0j})[f_1, f_2] = \frac{1}{n}\sum_{i=1}^n \mathbb{E}[f_1(U_i)f_2(U_i)] + \lambda_1\langle g_1, g_2\rangle_{\mathcal{H}_1} + \lambda_2\langle \beta_1, \beta_2\rangle_{\mathcal{H}_2}$$

$$= \mathbb{E}[f_1(U)f_2(U)] + \lambda_1\langle g_1, g_2\rangle_{\mathcal{H}_1} + \lambda_2\langle \beta_1, \beta_2\rangle_{\mathcal{H}_2} = \langle f_1, f_2\rangle.$$

The assertions of the proposition are therefore proved.

**S.1.3. Proof of Lemma B.2.** Instead of proving Lemma B.2, we directly prove the following lemma with stronger results.

**Lemma S.1.** Under Assumption 3.2, we have $\|K_u^\lambda\|^2 \leq c_K h^{-1}, \|P_x^\lambda\|^2 \leq c_K h_1^{-1}$ and $\|T_x^\lambda\|^2 \leq c_K h_2^{-1}$. Moreover, the norm $\|\cdot\|_{\sup}$ can be bounded by $\|\cdot\|$ such that for all $f = (g, \beta) \in \mathcal{H}^2$,

$$\|f\|_{\sup} \leq \min\Big\{c_K h^{-1/2}\|f\|,\, c_\mathcal{H}\|g\|_{\mathcal{H}_1}, c_\mathcal{H}\|\beta\|_{\mathcal{H}_2}\Big\}.$$

where $c_K$ is the universal upper bound of eigen-functions defined in (3.1) and $c_\mathcal{H} = \sup_{x \in \mathcal{X}} \sqrt{|K(x,x)|}$.

PROOF. We first bound the norm of $\|K_x^\lambda\|$. Using the reproducing property of $K^\lambda$, we have

$$\|K_u^\lambda\|^2 = P_x^\lambda(x) + zT_x^\lambda(x)$$
$$= \sum_{j \in \mathbb{Z}} \frac{|\phi_j(x)|^2}{1 + \lambda_1/\mu_{1,j}} + z \sum_{j \in \mathbb{Z}} \frac{|\phi_j(x)|^2}{\sigma_z^2 + \lambda_2/\mu_{2,j}}$$
$$\leq \sum_{j \in \mathbb{Z}} \frac{c_\phi^2}{1 + \lambda_1/\mu_{1,j}} + \sum_{j \in \mathbb{Z}} \frac{c_z c_\phi^2}{\sigma_z^2 + \lambda_2/\mu_{2,j}}$$
(S.5) $$= (1 + c_z)c_\phi^2 \gamma(\lambda) \leq (1 + c_z)c_\phi^2 c_h h^{-1}.$$

Use the reproducing property again:

$$|f(x,z)| = |\langle K_u^\lambda, f \rangle| \leq \|K_u^\lambda\| \cdot \|f\| \leq c_K h^{-1/2}\|f\|.$$

Similarly, we have

$$|f(x,z)| \leq |\langle K_{1,x}, g\rangle|_{\mathcal{H}_1} + |z| \cdot |\langle K_{2,x}, \beta\rangle|_{\mathcal{H}_2} \leq (1 + c_z)c_b\|f\|_{\mathcal{H}^2}.$$

The desired result is proved. Similarly, we can also show that $\|P_x^\lambda\|^2 \lesssim h_1^{-1}$ and $\|T_x^\lambda\|^2 \lesssim h_2^{-1}$. □

**S.1.4. Covering Number of RKHS.** In this section, we summarize several useful results on the covering number of RKHS with certain eigenvalue decaying rate, and one related result on entropy integral.

**Lemma S.2.** Suppose Assumption 3.2 holds, and let $\mathbb{B}_\mathcal{H}(1) = \{g \in \mathcal{H} : \|g\|_\mathcal{H} \leq 1\}$. For all $\epsilon \in (0, 1)$, we have

1. Finite rank (Carl and Triebel, 1980): if $\mu_k = 0$ for all $k > r$, then $\log N(\epsilon; \mathbb{B}_\mathcal{H}(1), \|\cdot\|_\mathcal{H}) = O(r \log(1/\epsilon))$;
2. Polynomially decaying (Carl and Stephani, 1990): if $\mu_k = O(k^{-2m})$, for some $m > 0$, then $\log N(\epsilon; \mathbb{B}_\mathcal{H}(1), \|\cdot\|_\mathcal{H}) = O\big((1/\epsilon)^{1/m}\big)$;

3. Exponentially decaying (Williamson et al., 2001; Guo et al., 2002): if $\mu_k = O(e^{-\alpha k^p})$, then $\log N(\epsilon; \mathbb{B}_{\mathcal{H}}(1), \|\cdot\|_{\mathcal{H}}) = O\big((\log(1/\epsilon))^{\frac{p+1}{p}}\big)$.

Next lemma gives a lower bound of the entropy integral $\xi_\lambda$ defined in (3.8).

**Lemma S.3.** Define $\mathcal{H}' = \{f = (g, \beta) \in \mathcal{H}^2 : \|f\|_{\sup} \leq 1, \|g\|_{\mathcal{H}_1}^2 \leq c_K^{-2} h \lambda_1^{-1}, \|\beta\|_{\mathcal{H}_2}^2 \leq c_K^{-2} h \lambda_2^{-1}\}$. We have

$$\xi_\lambda \geq \int_0^1 \sqrt{\log(1 + N(\mathcal{H}'; \|\cdot\|_{\sup}, \epsilon))} d\epsilon.$$

PROOF. We rescale the upper bounds in the definition of $\mathcal{H}'$ and obtain
$$\mathcal{H}' = (c_K^{-2} h \lambda^{-1})^{1/2} \{f \in \mathcal{H}^2 : \|f\|_{\sup} \leq (c_K^{-2} h \lambda^{-1})^{-1/2}, \|g\|_{\mathcal{H}_1} \leq 1, \|\beta\|_{\mathcal{H}_2} \leq 1\}$$
$$\subseteq (c_K^{-2} h \lambda^{-1})^{1/2} \mathcal{H}_0.$$

Therefore, the entropy integral of $\mathcal{H}'$ can be bounded by

$$\int_0^1 \sqrt{\log(1 + N(\mathcal{H}'; \|\cdot\|_{\sup}, \epsilon))} d\epsilon$$
$$\leq \int_0^1 \sqrt{\log(1 + N(\mathcal{H}_0; \|\cdot\|_{\sup}, (c_K^{-2} h \lambda^{-1})^{-1/2} \epsilon))} d\epsilon$$
$$= (c_K^{-2} h \lambda^{-1})^{1/2} \int_0^{(c_K^{-2} h \lambda^{-1})^{-1/2}} \sqrt{\log(1 + N(\mathcal{H}_0; \|\cdot\|_{\sup}, \epsilon))} d\epsilon$$
$$= (c_K^{-2} h \lambda^{-1})^{1/2} \omega((c_K^{-2} h \lambda^{-1})^{-1/2}) = \xi_\lambda.$$

□

## APPENDIX S.2: AUXILIARY RESULTS ON REMAINDER TERM

In this section, we present the proof of Lemma B.1 on the rate of the remainder term $R_N^{(j)} = (G_N, H_n^{(j)})$ defined in (B.1).

**S.2.1. Proof of Lemma B.1.** Denote $f_{n,\lambda}^{(j)} = \widehat{f}_{n,\lambda}^{(j)} - f_{0j} =: (g_{n,\lambda}, \beta_{n,\lambda}^{(j)})$ for all $j = 1, \ldots, s$. For notational simplicity, we also denote $f_{n,\lambda}^{(j)}$ as $f_{n,\lambda}$. Similar rule applies to other estimators when no confusion arises. We first study the tail probability of $\lambda \|f_{n,\lambda}\|_{\mathcal{H}}$. By Assumption 3.1, $\epsilon^2$ is a sub-exponential random variable with $\mathbb{E}[\epsilon^2] = \sigma^2$. According to Bernstein inequality for sub-exponential random variables, there exists $c' > \sigma^2$ such that

(S.6) $$\mathbb{P}\Big(\frac{1}{n} \sum_{i=1}^n \epsilon_i^2 > c'\Big) \leq 2 \exp(-c' n).$$

By the definition of $\widehat{f}_{n,\lambda}$, we have

$$(\lambda_1 \wedge \lambda_2)\|\widehat{f}_{n,\lambda}\|^2_{\mathcal{H}^2} \leq \frac{1}{n}\sum_{i=1}^n (\widehat{f}_{n,\lambda}(U_i) - Y_i)^2 + \lambda_1\|\widehat{g}_{n,\lambda}\|^2_{\mathcal{H}_1} + \lambda_2\|\widehat{\beta}_{n,\lambda}\|^2_{\mathcal{H}_2}$$

$$\leq \frac{1}{n}\sum_{i=1}^n \epsilon_i^2 + \lambda_1\|g_0\|^2_{\mathcal{H}_1} + \lambda_2\|\beta_0\|^2_{\mathcal{H}_2}.$$

Taking expectation on both sides yields

(S.7) $\qquad (\lambda_1 \wedge \lambda_2)\mathbb{E}\|\widehat{f}_{n,\lambda}\|^2_{\mathcal{H}^2} \leq \sigma^2 + \lambda_1\|g_0\|^2_{\mathcal{H}_1} + \lambda_2\|\beta_0\|^2_{\mathcal{H}_2} = O(1).$

Since $\lambda_1 \wedge \lambda_2 = o(1)$, we can choose some constant $c_0$ satisfying $c_0 \geq 2c' + 4\lambda_1\|g_0\|^2_{\mathcal{H}} + 4\lambda_2\|\beta_0\|^2_{\mathcal{H}}$ such that

$$\mathbb{P}((\lambda_1 \wedge \lambda_2)\|f_{n,\lambda}\|^2_{\mathcal{H}^2} \geq c_0)$$
$$\leq \mathbb{P}(\lambda_1\|g_{n,\lambda}\|^2_{\mathcal{H}_1} + \lambda_2\|\beta_{n,\lambda}\|^2_{\mathcal{H}_2} > c_0)$$
$$\leq \mathbb{P}(\lambda_1\|\widehat{g}_{n,\lambda}\|^2_{\mathcal{H}_1} + \lambda_2\|\widehat{\beta}_{n,\lambda}\|^2_{\mathcal{H}_2} > c' + \lambda_1\|g_0\|^2_{\mathcal{H}_1} + \lambda_2\|\beta_0\|^2_{\mathcal{H}_2})$$
(S.8) $\qquad \leq 2\exp(-c'n).$

Now, we define the remainder term on the $j$-th sub-population

$$R_n^{(j)} = \widehat{f}_{n,\lambda}^{(j)} - f_{0j} - \nabla\mathcal{L}_{n,\lambda}^{(j)}(f_{0j}) = \widehat{f}_{n,\lambda}^{(j)} - f_{0j}^u - \frac{1}{n}\sum_{i \in I_j}\epsilon_i K_{U_i}^\lambda.$$

We further denote $R_n^{(j)} = \left(G_n^{(j)}, H_n^{(j)}\right)$. According to the construction of $\bar{g}_{N,\lambda}$, we have

$$G_N := \bar{g}_{N,\lambda} - g_0^u - \frac{1}{N}\sum_{i=1}^N \epsilon_i P_{X_i}^\lambda$$

$$= \frac{1}{s}\sum_{j=1}^s \widehat{g}_{n,\lambda}^{(j)} - g_0^u - \frac{1}{s}\sum_{j=1}^s \sum_{i \in I_j}\epsilon_i P_{X_i}^\lambda$$

(S.9) $\qquad = \frac{1}{s}\sum_{j=1}^s \left(\widehat{g}_{n,\lambda}^{(j)} - g_0^u - \frac{1}{n}\sum_{i \in I_j}\epsilon_i P_{X_i}^\lambda\right) = \frac{1}{s}\sum_{j=1}^s G_n^{(j)}.$

Let $r_n = (nh)^{-1/2} + \lambda_1^{1/2} + \lambda_2^{1/2}$ and define the event $\mathcal{E}_j = \{f_{n,\lambda}^{(j)} \in \mathcal{H}, \|f_{n,\lambda}^{(j)}\| \leq c_1 r_n\}$ and $\mathcal{E} = \bigcap_{j=1}^s \mathcal{E}_j$.

We first study the remainder term in the $j$th sub-population. For simplicity,

we omit the superscript $(j)$ for $R_n^{(j)}$ as well.

$$\begin{aligned}
R_n &= f_{n,\lambda} - \nabla \mathcal{L}_{n,\lambda}(f_{0j}) \\
&= \mathbb{E}\nabla^2 \mathcal{L}_{n,\lambda}(f_{0j})[f_{n,\lambda}] - \nabla \mathcal{L}_{n,\lambda}(f_{0j}) \\
&= \mathbb{E}\nabla \mathcal{L}_{n,\lambda}(f_{n,\lambda} + f_{0j}) - \mathbb{E}\nabla \mathcal{L}_{n,\lambda}(f_{0j}) - \nabla \mathcal{L}_{n,\lambda}(f_{0j}) \\
&= (\nabla \mathcal{L}_{n,\lambda}(f_{0j}) - \mathbb{E}\nabla \mathcal{L}_{n,\lambda}(f_{0j})) \\
&\quad - (\nabla \mathcal{L}_{n,\lambda}(f_{n,\lambda} + f_{0j}) - \mathbb{E}\nabla \mathcal{L}_{n,\lambda}(f_{n,\lambda} + f_{0j}))
\end{aligned}$$

$$\text{(S.10)} \qquad = -\frac{1}{n} \sum_{i=1}^{n} \left( f_{n,\lambda}(U_i) K_{U_i}^\lambda - \mathbb{E}[f_{n,\lambda}(U) K_U^\lambda] \right).$$

The first equality is due to the definition of remainder term, the second and forth equalities are by Proposition A.2 and the third equality is by the taylor expansion of least square.

In order to bound the rate of $R_n$, we turn to study the empirical process in (S.10). We separate $R_n$ into two parts:

$$R_{n1} = \frac{1}{n} \sum_{i=1}^{n} \left( f_{n,\lambda}(U_i) \mathbb{1}_{\mathcal{E}_j} K_{U_i}^\lambda - \mathbb{E}[f_{n,\lambda}(U) \mathbb{1}_{\mathcal{E}_j} K_U^\lambda] \right), R_{n2} = R_n - R_{n1}.$$

We first bound $R_{n1}$. Let $d_n = c_1 c_K r_n h^{-1/2}$ and $\widetilde{f}_{n,\lambda} = (\widetilde{g}_{n,\lambda}, \widetilde{\beta}_{n,\lambda}) = d_n^{-1} f_{n,\lambda} \mathbb{1}_{\mathcal{E}_j}$. On the event $\mathcal{E} = \{\|f_{n,\lambda}^{(j)}\| \leq c_1 r_n \text{ for all } j\}$, it follows by Lemma S.1 that $\|\widetilde{f}_{n,\lambda}\|_{\sup} \leq 1$, $\|\widetilde{g}_{n,\lambda}\|_{\mathcal{H}_1}^2 \leq d_n^{-2} \lambda_1^{-1} r_n^2 \leq c_1^{-1} c_K^{-1} h \lambda_1^{-1}$ and $\|\widetilde{\beta}_{n,\lambda}\|_{\mathcal{H}_2}^2 \leq c_1^{-1} c_K^{-1} h \lambda_2^{-1}$.

Define an operator

$$\text{(S.11)} \qquad \Psi(f, U) = c_K^{-1} h^{1/2} f K_U^\lambda,$$

and we consider the following empirical process:

$$\mathbb{G}_n(\widetilde{f}_{n,\lambda}) = \frac{1}{\sqrt{n}} \sum_{i=1}^{n} \left( \Psi(\widetilde{f}_{n,\lambda}, U) - \mathbb{E}[\Psi(\widetilde{f}_{n,\lambda}, U)] \right).$$

We now apply Lemma S.6 for $\mathbb{G}_n(\widetilde{f}_{n,\lambda})$ with $d(\cdot) = \|\cdot\|_{\sup}$ and $\mathcal{F} = \mathcal{H}' = \{f = (g, \beta) \in \mathcal{H}^2 : \|f\|_{\sup} \leq 1, \|g\|_{\mathcal{H}_1}^2 \leq c_K^{-2} h \lambda_1^{-1}, \|\beta\|_{\mathcal{H}_2}^2 \leq c_K^{-2} h \lambda_2^{-1}\}$. To check the Lipschitz property, we have for any $f_1, f_2 \in \mathcal{H}^2$,

$$\text{(S.12)} \quad \|\Psi(f_1, U) - \Psi(f_2, U)\| \leq c_K^{-1} h^{1/2} \|K_U^\lambda\| \|f_1 - f_2\|_{\sup} \leq \|f_1 - f_2\|_{\sup},$$

where the last inequality follows (S.5). According to (S.46), we have for all

$j = 1, \ldots, s$,

$$\mathbb{P}\left(\frac{n\|R_{n1}^{(j)}\|}{d_n c_K h^{-1/2}(\sqrt{n}\xi_\lambda + 1)} \geq t\right)$$

$$= \mathbb{P}\left(\frac{\|\sum_{i \in I_j}[f_{n,\lambda}^{(j)}(U_i)\mathbb{1}_{\mathcal{E}_j}K_{U_i}^\lambda - \mathbb{E}(f_{n,\lambda}^{(j)}(U)\mathbb{1}_{\mathcal{E}_j}K_U^\lambda)]\|}{d_n c_K h^{-1/2}(\sqrt{n}\xi_\lambda + 1)} \geq t\right)$$

$$= \mathbb{P}\left(\frac{\sqrt{n}\|\mathbb{G}_n(\widetilde{f}_{n,\lambda}^{(j)})\|}{\sqrt{n}\xi_\lambda + 1} > t\right) \leq \mathbb{P}\left(\sup_{f \in \mathcal{H}} \frac{\sqrt{n}\|\mathbb{G}_n(f)\|}{\sqrt{n}\,\omega_d(1;\mathcal{H}) + 1} \geq t\right)$$

(S.13) $\quad \leq 2\left(\log(\sqrt{n}\xi_\lambda) + 2\right)\exp(-t^2/4)$,

where $\xi_\lambda$ is defined in (3.8). The first inequality follows from Lemma S.3 that $\xi_\lambda \geq \omega_d(1;\mathcal{H}')$, and the last inequality is due to (S.47) in Lemma S.6.

Next, we study the second term

$$R_{n2} = \frac{1}{n}\sum_{i=1}^n \left(f_{n,\lambda}(U_i)\mathbb{1}_{\mathcal{E}_j^c}K_{U_i}^\lambda - \mathbb{E}[f_{n,\lambda}(U)\mathbb{1}_{\mathcal{E}_j^c}K_U^\lambda]\right).$$

Recall that $f_{n,\lambda}^{(j)} = \widehat{f}_{n,\lambda}^{(j)} - f_{0j}$, and $f_{n,\lambda}$ is its simplified version. To study $R_{n2}$, we need the following lemma on the tail probability of $\|f_{n,\lambda}\|$.

**Lemma S.4.** For a sufficiently large $n$, there exist constants $C, c_0$ such that

$$\mathbb{P}\left(\|f_{n,\lambda}\| > C\left(\frac{\log N}{\sqrt{nh}} + \lambda_1^{1/2} + \lambda_2^{1/2}\right)\right) \leq 4n\exp(-c_0 \log^2 N).$$

We defer the proof of the lemma to Section S.2.2. Using Markov inequality, we have for all $j$,

$$\mathbb{P}\left(n\|R_{n2}^{(j)}\| \geq t\right) = \mathbb{P}\left(\sum_{i=1}^n \left(f_{n,\lambda}(U_i)\mathbb{1}_{\mathcal{E}_j^c}K_{U_i}^\lambda - \mathbb{E}[f_{n,\lambda}(U)\mathbb{1}_{\mathcal{E}_j^c}K_U^\lambda]\right) > t\right)$$

$$\leq \frac{2n}{t}\mathbb{E}\left\|f_{n,\lambda}(U)\mathbb{1}_{\mathcal{E}_j^c}K_U^\lambda\right\|$$

$$\leq 2t^{-1}h^{-1/2}n\sqrt{\mathbb{E}\|f_{n,\lambda}\|_{\mathcal{H}^2}^2 \mathbb{P}(\mathcal{E}_j^c)}$$

(S.14) $\quad \leq Ct^{-1}n^{3/2}(\lambda_1 \wedge \lambda_2)^{-1}h^{-1/2}\exp(-c_0\log^2 N)$.

The first and second inequalities are simply Markov inequality and Cauchy-Schwartz inequality, respectively, while the last one is by (S.7), Lemma S.1 and Lemma S.4.

Let $a_n = n^{-1}d_n c_K h^{-1/2}(\sqrt{n}\xi_\lambda + 1) = c_1 c_K^2 (nh)^{-1}r_n(\sqrt{n}\xi_\lambda + 1)$. For every $j \in 1, \ldots, s$. Combining (S.13) and (S.14), we can conclude that $\|R_n\| = o_P(a_n \log n)$ and accordingly we have $\|H_n\| = o_P(a_n \log n)$ and $\|G_n\| = o_P(a_n \log n)$. This completes the first part of the proof.

Next, we study $G_N = s^{-1} \sum_{j=1}^{s} G_n^{(j)}$. According to (S.13) and (S.14) again, we have

$$(S.15) \quad \mathbb{P}\bigg(\frac{n\|G_{n1}^{(j)}\|}{d_n c_K h^{-1/2}(\sqrt{n}\xi_\lambda + 1)} \geq t\bigg) \leq 2\big(\log(\sqrt{n}\xi_\lambda) + 2\big)\exp(-t^2/4),$$

and

$$(S.16) \quad \mathbb{P}\big(n\|G_{n2}^{(j)}\| \geq t\big) \leq C t^{-1} n^{3/2}(\lambda_1 \wedge \lambda_2)^{-1} h^{-1/2} \exp(-c_0 \log^2 N).$$

For any $1 \leq j \leq s$, define the event

$$\mathcal{A}_j = \big\{\|G_{n1}^{(j)}\| \leq a_n \log N\big\} \text{ and } \mathcal{A} = \cap_{j=1}^{s} \mathcal{A}_j.$$

By (S.15) and (S.16), we have

$$\mathbb{P}(\mathcal{A}_j^c) \leq \mathbb{P}\bigg(\frac{n\|G_{n1}^{(j)}\|}{a_n} \geq \frac{1}{2}\log N\bigg) + \mathbb{P}\bigg(n\|G_{n2}^{(j)}\| \geq \frac{1}{2}a_n \log N\bigg)$$

$$\leq 2\big(\log(\sqrt{n}\xi_\lambda) + 2\big)\exp(-\log^2 N/4)$$

$$(S.17) \qquad + 2C(a_n \log N)^{-1} n^{3/2}(\lambda_1 \wedge \lambda_2)^{-1} h^{-1/2} \exp(-c_0 \log^2 N).$$

According to the definition of $\mathcal{A}_j$, we have $\|G_n^{(j)} \mathbb{1}_{\mathcal{A}_j}\| \leq a_n \log N$ for all $j$. We again separate $G_n^{(j)}$ into two parts: $G_n^{(j)} = G_n^{(j)} \mathbb{1}_{\mathcal{A}_j} + G_n^{(j)} \mathbb{1}_{\mathcal{A}_j^c}$ for all $j$. Similar to (S.16), it follows from (S.7), Lemma S.1 and Lemma S.4 that

$$\mathbb{E}\|G_n^{(j)} \mathbb{1}_{\mathcal{A}_j^c}\| \leq 2\mathbb{E}\left\|g_{n,\lambda}^{(j)}(U)\mathbb{1}_{\mathcal{A}_j^c} K_U^\lambda\right\|$$

$$(S.18) \qquad \leq 2h^{-1/2}\sqrt{\mathbb{E}\|g_{n,\lambda}^{(j)}\|_{\mathcal{H}_1}^2 \mathbb{P}(\mathcal{A}_j^c)} \leq C(\lambda_1 h)^{-1/2} \mathbb{P}(\mathcal{A}_j^c).$$

Combining (S.13) and (S.14), we have $\mathbb{E}[\|G_n^{(j)} \mathbb{1}_{\mathcal{A}_j}\|] \leq \mathbb{E}\|G_n^{(j)}\| \leq 4a_n \log N$ for sufficiently large $N$. We next apply Lemma S.8 to $\big\{G_n^{(j)} \mathbb{1}_{\mathcal{A}_j} - \mathbb{E}[G_n^{(j)} \mathbb{1}_{\mathcal{A}_j}]\big\}_{j=1}^{s}$ since they are i.i.d. mean zero bounded variables. Note that $\|G_n^{(j)} \mathbb{1}_{\mathcal{A}_j} - \mathbb{E}[G_n^{(j)} \mathbb{1}_{\mathcal{A}_j}]\| \leq 2a_n \log N$ by (3.10) with (S.17), (S.18), and for sufficiently

large $N$, $\mathbb{E}\|G_n^{(j)} \mathbb{1}_{\mathcal{A}_j^c}\| \leq s^{-1/2} a_n \log^2 N$. Therefore, we have

$$\mathbb{P}\Big(\Big\|\frac{1}{s}\sum_{j=1}^s G_n^{(j)} \mathbb{1}_{\mathcal{A}_j}\Big\| > 6 a_n \log N\Big)$$

$$= \mathbb{P}\Big(\Big\|\frac{1}{s}\sum_{j=1}^s G_n^{(j)} \mathbb{1}_{\mathcal{A}_j} - \mathbb{E}[G_n^{(j)} \mathbb{1}_{\mathcal{A}_j}] + \mathbb{E}[G_n^{(j)} \mathbb{1}_{\mathcal{A}_j}]\Big\| > 6 a_n \log N\Big)$$

$$\leq \mathbb{P}\Big(\Big\|\frac{1}{s}\sum_{j=1}^s G_n^{(j)} \mathbb{1}_{\mathcal{A}_j} - \mathbb{E}[G_n^{(j)} \mathbb{1}_{\mathcal{A}_j}]\Big\| > 2 s^{-1/2} a_n \log^2 N\Big).$$

(S.19) $\quad \leq 2 \exp(-\log^2 N /2)$

Therefore, we can bound the tail probability of $G_N$:

$$\mathbb{P}\bigg(\Big\|\frac{1}{s}\sum_{j=1}^s G_n^{(j)}\Big\| > 12 a_n \log N\bigg)$$

$$\leq \mathbb{P}\bigg(\Big\|\frac{1}{s}\sum_{j=1}^s G_n^{(j)} \mathbb{1}_{\mathcal{A}_j}\Big\| > 6 a_n \log N\bigg) + \mathbb{P}\bigg(\Big\|\frac{1}{s}\sum_{j=1}^s G_n^{(j)} \mathbb{1}_{\mathcal{A}_j^c}\Big\| > 6 a_n \log N\bigg)$$

(S.20)
$$\leq 2 \exp(-\log^2 N/2) + (6 a_n \log N)^{-1} \mathbb{E}\|G_n^{(j)} \mathbb{1}_{\mathcal{A}_j^c}\| = o(1)$$

where the second inequality follows from (S.19) and (S.18), and the last equality follows by applying (3.10) to (S.17) and (S.18). Therefore, we complete the proof of the lemma.

**S.2.2. Proof of Lemma S.4.** We now prove a sharper rate on $f_{n,\lambda}$. Using the zero-order optimality condition (Luenberger, 1969),

(S.21)
$$\mathcal{L}_{n,\lambda}(f_{0j} + f_{n,\lambda}) - \mathcal{L}_{n,\lambda}(f_{0j})$$
$$= D\mathcal{L}_{n,\lambda}(f_{0j})[f_{n,\lambda}] + \frac{1}{2} D^2 \mathcal{L}_{n,\lambda}(f_{0j})[f_{n,\lambda}, f_{n,\lambda}] \leq 0.$$

We next bound the above two terms separately.

For the first term in the RHS of (S.21), Proposition A.2 gives the representation

$$\nabla \mathcal{L}_{n,\lambda}(f_{0j}) = \frac{1}{n}\sum_{i=1}^n \epsilon_i K_{U_i}^\lambda - W_\lambda f_{0j}.$$

For the regularization term $W_\lambda f_{0j}$, according to Proposition A.1, we have

$$\|W_\lambda f_{0j}\|^2 = \sum_{k=0}^\infty \frac{\lambda_1^2/\mu_{1,k}^2}{1+\lambda_1/\mu_{1,k}} \langle g_0, \phi_k\rangle_{L^2(\mathbb{P})}^2 + \sum_{k=0}^\infty \frac{\lambda_2^2/\mu_{2,k}^2}{\sigma_z^2 + \lambda_2/\mu_{2,k}} \langle \beta_{0j}, \phi_k\rangle_{L^2(\mathbb{P})}^2.$$

It follows from $g_0 \in \mathcal{H}_1$ that $\sum_{k=0}^{\infty} \mu_{1,k}^{-1} \langle g_0, \phi_k \rangle_{L^2(\mathbb{P})}^2 < \infty$. Since

$$\Big| \frac{\lambda_1/\mu_{1,k}^2}{1 + \lambda_1/\mu_{1,k}} \langle g_0, \phi_k \rangle_{L^2(\mathbb{P})}^2 \Big| \leq |\mu_{1,k}^{-1} \langle g_0, \phi_k \rangle_{L^2(\mathbb{P})}^2|$$

for any $k \in \mathbb{N}$, by dominate convergence theorem, we have

$$\lim_{N \to \infty} \sum_{k=0}^{\infty} \frac{\lambda_1/\mu_{1,k}^2}{1 + \lambda_1/\mu_{1,k}} \langle g_0, \phi_k \rangle_{L^2(\mathbb{P})}^2 = 0.$$

Similar result can be obtained for $\beta_{0j}$. Hence,

(S.22) $$\|W_\lambda f_{0j}\|^2 = o(\lambda_1 + \lambda_2).$$

Similarly, we can prove

(S.23) $$\|W_{\lambda,1} g_0\|^2 = o(\lambda_1) \text{ and } \|W_{\lambda,2} \beta_{0j}\|^2 = o(\lambda_2).$$

Define two events $\mathcal{B}_1 = \{\|f_{n,\lambda}\|_{\mathcal{H}^2} \leq c_0^{1/2} (\lambda_1 \wedge \lambda_2)^{-1/2}\}$ and

$$\mathcal{B}_2 = \{\|\sum_{i \in I_j} \epsilon_i K_{U_i}^\lambda\| \leq c_K (n/h)^{1/2} \log N\}.$$

According to (S.8), we have $\mathbb{P}(\mathcal{B}_1^c) \leq 2\exp(-c'n)$ for some constant $c'$. We denote $\epsilon_i^{(1)} = \epsilon_i \mathbb{1}_{\{\epsilon_i \leq \log N\}}$ and $\epsilon_i^{(2)} = \epsilon_i \mathbb{1}_{\{\epsilon_i > \log N\}}$. By (3.10), we have $\sqrt{n^2 h} N^{-c} = o(1)$ for some sufficiently large $c > 0$. On the other hand, since $\exp(-\log^2 N) = N^{-\log N}$, we have for sufficiently large $N$,

$$\frac{\sigma n^{1/2}}{c_K \log N} \exp(-\log^2 N) < \frac{1}{4} c_K (nh)^{-1/2} \log N.$$

Denoting $p_n = c_K (nh)^{-1/2} \log N$, we consider the decomposition

$$\mathbb{P}(\mathcal{B}_2^c) \leq \underbrace{\mathbb{P}\Big(\|\frac{1}{n} \sum_{i \in I_j} \epsilon_i^{(1)} K_{U_i}^\lambda\| > \frac{1}{2} p_n\Big)}_{P_1} + \underbrace{\mathbb{P}\Big(\|\frac{1}{n} \sum_{i \in I_j} \epsilon_i^{(2)} K_{U_i}^\lambda\| > \frac{1}{2} p_n\Big)}_{P_2}.$$

We first control $P_2$ by

$$P_2 \leq \frac{2\mathbb{E}\|\epsilon^{(2)} K_U^\lambda\|}{p_n} \leq 2 \frac{\sigma n^{1/2}}{c_K \log N} \exp(-\log^2 N),$$

where the first inequality is by Chebyshev's inequality and the second is due to the fact that $\|K_U^\lambda\| \leq c_K h^{-1/2}$ in (S.5) and $\epsilon$ is subgaussian. We next

control $P_1$. By Lemma S.8, there exists some constants $C_0$ such that

$$P_1 \leq \mathbb{P}\Big(\|\frac{1}{n}\sum_{i \in I_j} \epsilon_i^{(1)} K_{U_i}^\lambda - \mathbb{E}[\epsilon^{(1)} K_U^\lambda]\| > \frac{1}{2}p_n - \mathbb{E}\|\epsilon^{(1)} K_U^\lambda\|\Big)$$

$$\leq \mathbb{P}\Big(\|\frac{1}{n}\sum_{i \in I_j} \epsilon_i^{(1)} K_{U_i}^\lambda - \mathbb{E}[\epsilon^{(1)} K_U^\lambda]\| > \frac{1}{4}p_n\Big)$$

$$\leq 2\exp\left(-C_0 \log^2 N\right).$$

Combining the upper bounds of $P_1$ and $P_2$, we know that there exists some constant $C'$

$$\mathbb{P}(\mathcal{B}_2^c) \leq 2\exp\left(-C_0 \log^2 N\right) + 2\frac{\sigma n^{1/2}}{c_K \log N}\exp(-\log^2 N)$$

(S.24) $$\leq C' n^{1/2} \exp\left(-C_0 \log^2 N\right),$$

where the first inequality follows $\|n^{-1}\sum_{i \in I_j} \epsilon_i^{(1)} K_{U_i}^\lambda - \mathbb{E}[\epsilon^{(1)} K_U^\lambda]\| \leq c_K \log n h^{-1/2}$ when $\max_j \epsilon_j \leq \log n$ and $\epsilon$ is subgaussian. Therefore combining (S.22) and (S.24), we derive that on $\mathcal{B}_2$, there exists some constant $C_1$,

(S.25) 
$$|D\mathcal{L}_{n,\lambda}(f_{0j})[f_{n,\lambda}]| \leq \|\nabla \mathcal{L}_{n,\lambda}(f_{0j})\|\|f_{n,\lambda}\|$$
$$\leq C_1(\log N(nh)^{-1/2} + \lambda_1^{1/2} + \lambda_2^{1/2})\|f_{n,\lambda}\|.$$

For the second Fréchet derivative, we have

$$D^2\mathcal{L}_{n,\lambda}(f_{0j})[f_{n,\lambda}, f_{n,\lambda}]$$
$$= \langle \nabla^2 \mathcal{L}_{n,\lambda}(f_{0j})[f_{n,\lambda}], f_{n,\lambda}\rangle$$
$$= \frac{1}{n}\sum_{i \in I_j} f_{n,\lambda}(U_i)\langle K_{U_i}^\lambda, f_{n,\lambda}\rangle + \langle W_\lambda f_{n,\lambda}, f_{n,\lambda}\rangle$$
$$= \frac{1}{n}\sum_{i \in I_j} f_{n,\lambda}(U_i)\langle K_{U_i}^\lambda, f_{n,\lambda}\rangle - \mathbb{E}[f_{n,\lambda}(U)\langle K_U^\lambda, f\rangle]$$
$$\quad + \mathbb{E}[f_{n,\lambda}(U)\langle K_U^\lambda, f\rangle] + \langle W_\lambda f_{n,\lambda}, f\rangle$$

(S.26) $$= \frac{1}{n}\sum_{i \in I_j}\langle f_{n,\lambda}(U_i) K_{U_i}^\lambda - \mathbb{E}[f_{n,\lambda}(U) K_U^\lambda], f_{n,\lambda}\rangle + \|f_{n,\lambda}\|^2.$$

The empirical process in (S.26) is in the same form as the one in (S.10) with the only difference that $f_{n,\lambda}$ is in a different Hilbert space. Here we focus on bounding the supremum of the empirical process indexed by $\mathcal{H}_0 = \{f = (g, \beta) \in \mathcal{H} : \|f\|_{\sup} \leq 1, \|g\|_{\mathcal{H}}^2 \leq 1, \|\beta\|_{\mathcal{H}}^2 \leq 1\}$. We denote $q_n = (c_\mathcal{H} \wedge c_0^{1/2})(\lambda_1 \wedge \lambda_2)^{-1/2}$, where $c_\mathcal{H}$ is defined in Lemma S.1. Let $\mathring{f}_{n,\lambda} = q_n^{-1} f_{n,\lambda}$, and it follows from Lemma S.1 again that under the event $\mathcal{B}_1 = \{\lambda\|f_{n,\lambda}\|_{\mathcal{H}^2} \leq c_0\}$, we have $\|\mathring{f}_{n,\lambda}\|_{\sup} \leq 1$ and $\|\mathring{f}_{n,\lambda}\|_{\mathcal{H}^2} \leq 1$ and therefore $\mathring{f}_{n,\lambda} \in \mathcal{H}_0$. By the

definition of $\Psi(f, X)$ given in (S.11), we have
$$\Psi(\mathring{f}_{n,\lambda}, X) = c_K^{-1} h^{1/2} \mathring{f}_{n,\lambda} K_U^\lambda(x)$$
whose Lipchitz property is the same as (S.12).

By applying Lemma S.6 to
$$\mathbb{G}_n(\mathring{f}_{n,\lambda}) := \frac{1}{\sqrt{n}} \sum_{i \in I_j} \left( \Psi(\mathring{f}_{n,\lambda}, U) - \mathbb{E}[\Psi(\mathring{f}_{n,\lambda}, U)] \right).$$

with $d(\cdot) = \|\cdot\|_{\sup}$ and $\mathcal{F} = \mathcal{H}_0$, we have

$$\mathbb{P}\left( \mathcal{B}_1 \cap \left\{ \frac{\|\sum_{i \in I_j}[f_{n,\lambda}(U_i)K_{U_i}^\lambda - \mathbb{E}(f_{n,\lambda}(U)K_U^\lambda)]\|}{q_n c_K h^{-1/2}(\sqrt{n}\,\omega(\|q_n^{-1} f_{n,\lambda}\|_{\sup}) + 1)} \geq t \right\} \right)$$
$$= \mathbb{P}\left( \mathcal{B}_1 \cap \left\{ \frac{\sqrt{n}\|\mathbb{G}_n(\mathring{f}_{n,\lambda}^{(j)})\|}{\sqrt{n}\,\omega(\|\mathring{f}_{n,\lambda}\|_{\sup}) + 1} > t \right\} \right)$$
$$\leq \mathbb{P}\left( \sup_{f \in \mathcal{H}_0} \frac{\sqrt{n}\|\mathbb{G}_n(f)\|}{\sqrt{n}\,\omega(\|f\|_{\sup}) + 1} \geq t \right)$$
(S.27) $$\leq 2 \left( \log(\sqrt{n}\,\omega(1)) + 2 \right) \exp(-t^2/4),$$

where the last inequality comes from (S.47).

For $q'_n := c_K q_n h^{-1/2} \log N(\sqrt{n}\,\omega(\|q_n^{-1} f_{n,\lambda}\|_{\sup})) + 1)$, we define the event $\mathcal{B}_3 = \left\{ \|\sum_{i=1}^n [f_{n,\lambda}(U_i)K_{U_i}^\lambda - \mathbb{E}(f_{n,\lambda}(U)K_U^\lambda)]\| \leq q'_n \right\}$. By (S.27), we have $\mathbb{P}(\mathcal{B}_1 \cap \mathcal{B}_3^c) \leq 2 \left( \log(\sqrt{n}\,\omega(1)) + 2 \right) \exp(-\log^2 N/4)$. Therefore, under the event $\mathcal{B}_3$, we can bound the first term in (S.26) as

$$\left| \frac{1}{n} \sum_{i=1}^n \langle f_{n,\lambda}(U_i)K_{U_i}^\lambda - \mathbb{E}[f_{n,\lambda}(U)K_U^\lambda], f_{n,\lambda} \rangle \right|$$
$$\leq \left\| \frac{1}{n} \sum_{i=1}^n f_{n,\lambda}(U_i)K_{U_i}^\lambda - \mathbb{E}[f_{n,\lambda}(U)K_U^\lambda] \right\| \cdot \left\| f_{n,\lambda} \right\|$$
$$\leq c_K n^{-1} q_n h^{-1/2} \log N(\sqrt{n}\,\omega(\|q_n^{-1} f_{n,\lambda}\|_{\sup}) + 1)\|f_{n,\lambda}\|$$
(S.28) $$\leq c_K n^{-1} q_n h^{-1/2} \log N(\sqrt{n}\,\omega(q_n^{-1} h^{-1/2} \|f_{n,\lambda}\|) + 1)\|f_{n,\lambda}\|,$$

where the last inequality follows from Lemma S.1 and the fact that $\omega(\cdot)$ is non-decreasing.

Combining (S.21),(S.25) and (S.26), it yields that under the events $\mathcal{B}_2$

and $\mathcal{B}_3$ there exists a constant $C''$,

$$\|f_{n,\lambda}\|^2 \leq |D\mathcal{L}_{n,\lambda}(f_{0j})[f_{n,\lambda}]| + \Big|\frac{1}{n}\sum_{i\in I_j}\langle f_{n,\lambda}(U_i)K_{U_i}^\lambda - \mathbb{E}[f_{n,\lambda}(U)K_U^\lambda], f_{n,\lambda}\rangle\Big|.$$

$$\leq C_1(\log N(nh)^{-1/2} + \lambda^{1/2})\|f_{n,\lambda}\|$$
$$+ c_K n^{-1} q_n h^{-1/2} \log N(\sqrt{n}\,\omega(q_n^{-1}h^{-1/2}\|f_{n,\lambda}\|) + 1)\|f_{n,\lambda}\|.$$

The inequality above implies that either $\|f_{n,\lambda}\| \leq 2C_1(\log N(nh)^{-1/2} + \lambda_1^{1/2} + \lambda_2^{1/2})$ or

(S.29) $\qquad \|f_{n,\lambda}\| \leq 2c_K q_n h^{-1/2} n^{-1/2} \omega(q_n^{-1}h^{-1/2}\|f_{n,\lambda}\|) \log N$

According to (3.12), there exists a constant $c_r$ and $r_n = c_r(\log N(nh)^{-1/2} + \lambda_1^{1/2} + \lambda_2^{1/2})$ such that

$$r_n/\omega(q_n h^{-1/2} r_n) \geq 2c_K q_n n^{-1/2} h^{-1/2} n^{1/2} \log N.$$

It derives from (S.29) that $r_n/\omega(q_n h^{-1/2} r_n) \geq \|f_{n,\lambda}\|/\omega(q_n^{-1}h^{-1/2}\|f_{n,\lambda}\|)$. Since $\bar\omega_a$ in Assumption 3.3 is non-increasing, we have $\|f_n\| \leq r_n$.

We conclude that on the events $\mathcal{B}_2$ and $\mathcal{B}_3$, for the constant $C = c_r \vee 2C_1$, we have $\|f_{n,\lambda}\| \leq C\left(\frac{\log N}{\sqrt{nh}} + \lambda_1^{1/2} + \lambda_2^{1/2}\right)$. Therefore, for sufficiently large $n$, there exists a constant $c_0$,

$$\mathbb{P}\left(\|f_{n,\lambda}\| > C\left(\frac{\log N}{\sqrt{nh}} + \lambda_1^{1/2} + \lambda_2^{1/2}\right)\right)$$
$$\leq \mathbb{P}(\mathcal{B}_2^c) + \mathbb{P}(\mathcal{B}_3^c) \leq \mathbb{P}(\mathcal{B}_1^c) + \mathbb{P}(\mathcal{B}_1 \cap \mathcal{B}_3^c) + \mathbb{P}(\mathcal{B}_2^c)$$
$$\leq 2\exp(-c'n) + 2\left(\log(\sqrt{n}\,\omega(1)) + 2\right)\exp(-\log^2 N/4)$$
$$+ C'n^{1/2}\exp\left(-C_0 \log^2 N\right)$$

(S.30) $\qquad \leq 4n\exp(-c_0 \log^2 N).$

We therefore obtain the tail bound in Lemma S.4.

## APPENDIX S.3: RESULTS ON LIKELIHOOD RATIO TEST

We first prove Lemma C.1 on the rate of remainder term of $\text{APLRT}_{N,\lambda}$ in subsection S.3.1 and then prove Lemma C.2 on the power of $\text{APLRT}_{N,\lambda}$ in subsection S.3.2. In the last subsection S.3.3, we complete the remaining proofs of Theorems 4.1 and 4.2, i.e., on $\text{PLRT}_{n,\lambda}$ and $\text{PLRT}_{n,\lambda}^{\text{const}}$.

We start from a preliminary lemma.

**Lemma S.5.** Suppose (3.10), (3.11) and (3.12) hold and $\lambda = \lambda_1 \wedge \lambda_2 = o(1)$. Define the full-sample estimate $\widehat{g}_{N,\lambda}$ as $\bar{g}_{N,\lambda}$ in the case $s = 1$. We have $\|\bar{g}_{N,\lambda} - \widehat{g}_{N,\lambda}\| = O_P(N^{-1/2})$ and $\|\bar{g}_{N,\lambda} - g_0\| = O_P(r_N)$, where recall that $r_N = \log N(Nh_1)^{-1/2} + \lambda_1^{1/2} + \lambda_2^{1/2}$.

PROOF. We first prove that $\|\bar{g}_{N,\lambda} - \widehat{g}_{N,\lambda}\| = o_P(N^{-1/2})$. Denote $g_0^u = g_0 + W_{1,\lambda} g_0$. According to Lemma S.4, we have $\|\bar{g}_{N,\lambda} - g_0^u - \nabla \mathcal{L}_{N,\lambda}(g_0)\| = o_P(N^{-1/2})$. By the definition of $\widehat{g}_{N,\lambda}$ (noting that (3.10), (3.11) and (3.12) satisfy for $s=1$), it holds that $\|\widehat{g}_{N,\lambda} - g_0^u - \nabla \mathcal{L}_{N,\lambda}(g_0)\| = o_P(N^{-1/2})$. Therefore, we have

$$\|\bar{g}_{N,\lambda} - \widehat{g}_{N,\lambda}\| \le \|\bar{g}_{N,\lambda} - g_0^u - \nabla \mathcal{L}_{N,\lambda}(g_0)\| + \|\widehat{g}_{N,\lambda} - g_0^u - \nabla \mathcal{L}_{N,\lambda}(g_0)\|$$
$$(\text{S.31}) \qquad = o_P(N^{-1/2}).$$

Actually, by the proof of Lemma S.4, we can bound the tail probability that for some constant $C, C'$,

$$\mathbb{P}\left(\|\bar{g}_{N,\lambda} - \widehat{g}_{N,\lambda}\| \ge CN^{-1/2}\right)$$
$$\le \mathbb{P}\left(\|\bar{g}_{N,\lambda} - \widehat{g}_{N,\lambda}\| \ge 2s^{-1/2} a_n \log N\right)$$
$$\le \mathbb{P}\left(\|\bar{g}_{N,\lambda} - g_0^u - \nabla \mathcal{L}_{n,\lambda}(g_0)\| \ge s^{-1/2} a_n \log n\right)$$
$$+ \mathbb{P}\left(\|\widehat{g}_{N,\lambda} - g_0^u - \nabla \mathcal{L}_{n,\lambda}(g_0)\| \ge s^{-1/2} a_n \log n\right) = o(1),$$

where the last equality follows from Lemma B.1. Therefore, we have

$$\|\bar{g}_{N,\lambda} - g_0\| \le \|\widehat{g}_{N,\lambda} - g_0\| + \|\bar{g}_{N,\lambda} - \widehat{g}_{N,\lambda}\| = O_P(r_N + N^{-1/2}),$$

where the last equality is due to Lemma S.4 for $s=1$. Since $N^{-1/2} = O(r_N)$, the lemma is proven. $\square$

**S.3.1. Proof of Lemma C.1.** Denote $f_{0j} = (g_0, \beta_{0j})$, $\widehat{f}_0^{(j)} = (g_0, \widehat{\beta}_{n,\lambda}^{(j)})$, $g_{N,\lambda} = \bar{g}_{N,\lambda} - g_0$, $\check{f} = (g_{N,\lambda}, 0)$, $f_0^{(j)} = (0, \widehat{\beta}_{n,\lambda}^{(j)} - \beta_{0j})$. By applying the Taylor expansion to $\mathcal{L}_{n,\lambda}^{(j)}$ at $\widehat{f}_{n,\lambda}^{(j)}$, we have

$$(\text{S.32}) \quad \text{APLRT}_{N,\lambda} = \frac{1}{s} \sum_{j=1}^{s} \left( \mathcal{L}_{n,\lambda}^{(j)}(\widehat{f}_{n,\lambda}^{(j)}) - \mathcal{L}_{n,\lambda}^{(j)}(\widehat{f}_0^{(j)}) \right) = I_1 + I_2 + I_3,$$

where

$$I_1 := \frac{1}{s} \sum_{j=1}^{s} \langle \nabla \mathcal{L}^{(j)}(f_{0j}), \check{f} \rangle,$$
$$I_2 := \frac{1}{s} \sum_{j=1}^{s} \nabla^2 \mathcal{L}^{(j)}(f_{0j})[f_0^{(j)}, \check{f}],$$
$$I_3 := \frac{1}{s} \sum_{j=1}^{s} \frac{1}{2} \nabla^2 \mathcal{L}_{n,\lambda}^{(j)}(\widehat{f}_0^{(j)})[\check{f}, \check{f}].$$

We first decompose $I_3$ into two terms as follows

$$I_{31} := \frac{1}{s}\sum_{j=1}^{s}\frac{1}{2}\Big(\nabla^2\mathcal{L}_{n,\lambda}^{(j)}(\widehat{f}_0^{(j)})[\check{f},\check{f}] - \mathbb{E}\{\nabla^2\mathcal{L}_{n,\lambda}^{(j)}(\widehat{f}_0^{(j)})\}[\check{f},\check{f}]\Big),$$

$$I_{32} := \frac{1}{s}\sum_{j=1}^{s}\frac{1}{2}\mathbb{E}\{\nabla^2\mathcal{L}_{n,\lambda}^{(j)}(f_{0j})\}[\check{f},\check{f}].$$

According to Proposition A.2, we have

$$|I_{31}| \leq \frac{1}{s}\sum_{j=1}^{s}n^{-1}\Big\|\sum_{i=1}^{n}\check{f}(U_i)K_{U_i}^{\lambda} - \mathbb{E}[\check{f}(U)K_U^{\lambda}]\Big\|\|g_{N,\lambda}\|.$$

Lemma S.5 gives the rate of $\|g_{N,\lambda}\| = O_P(N^{-1/2})$ and (S.20) gives the rate of $s^{-1}\sum_{j=1}^{s}n^{-1}\|\sum_{i=1}^{n}\check{f}(U_i)K_{U_i}^{\lambda} - \mathbb{E}[\check{f}(U)K_U^{\lambda}]\| = O_P(12a_n r_N \log N)$. We can therefore control $|I_{31}|$ as

(S.33) $$\mathbb{P}\left(|I_{31}| > 12a_n r_N \log N\right) = o(1).$$

Therefore, $|I_{31}| = O_P(a_n r_N \log N)$. Applying Proposition A.2 again, the second term $I_{32} = \|\check{f}\|^2/2 = \|g_{N,\lambda}\|^2/2$. Since $a_n \log N = o(N^{-1/2})$ and $\sqrt{Nh} = o(1)$ yields $Na_n r_N \log N = o(h^{-1/2})$, we derive that

(S.34) $$2N \cdot I_3 = N\|g_{N,\lambda}\|^2 + o_P(h^{-1/2}).$$

Recall that $G_N = g_{N,\lambda} + s^{-1}\sum_{j=1}^{s}\nabla\mathcal{L}_1^{(j)}(f_{0j})$, where $\nabla\mathcal{L}_1^{(j)}(f_{0j})$ is the first coordinate of $\nabla\mathcal{L}_{n,\lambda}^{(j)}(f_{0j}) \in \mathcal{H}^2$, and satisfies

$$\nabla\mathcal{L}_1^{(j)}(f_{0j}) = -\frac{1}{n}\sum_{i\in I_j}\epsilon_i P_{X_i}^{\lambda} + W_{1,\lambda}g_0 \in \mathcal{H}_1.$$

By (S.20), we have $\|G_N\| = O_P(a_n \log N) = o_P(N^{-1/2})$ and

(S.35) $$2NI_3 = N\|g_{N,\lambda}\|^2 + o_P(h^{-1/2}) = N\Big\|\frac{1}{s}\sum_{j=1}^{s}\nabla\mathcal{L}_1^{(j)}(f_{0j})\Big\|^2 + o_P(h^{-1/2}).$$

Next we study the leading term $N\|s^{-1}\sum_{j=1}^{s}\nabla\mathcal{L}_1^{(j)}(f_{j,0})\|^2$ in the above equation. By Proposition A.2,

$$N\|s^{-1}\textstyle\sum_{j=1}^{s}\nabla\mathcal{L}_1^{(j)}(f_{j,0})\|^2 = N^{-1}\|\sum_{i=1}^{N}\epsilon_i P_{X_i}^{\lambda}\|^2 \\ -2\sum_{i=1}^{N}\epsilon_i(W_{1,\lambda}g_0)(X_i) + N\|W_{1,\lambda}g_0\|^2.$$

We bound the second term on the RHS above. It follows by the Fourier expansions of $W_{1,\lambda}g_0$ in Proposition A.1 that

$$\mathbb{E}\left\{|\sum_{i=1}^{N}\epsilon_i(W_{1,\lambda}g_0)(X_i)|^2\right\} = N\mathbb{E}\{\epsilon^2|W_{1,\lambda}g_0(X)|^2\}$$

(S.36) $$= \sigma^2 N \sum_{\nu \in \mathbb{Z}} |\langle g_0, \phi_\nu \rangle_{L^2\mathbb{P}}|^2 \left(\frac{\lambda_1/\mu_\nu}{1+\lambda_1/\mu_\nu}\right)^2 = o(N\lambda_1),$$

where the last equality is derived based on the similar arguments in the derivation of (S.22). So $\sum_{i=1}^{N}\epsilon_i(W_{1,\lambda}g_0)(X_i) = o_P((N\lambda_1)^{1/2}) = o_P(h^{-1/2})$. According to (S.22) again, $\|W_{1,\lambda}g_0\|^2 = o(\lambda_1)$. Consequently,

$$N\|s^{-1}\sum_{j=1}^{s}\nabla\mathcal{L}_1^{(j)}(f_{j,0})\|^2 = N^{-1}\|\sum_{i=1}^{N}\epsilon_i P_{X_i}^\lambda\|^2 + o(N\lambda) + o_P(h^{-1/2}).$$

In what follows, we study the statistical rate of $N^{-1}\|\sum_{i=1}^{N}\epsilon_i P_{X_i}^\lambda\|^2$. We write

$$\frac{1}{N}\left\|\sum_{i=1}^{N}\epsilon_i P_{X_i}^\lambda\right\|^2 = \frac{1}{N}\sum_{i=1}^{N}\epsilon_i^2 P^\lambda(X_i, X_i) + W(N),$$

where $W(N) = \sum_{i\neq j}\epsilon_i\epsilon_j P^\lambda(X_i, X_j)$. Here the $W(N)$ term is the leading term of (C.13).

Finally, we bound $N^{-1}\sum_{i=1}^{N}\epsilon_i^2 P^\lambda(X_i, X_i)$ by studying the moment

$$\mathbb{E}\left[\frac{1}{N}\sum_{i=1}^{N}\epsilon_i^2 P^\lambda(X_i, X_i)\right]^2$$

$$= \mathbb{E}\left[\frac{1}{N}\sum_{i=1}^{N}\epsilon_i^2 P^\lambda(X_i, X_i) - \mathbb{E}[\epsilon^2 P^\lambda(X, X)]\right]^2 + (\mathbb{E}[\epsilon^2 P^\lambda(X, X)])^2$$

$$\leq 2N^{-1}\mathbb{E}[\epsilon^4 P^\lambda(X, X)^2] + (\mathbb{E}[\epsilon^2 P^\lambda(X, X)])^2 = h^{-1}\sigma_P^2 + o(1),$$

where the last equality is due to $\mathbb{E}[\epsilon^4 P^\lambda(X, X)^2] \leq O(h^{-2})$ and the definition of $\sigma_P^2$. Then, as $N^{-1}h^{-2} = o(1)$, we can derive that

(S.37) $$\frac{1}{N}\sum_{i=1}^{N}\epsilon_i^2 P^\lambda(X_i, X_i) = h^{-1}\sigma_P^2 + o_P(1).$$

According to (S.36) and (S.37), we have

(S.38) $$N\left\|s^{-1}\sum_{j=1}^{s}\nabla\mathcal{L}_1^{(j)}(f_{j,0})\right\|^2 = O_P(h^{-1} + N\lambda + h^{-1/2}) = O_P(h^{-1}).$$

Now we turn to bound $I_1$. By the definition of $\check{f} = (g_{N,\lambda}, 0)$, we have

$$I_1 = -\left\|\frac{1}{s}\sum_{j=1}^{s}\nabla\mathcal{L}_1^{(j)}(f_{0j})\right\|^2 + \frac{1}{s}\sum_{j=1}^{s}\langle\nabla\mathcal{L}_1^{(j)}(f_{0j}), G_N\rangle$$

(S.39) $$= -\left\|\frac{1}{s}\sum_{j=1}^{s}\nabla\mathcal{L}_1^{(j)}(f_{0j})\right\|^2 + o_P((N^2h)^{-1/2}),$$

where the last inequality follows from (S.38) and (S.20). Finally, we calculate $I_2$ by

(S.40) $$I_2 = \frac{1}{s}\sum_{j=1}^{s}\langle f_0^{(j)}, \check{f}\rangle = \frac{1}{s}\sum_{j=1}^{s}\left\langle(0, \widehat{\beta}_{n,\lambda}^j - \beta_{0j}), (\bar{g}_{N,\lambda} - g_0, 0)\right\rangle = 0.$$

Combining (S.35), (S.39) and (S.40), it follows that

(S.41) $$-2N \cdot \text{APLRT}_{N,\lambda} = N\left\|\frac{1}{s}\sum_{j=1}^{s}\nabla\mathcal{L}_1^{(j)}(f_{0j})\right\|^2 + o_P(h^{-1/2}).$$

Moreover, using (S.38), (S.37) and (S.36), we have

$$-2N \cdot \text{APLRT}_{N,\lambda} = W(N) + h^{-1}\sigma_P^2 + N\|W_{1,\lambda}g_0\|^2 + o_P(h^{-1/2})$$

Therefore, we complete the proof of the lemma.

**S.3.2. Proof of Lemma C.2.** We separate the $\text{APLRT}_{N,\lambda}$ into two parts: $\text{APLRT}_{N,\lambda} =: T_1 - T_2$, where

$$T_1 := \frac{1}{s}\sum_{j=1}^{s}\left\{\mathcal{L}_{n,\lambda}^{(j)}\big((\bar{g}_{N,\lambda}, \widehat{\beta}_{n,\lambda}^{(j)})\big) - \mathcal{L}_{n,\lambda}^{(j)}\big((g_N, \widehat{\beta}_{n,\lambda}^{(j)})\big)\right\},$$

$$T_2 := \frac{1}{s}\sum_{j=1}^{s}\left\{\mathcal{L}_{n,\lambda}^{(j)}\big((g_0^*, \widehat{\beta}_{n,\lambda}^{(j)})\big) - \mathcal{L}_{n,\lambda}^{(j)}\big((g_N, \widehat{\beta}_{n,\lambda}^{(j)})\big)\right\}.$$

We next bound $T_1$ and $T_2$. Note that $T_1$ can be viewed as APLRT statistic for testing $H_0 : g_0 = g_N$ v.s. $H_1 : g_0 \neq g_N$. We next check the proof of Theorem 4.5 to show that $(2u_N)^{-1/2}(r_P T_1 - N\|W_{1,\lambda}g_N\|^2 - u_N) = O_P(1)$ uniformly over $\beta_{01}, \ldots, \beta_{0s}, g_N \in \mathcal{F}_\zeta$. In the proofs of Lemma S.13 and Theorem 4.5, the results are only related to the magnitudes of $\|g_N\|_\mathcal{H}, \|W_{1,\lambda}g_N\|$ and $\|\beta_{0j}\|_\mathcal{H}, \|W_{2,\lambda}\beta_{0j}\|$ for $1 \leq j \leq s$. Other constants are only dependent on the kernel function, eigen-functions ord $\epsilon$. By the definition of $\mathcal{F}_\zeta$, we can uniformly bound $\|g_N\|_\mathcal{H} \leq \zeta$ and $\|W_{1,\lambda}g_N\|^2 \leq \lambda\|g_N\|^2 \leq \zeta^2$ by (S.22) and similarly for $\beta_{0j}$ for all $j = 1, \ldots, s$.

We next give a uniform bound of $T_2$. Denoting $f_N = (g_N, \widehat{\beta}_{n,\lambda}^{(j)})$, $\Delta g =$

$g_N - g_0^*$ and $f_{N0} = (g_N - g_0^*, 0)$, we compute the expectation

$$\mathbb{E}\left[\mathcal{L}_{n,\lambda}^{(j)}((g_0^*, \widehat{\beta}_{n,\lambda}^{(j)})) - \mathcal{L}_{n,\lambda}^{(j)}((g_N, \widehat{\beta}_{n,\lambda}^{(j)}))\right]$$

$$= \mathbb{E}\left[\nabla \mathcal{L}_{n,\lambda}^{(j)}(f_N)[f_{N0}] + \frac{1}{2}\nabla^2 \mathcal{L}_{n,\lambda}^{(j)}(f_N)[f_{N0}, f_{N0}]\right]$$

(S.42) $\qquad = \mathbb{E}\langle \nabla \mathcal{L}_{n,\lambda}^{(j)}(f_N), f_{N0}\rangle + \lambda_1 \langle \Delta g, g_N\rangle_{\mathcal{H}} + \frac{1}{2}\|\Delta g\|^2.$

Recall that $H_n^{(j)} = \widehat{\beta}_{n,\lambda}^{(j)} - \beta_{j0}^u + n^{-1}\sum_{i\in I_j} \epsilon_i Z_i T_{X_i}^\lambda$. By Cauchy-Schwarz inequality, Lemma S.1 and Lemma S.4, we have

(S.43)
$$\mathbb{E}\langle \nabla \mathcal{L}_{n,\lambda}^{(j)}(f_N), f_{N0}\rangle$$
$$= \mathbb{E}\langle \nabla \mathcal{L}_{n,\lambda}^{(j)}((g_N, \beta_{j0})), f_{N0}\rangle + \langle \mathbb{E}\nabla^2 \mathcal{L}_{n,\lambda}^{(j)}((g_N, \beta_{j0}))[\widehat{\beta}_{n,\lambda}^{(j)} - \beta_{j0}], \Delta g\rangle$$
$$\leq \langle W_{\lambda,2}\beta_{j0}, \Delta g\rangle + \mathbb{E}\|H_n^{(j)}\|\|\Delta g\| \leq o(\sqrt{\lambda_2})\|\Delta g\|,$$

where the last inequality is due to (S.22), Lemma S.4 and $a_n \log N = o_P(\sqrt{\lambda_2})$. Since $\Delta g \in \mathcal{F}_\zeta$, we also have $\langle \Delta g, g_N\rangle_{\mathcal{H}} \leq |\langle g_0^*, \Delta g\rangle_{\mathcal{H}}| + \langle \Delta g, \Delta g\rangle_{\mathcal{H}} \leq \|g_0^*\|_{\mathcal{H}}\zeta^{1/2} + \zeta$.

Denote $\Delta \beta_{n,\lambda}^{(j)} = \widehat{\beta}_{n,\lambda}^{(j)} - \beta_{0j}$, the variance of $T_2$ is

$$\mathbb{E}\left[|T_2 - \mathbb{E}[T_2]|^2\right] \leq \frac{1}{s^2}\sum_{j=1}^{s}\mathbb{E}\left[\mathcal{L}_{n,\lambda}^{(j)}((g_0^*, \widehat{\beta}_{n,\lambda}^{(j)})) - \mathcal{L}_{n,\lambda}^{(j)}((g_N, \widehat{\beta}_{n,\lambda}^{(j)}))\right]^2$$

$$\leq \frac{n}{N^2}\sum_{j=1}^{s}\mathbb{E}\left[|2\epsilon_1\Delta g(X_1) + 2Z_1\Delta\beta_{n,\lambda}^{(j)}(X_1)\Delta g(X_1) + Z_1^2\Delta g(X_1)^2|^2\right]$$

(S.44) $\quad \leq O(\|\Delta g\|^2 + r_n^2\|\Delta g\|^2 + \|\Delta g\|^4),$

where the last inequality holds uniformly according to the definition of $\mathcal{F}_\zeta$ that $\text{Var}(\Delta g^2) \leq \zeta \mathbb{E}^2[g_N^2]$. Combining (S.42), (S.43) and (S.44), we have

$$\begin{aligned}-2NT_2 &= -N\|\Delta g\|^2 + O_P(N\lambda_1 + N\sqrt{\lambda_2}\|\Delta g\| \\ &\quad + N^{1/2}\|\Delta g\| + N^{1/2}r_n\|\Delta g\| + N^{1/2}\|\Delta g\|^2),\end{aligned}$$

uniformly for $\beta_{01}, \ldots, \beta_{0s}, \Delta g \in \mathcal{F}_\zeta$.

With the two parts ready, the normalized averaged statistic

$$(2u_N)^{-1/2}(-2Nr_P \cdot \text{APLRT}_{N,\lambda} - u_N)$$
$$\geq N(2u_N)^{-1/2}\|\Delta g\|^2(1 + O_P(\lambda_1\|\Delta g\|^{-2} + \sqrt{\lambda_2}\|\Delta g\|^{-1}$$
$$+ N^{-1/2} + N^{-1/2}\|\Delta g\|^{-1})) + (2u_N)^{-1/2}N\|W_{1,\lambda}g_0\|^2 + O_P(1)$$
$$\overset{(\text{S.22})}{\geq} N(2u_N)^{-1/2}\|\Delta g\|^2(1 + O_P(\lambda_1\|\Delta g\|^{-2} + \sqrt{\lambda_2}\|\Delta g\|^{-1}$$
$$+ N^{-1/2} + N^{-1/2}\|\Delta g\|^{-1})) + O_P(1),$$

uniformly for $\beta_{01}, \ldots, \beta_{0s}, \Delta g \in \mathcal{F}_\zeta$.

**S.3.3. Remaining Proofs of Theorems 4.1 and 4.2.** To complete the proofs of Theorems 4.1 and 4.2, we prove the properties of $\text{PLRT}_{n,\lambda}$. Also, we prove the properties of $\text{PLRT}_{n,\lambda}^{\text{const}}$ defined in Remark 4.3.

• **Analysis of $\text{PLRT}_{n,\lambda}$.** We follow a similar proof strategy as Theorem 4.5 and only outline their difference here. Recall that the first coordinate of the derivative $\nabla\mathcal{L}_{n,\lambda}^{(1)}((g_0,\beta_{0j}))$ is denoted as

$$\nabla\mathcal{L}_2^{(1)}(f_{0j}) = -\frac{1}{n}\sum_{i\in I_j}\epsilon_i Z_i T_{X_i}^\lambda + W_{2,\lambda}\beta_{0j} \in \mathcal{H}.$$

By the Taylor expansion, we have

$$-2n \cdot \text{PLRT}_{n,\lambda} = n\|\nabla\mathcal{L}_2^{(1)}(f_{0j})\|^2 + o_P(h_2^{-1/2}).$$

Expanding $\|\nabla\mathcal{L}_2^{(1)}(f_{0j})\|^2$ yields the U-statistic

$$W_3^{(j)}(n) = \sum_{u<v, u,v\in I_j} 2\epsilon_u\epsilon_v Z_u Z_v T^\lambda(X_u, X_v)$$

with $\text{Var}(W_3^{(j)}(n)) \sim 2h_2^{-1}n^2\rho_T^2$. Comparing with $W_3(2n)$ in (C.7), $W_3^{(j)}(n)$ is of the same type but with half sample size of $W_3(2n)$. Therefore, applying the same procedure as the analysis of $\text{PLRT}_{n,\lambda}^{(j,j')}$. We have the following asymptotic normality

$$(2h_1^{-1}N\rho_P)^{-1/2}W_3^{(j)}(n) \rightsquigarrow N(0,1).$$

Based on $\mathbb{E}[\epsilon^2 Z_i^2 T^\lambda(X,X)] = h^{-1}\sigma_T^2$, we have $n^{-1}\sum_{i\in I_j}\epsilon_i^2 Z_i^2 T^\lambda(X_i, X_i) = h^{-1}\sigma_T^2 + O_P(1)$ In summary, we have

$$(2h_2^{-1}\sigma_T^4/\rho_T^2)^{-1/2}(-2nr_T \cdot \text{PLRT}_{n,\lambda} - h_2^{-1}\sigma_T^4/\rho_T^2) \rightsquigarrow N(0,1).$$

The power analysis follows the proof of Theorem 4.6 and we can obtain

$$\left|\frac{-2nr_T \cdot \text{PLRT}_{n,\lambda} - u_n}{(2u_n)^{1/2}}\right| \geq \frac{n\|\beta_n\|^2}{(2u_n)^{1/2}}\left[1 + O_P\left(\frac{\lambda_1}{\|\beta_n\|^2} + \frac{(r_N + n^{-1/2})}{\|\beta_n\|} + \frac{1}{n^{1/2}}\right)\right] + O_P(1),$$

where $O_P$ is uniformly for $\beta_n \in \mathcal{F}_\zeta$. When $\|\beta_n\| \geq c_n\sqrt{\lambda_2 + r_N^2 + (nh_2^{1/2})^{-1}}$, the probability is large enough.

• **Analysis of** $\text{PLRT}_{n,\lambda}^{\text{const}}$. We first do the following decomposition: $\text{PLRT}_{n,\lambda}^{\text{const}} = L_1 - L_2$, where

$$L_1 := \mathcal{L}_{n,\lambda}^{(1)}\big((\bar{g}_{N,\lambda}, \widehat{\beta}_{n,\lambda}^{(1)})\big) - \mathcal{L}_{n,\lambda}^{(1)}\big((g_0, c)\big),$$
$$L_2 := \mathcal{L}_{n,\lambda}^{(1)}\big((\widehat{g}_c^{(1)}, \widehat{c})\big) - \mathcal{L}_{n,\lambda}^{(1)}\big((g_0, c)\big).$$

Cheng and Shang (2015) showed that

$$(2h_2^{-1}\rho_P^2)^{-1/2}\left(-2nL_2 - h^{-1}\sigma_P^2\right) \rightsquigarrow N(0,1).$$

As for $L_1$, we obtain

$$-2nL_1 = n\|\nabla\mathcal{L}_{n.\lambda}^{(1)}(f_{0j})\|^2 + o_P(h^{-1/2})$$

by following similar analysis used in $\text{PLRT}_{n,\lambda}$. The U-statistic after expanding $\|\nabla\mathcal{L}_{n.\lambda}^{(1)}(f_{0j})\|^2$ turns out to be

$$W_4(n) = \sum_{i<j} 2\epsilon_i\epsilon_j(P^\lambda(X_i, X_j) + Z_iZ_jT^\lambda(X_i, X_j))$$

with $\text{Var}(W_4(n)) = 2n^2h^{-1}(\rho_P^2 + \rho_T^2)$. Combining that $h\mathbb{E}[\epsilon^2(P^\lambda(X,X) + Z^2T^\lambda(X,X))] = \sigma_P^2 + \sigma_T^2$, we have

$$(2h_2^{-1}(\rho_P^2 + \rho_T^2))^{-1/2}\left(-2nL_1 - h_2^{-1}(\sigma_P^2 + \sigma_T^2)\right) \rightsquigarrow N(0,1).$$

Combining with the limiting distribution of $L_2$, we complete the whole proof.

## APPENDIX S.4: AUXILIARY LEMMAS

Define an empirical process

(S.45) $$\mathbb{G}_n(f) = \frac{1}{\sqrt{n}}\sum_{i=1}^n[\Psi(f, X_i) - \mathbb{E}(\Psi(f, X))],$$

where the operator $\Psi : \mathcal{H} \times \mathbb{R} \to \mathcal{H}$.

**Lemma S.6.** Let a function space $\mathcal{F} = \{f \in \mathcal{H}' \subseteq \mathcal{H} : d(f) \leq 1\}$, where $d(\cdot) : \mathcal{F} \to \mathbb{R}_+$ is the norm of interest. Define

$$\omega_d(\delta; \mathcal{F}) = \int_0^\delta \sqrt{\log(1 + \mathcal{N}(\mathcal{F}; d, \epsilon))}d\epsilon.$$

If $\Psi$ is 1-Lipschitz with respect to $d$, i.e., $\|\Psi(f,X) - \Psi(g,X)\| \leq d(f-g)$ for any $f, g \in \mathcal{F}$, we have,

(S.46)
$$\mathbb{P}\left(\sup_{f \in \mathcal{F}} \frac{\sqrt{n}\|\mathbb{G}_n(f)\|}{\sqrt{n}\,\omega_d(d(f);\mathcal{F}) + 1} > t\right) \leq 2\left(\log\left(\sqrt{n}\omega_d(1;\mathcal{F})\right) + 2\right)\exp(-t^2/4).$$

In particular, we have

(S.47)
$$\mathbb{P}\left(\sup_{f \in \mathcal{F}} \frac{\sqrt{n}\|\mathbb{G}_n(f)\|}{\sqrt{n}\,\omega_d(1;\mathcal{F}) + 1} > t\right) \leq 2\left(\log\left(\sqrt{n}\omega_d(1;\mathcal{F})\right) + 2\right)\exp(-t^2/4).$$

**Remark S.7.** The above concentration inequality is similar as the tail inequality of continuity moduli of empirical processes, which was initiated by Giné et al. (2003); Giné and Koltchinskii (2006) in a measurable function class. Our result here focuses on Hilbert spaces, and fits into the RKHS framework in consideration.

PROOF. For simplicity, we will write $\omega_d(\delta;\mathcal{F})$ as $\omega(\delta)$ in the proof. Since for any $f, g \in \mathcal{F}$, $\|\Psi(f,X) - \Psi(g,X)\| \leq d(f-g)$, it follows Theorem 3.5 of Pinelis (1994), for any $t > 0$,

$$\mathbb{P}\left(\frac{\|\mathbb{G}_n(f) - \mathbb{G}_n(g)\|}{d(f-g)} \geq t\right) \leq 2\exp\left(-\frac{t^2}{8}\right).$$

Therefore, $\mathbb{G}_n(f)$ is a subgaussian empirical process such that by Lemma 8.1 in Kosorok (2008), we have $\|\|\mathbb{G}_n(g) - \mathbb{G}_n(f)\|\|_{\phi_2} \leq 8d(f-g)$, where $\|\cdot\|_{\phi_2}$ is the Orlicz norm associated with $\phi_2(s) \equiv \exp(s^2) - 1$. Define $\mathcal{F}_\delta = \{f \in \mathcal{F} \mid d(f) \leq \delta\}$. Theorem 8.4 of Kosorok (2008) implies that there exists a constant $C_1$, for an arbitrary $\delta > 0$,

$$\left\|\sup_{f-g \in \mathcal{F}_\delta} \|\mathbb{G}_n(f) - \mathbb{G}_n(f)\|\right\|_{\phi_2} \leq C_1 \int_0^\delta \sqrt{\log(1 + N(\mathcal{F};d,\epsilon))}d\epsilon = C_1\omega(\delta).$$

We apply Lemma 8.1 in Kosorok (2008) again and obtain,

(S.48)
$$P\left(\sup_{\substack{f \in \mathcal{F} \\ d(f) \leq \delta}} \|\mathbb{G}_n(g)\| \geq t\right) \leq 2\exp\left(-\frac{t^2}{C_1^2\omega(\delta)^2}\right).$$

We next apply the peeling argument. Define the set $V_0 := \{f \in \mathcal{F} : \sqrt{n}\omega(d(f)) \leq 2\}$ and

$$V_m = \{f \in \mathcal{F} : 2^{m-1} \leq \sqrt{n}\omega(d(f)) \leq 2^m\} \text{ for } 1 \leq m \leq M,$$

where the index $m$ ranges up to $M := \lceil \log \sqrt{n}\omega(1) \rceil$ as $\omega(\cdot)$ is increasing and

$d(f) \leq 1$ for all $f \in \mathcal{F}$. Using union bound, we have

$$\mathbb{P}\left(\sup_{f \in \mathcal{F}} \frac{\sqrt{n}||\mathbb{G}_n(f)||}{\sqrt{n}\omega(d(f)) + 1} > t\right)$$

$$\leq \sum_{m=0}^{M} \mathbb{P}\left(\sup_{f \in V_m} \frac{\sqrt{n}||\mathbb{G}_n(f)||}{\sqrt{n}\omega(d(f)) + 1} \geq t\right)$$

$$\leq \mathbb{P}\left(\sup_{d(f) \leq \omega^{-1}(2/\sqrt{n})} \sqrt{n}||\mathbb{G}_n(f)|| \geq t\right)$$

$$+ \sum_{m=1}^{M} \mathbb{P}\left(\sup_{d(f) \leq \omega^{-1}(2^m/\sqrt{n})} \sqrt{n}||\mathbb{G}_n(f)|| \geq (1 + 2^{m-1})t\right)$$

$$\leq 2\exp\left(-\left(\frac{t/\sqrt{n}}{2C_1/\sqrt{n}}\right)^2\right) + \sum_{m=1}^{M} 2\exp\left(-\left(\frac{(1 + 2^{m-1})t/\sqrt{n}}{C_1 2^m/\sqrt{n}}\right)^2\right)$$

$$\leq 2(M+1)\exp\left(-\frac{t^2}{4C_1^2}\right).$$

The second last inequality follows from (S.48). We can therefore derive (S.47) by the fact that $d(f) \leq 1$ for all $f \in \mathcal{F}$ and $\omega(\cdot)$ is an increasing function. □

**Lemma S.8.** (Pinelis, 1994) If $X_1, \ldots, X_n$ are zero mean independent random variables in a separable Hilbert space and $||X_i|| \leq C$ for $i = 1, \ldots, n$, then

$$\mathbb{P}\left(\left\|\frac{1}{n}\sum_{i=1}^{n} X_i\right\| > C\sqrt{\frac{2t}{n}}\right) < 2e^{-t}.$$

**Lemma S.9.** Define $\ell = (\sqrt{5} - 1)/2$, $I_1 = \pi/(2m\sin(\pi/(2m)))$ and $I_2 = (2m-1)\pi/(4m^2\sin(\pi/(2m)))$. Under $\mathcal{H}_0^m[0,1]$, we choose $\lambda_i = h_i^{2m}$ for $i = 1, 2$, while under the Gaussian RKHS $\mathcal{H}_G$ we choose $h_i = (\log(1/\lambda_i))^{-1/p}$ for $i = 1, 2$. We obtain the following constants needed in the inference results:

Sobolev space : $\sigma_g^2 \sim \dfrac{\sigma^2 I_2}{\pi}, \sigma_\beta^2 \sim \dfrac{\sigma^2 I_2}{\sigma_z^{2-1/m}\pi}, \rho_P^2 \sim \dfrac{\sigma^4 I_2}{\pi}, \sigma_P^2 \sim \dfrac{\sigma^2 I_1}{\pi},$

$\rho_T^2 \sim \dfrac{\sigma^4 \sigma_z^{1/m} I_2}{\pi}, \sigma_T^2 \sim \dfrac{\sigma^2 \sigma_z^{1/m} I_1}{\pi};$

Gaussian RKHS : $\sigma_g^2 \leq 1.7178\sigma^2, \sigma_\beta^2 \leq 1.7178(\sigma/\sigma_z)^2, \rho_P^2 \sim \sigma^4 \log(1/\ell),$

$\sigma_P^2 \sim \sigma^2 \log(1/\ell), \rho_T^2 \sim \sigma^4 \log(1/\ell)$ and $\sigma_T^2 \sim \sigma^2 \log(1/\ell).$

PROOF. For the Sobolev space, by definitions, we have

$$\sigma_g^2 = \lim_{N\to\infty} \sigma^2 h_1\Big(1 + 2\sum_{k=0}^{\infty}(1 + \lambda_1(2\pi k)^{2m})^{-2}\Big),$$

$$\rho_P^2 = \lim_{N\to\infty} \sigma^4 h_1\Big(1 + 2\sum_{k=0}^{\infty}(1 + \lambda_1(2\pi k)^{2m})^{-2}\Big),$$

$$\sigma_P^2 = \lim_{N\to\infty} \sigma^2 h_1\Big(1 + 2\sum_{k=0}^{\infty}(1 + \lambda_1(2\pi k)^{2m})^{-1}\Big).$$

We calculate the infinite summation

$$\sum_{k=1}^{\infty} \frac{2\pi h_1}{(1+(2\pi h_1 k)^{2m})^l} \leq \sum_{k=0}^{\infty} \int_{2\pi h_1 k}^{2\pi h_1(k+1)} \frac{1}{(1+x^{2m})^l}dx \leq \int_0^{\infty} \frac{1}{(1+x^{2m})^l}dx$$

and similarly, we have

$$\sum_{k=0}^{\infty} \frac{2\pi h_1}{(1+(2\pi h_1 k)^{2m})^l} \geq \int_0^{\infty} \frac{1}{(1+x^{2m})^l}dx.$$

Therefore we have the asymptotic variance

$$\sigma_g^2 \sim \frac{\sigma^2}{\pi}\int_0^{\infty}\frac{1}{(1+x^{2m})^2}dx = \frac{\sigma^2 I_2}{\pi}.$$

Similar arguments apply to $\rho_P^2$ and $\sigma_P^2$.

Recall the definitions of the remaining three constants:

$$\sigma_\beta^2 = \lim_{N\to\infty} \sigma^2 \sigma_z^2 h_2\Big(1 + 2\sum_{k=0}^{\infty}(\sigma_z^2 + \lambda_2(2\pi k)^{2m})^{-2}\Big),$$

$$\rho_T^2 = \lim_{N\to\infty} \sigma^4 \sigma_z^4 h_2\Big(1 + 2\sum_{k=0}^{\infty}(\sigma_z^2 + \lambda_2(2\pi k)^{2m})^{-2}\Big),$$

$$\sigma_T^2 = \lim_{N\to\infty} \sigma^2 \sigma_z^2 h_2\Big(1 + 2\sum_{k=0}^{\infty}(\sigma_z^2 + \lambda_2(2\pi k)^{2m})^{-1}\Big).$$

The computation of the above constants can be reduced to those of the previous ones by the change of variable. For example, we rewrite $\sigma_\beta^2$ by

$$\sigma_\beta^2 \sim \sum_{k=1}^{\infty} \frac{2\pi h_2 \sigma^2 \sigma_z^{-2}}{(1+\sigma_z^{-2}(2\pi h_2 k)^{2m})^l} \sim \frac{\sigma^2}{\pi\sigma_z^2}\int_0^{\infty}\frac{dx}{(1+\sigma_z^{-2}x^{2m})^2} = \frac{\sigma^2 I_2}{\pi\sigma_z^{2-1/m}}$$

We can evaluate $\sigma_T^2$ and $\rho_T^2$ in the same way.

As for the Gaussian RKHS, we have

$$\sigma_g^2 = \lim_{N\to\infty} \sigma^2 h_1 \sum_{k=0}^{\infty} \left(1 + \lambda_1 \ell \exp(-2\ln\ell \cdot k)\right)^{-2},$$

$$\rho_P^2 = \lim_{N\to\infty} \sigma^4 h_1 \sum_{k=0}^{\infty} \left(1 + \lambda_1 \ell \exp(-2\ln\ell \cdot k)\right)^{-2},$$

$$\sigma_P^2 = \lim_{N\to\infty} \sigma^2 h_1 \sum_{k=0}^{\infty} \left(1 + \lambda_1 \ell \exp(-2\ln\ell \cdot k)\right)^{-1}.$$

To be more general, we set $h_1 = (\log(1/\lambda_1))^{-1/p}$ and will show that

$$\sum_{k=0}^{\infty} \frac{h_1}{(1 + \lambda_1 \beta \exp(\alpha k^p))^l} \sim \alpha^{-1/p}.$$

For $p = 1$, it gives the result for Gaussian RKHS.

On one hand, according to the convexity, we have

$$\sum_{k=1}^{\infty} \frac{h_1}{(1 + \lambda_1 \beta \exp(\alpha k^p))^l} \leq \int_0^{\infty} \frac{dx}{(1 + \lambda_1 \beta \exp(\alpha x^p))^l}$$

We can further study the integration

$$\int_0^{\infty} \frac{dx}{(1 + \lambda_1 \beta \exp(\alpha x^p))^l}$$
$$= \int_0^{(\alpha^{-1}\log(1/\lambda_1))^{1/p}} + \int_{(\alpha^{-1}\log(1/\lambda_1))^{1/p}}^{\infty} \frac{dx}{(1 + \lambda_1 \beta \exp(\alpha x^p))^l}$$
$$\leq (\alpha^{-1}\log(1/\lambda_1))^{1/p} + \int_{(\alpha^{-1}\log(1/\lambda_1))^{1/p}}^{\infty} (\lambda_1 \beta)^{-l} \exp(-\alpha l x^p) dx$$
$$\leq (\alpha^{-1}\log(1/\lambda_1))^{1/p} + \alpha p^{-1} \beta^{-l} p (\alpha^{-1}\log(1/\lambda_1))^{-(p-1)/p}.$$

On the other hand, for any $t \in (0, 1)$,

$$\sum_{k=1}^{\infty} \frac{h_1}{(1 + \lambda_1 \beta \exp(\alpha k^p))^l} \geq \int_1^{((t/\alpha)\log(1/\lambda_1))^{1/p}} \frac{dx}{(1 + \lambda_1 \beta \exp(\alpha x^p))^l}$$
$$\geq \frac{1}{(1 + \lambda_1^{1-t}\beta)^l} \left[((t/\alpha)\log(1/\lambda_1))^{1/p} - 1\right].$$

Therefore,

$$\left(\frac{t}{\alpha}\right)^{1/p} \leq \lim_{N\to\infty} \sum_{k=0}^{\infty} \frac{h_1}{(1 + \lambda_1 \beta \exp(\alpha x^p))^l} \leq \left(\frac{1}{\alpha}\right)^{1/p},$$

and letting $t \to 1$, we achieve the desired result.

The eigen-functions of the Gaussian RKHS are (Zhu et al., 1998):

$$\phi_k(x) = (\sqrt{5}/4)^{1/4}(2^{k-1}k!)^{-1/2}e^{-(\sqrt{5}-1)x^2/4}H_k\big((\sqrt{5}/2)^{1/2}x\big),$$

where $H_k(\cdot)$ is the $k$-th Hermite polynomial, for $k = 0, 1, 2, \cdots$. The eigen-values are $\mu_k = \ell^{2k+1}$, where $\ell = (\sqrt{5}-1)/2$, can also be re-written as an exponential form: $\mu_k = \ell\exp(2\ln\ell \cdot k)$. According to Krasikov (2004), we can get that

$$c_K = \sup_{k \in \mathbb{N}} \|\phi_k\|_{\sup} \leq \frac{2e^{15/8}(\sqrt{5}/4)^{1/4}}{3\sqrt{2\pi}2^{1/6}} \leq 1.336.$$

Therefore, we have

$$\sigma_g^2 \leq \lim_{N \to \infty} \sigma^2 c_K^2 h_1 \sum_{k=0}^{\infty}(1+\lambda\ell\exp(-2\ln\ell\cdot k))^{-2} \leq c_K^2 \cdot 2\sigma^2 \log(1/\ell) \leq 1.7178\sigma^2,$$

The values of $\sigma_\beta^2, \rho_T^2$ and $\sigma_T^2$ can be obtained by the rescaling as well. □